# Arbitrary Lagrangian-Eulerian Discontinuous Galerkin Methods

with Application to Compressible Flows

*Author*
Jayesh Vinay Badwaik

# Contents

# Chapter 1 Introduction

In the past few decades, mesh based methods like finite volume and discontinuous Galerkin finite element methods have been exceedingly successful in the task of modelling compressible fluid flows. Due to the upwind nature of the finite volume schemes, dissipation is implicitly built in the schemes. This allows the schemes to be able to compute discontinuous solutions in a stable manner. Discontinuous Galerkin methods can be considered to be higher order generalizations of the finite volume methods which also make use of the upwinding technology but instead of cell averages, evolve a polynomial inside each cell.

In the most basic form of these methods, the mesh used for simulations is kept fixed over the course of simulations. This can sometimes be an issue for unresolved solutions, where the truncation error from the simulation can dominate the behavior of fluids. A striking example of such a situation is the growth of the Kelvin-Helmholtz instabilities between two fluids of different densities. The problem setup consists three scenarios each having an unresolved discontinuity (a slip-surface) between two fluids moving in opposite directions, but in each scenario, there is a boost Fig. 1.0.1.

Due to the principle of Galilean invariance, the fluid should behave identically irrespective of the value of $V_0$. However, using AREPO code [31] for example, one can see in Fig. 1.0.2 that the instabilities disappear in the moving frames. The disappearance of the instabilities is due to dampening of the unstable modes in the fluid by the numerical dissipation of method.

The problem of numerical dissipation is not limited to the moving frames. The effect is observable even in rest frames where the instabilities disappear at lower resolutions. So, even though the methods themselves are Galilean invariant, the corresponding truncation error is not Galilean invariant. Apart from suppression of instabilities, the violation of Galilean invariance of truncation error is a major contributor to smearing out of shocks and contact discontinuities, especially in higher-order methods. Even for smooth solutions, the accuracy of the solution can be reduced due to the effects of truncation error.

One of the ways to deal with these challenges is to move the mesh along with the fluid, thereby reducing the truncation error in the simulations. A natural choice is to move the mesh with the velocity of the fluid and thereby exploit the Lagrangian formulation for fluid mechanics. The resulting method is called the Lagrangian method and was first developed in 1950s and 60s [33, 34]. A general review of these methods can be found in [5].

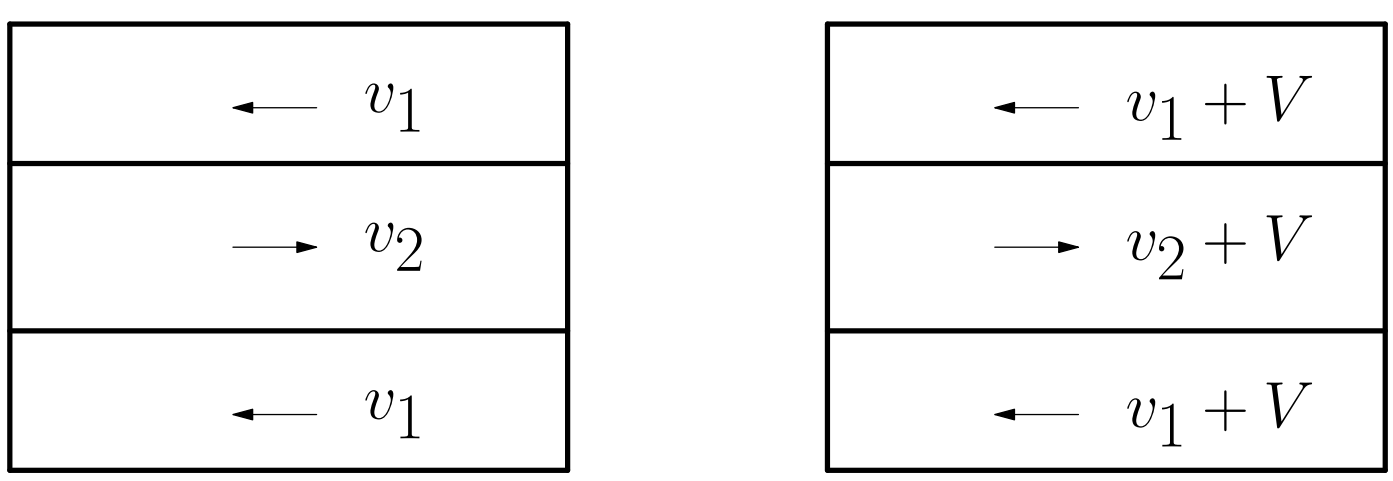

Figure 1.0.1: Problem Setup to Study Kelvin-Helmholtz Instabilities

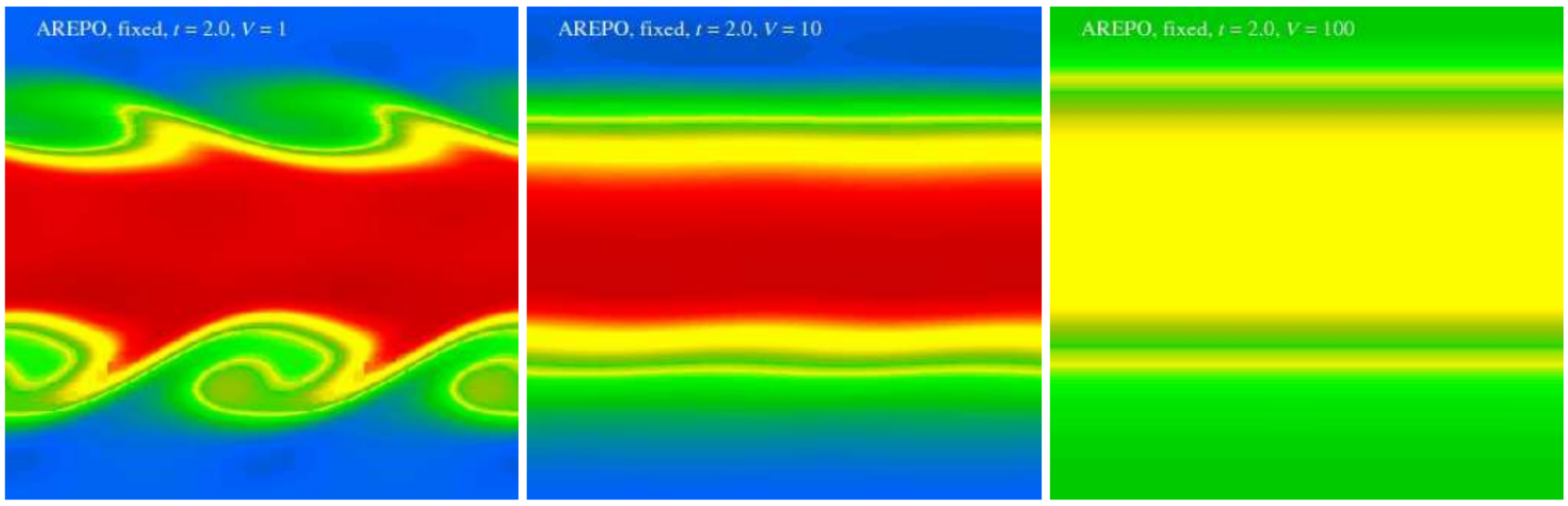

Figure 1.0.2: Kelvin-Helmholtz Instabilities in Fluid Flow

In Lagrangian simulations, a given region is always approximated by the same number of mesh points. This ensures that the accuracy of the simulation is in general maintained to the accuracy of the initial approximation. Furthermore, compared to the fixed mesh methods, the number of points required to get a similar accuracy in Lagrangian simulations is often surprisingly small. This makes Lagrangian methods quite attractive for fluid simulations.

Unfortunately, the very features of the Lagrangian method which make it so effective are also the ones which make it totally unsatisfactory for simulations involving large shear flows. The property of mesh points to exactly follow the fluid flow causes the underlying mesh to become highly distorted, if not degenerate. This results in the reduction of accuracy of approximation and consequently of the simulation, with the simulation not being able to continue till completion in some cases.

In case of small deformations, one can add corrections to the mesh velocity in order to maintain the quality of the mesh for a longer time and often till the completion of the simulation. At this stage however, the mesh velocity is not equal the fluid velocity, and Lagrangian formulation is no longer an appropriate setting for our computations. We must instead use the Arbitrary Lagrangian-Eulerian (ALE) formulation where the vertices of the computational mesh may be moved in an arbitrary fashion. Methods using ALE formulation were originally developed in [29, 18, 32, 23] in the context of finite difference and finite volume methods. The method was subsequently adopted in the finite element context and early applications can found in [17, 25, 4, 24].

In theory, the mesh velocities in the ALE formulation can be any arbitrary value but in practice, it is desirable to keep the mesh velocities as close to fluid velocity as possible. Methods which satisfy this property are called **almost-Lagrangian** methods. The methods are practical in cases of small deformations of the meshes, but for large deformations, it becomes impractical to maintain a good quality mesh with mesh velocities close to fluid velocities and it becomes necessary to make changes to the mesh topology.

An approach to handle the problem of the mesh quality is to move the mesh for as long as possible with the fluid velocity, and remesh when the mesh quality becomes unacceptable [3, 22]. For each remeshing, the simulation has to be stopped, a new mesh must be generated and the solution must be transferred to the new mesh. This approach can be efficient if the number of remeshing required are very low, since the simulation is fully ALE and free from interpolation errors.

The approach described above is sometimes also known as remap ALE methods or indirect ALE methods. Remap ALE schemes have been used extensively in some problems and some of the recent work on the topic can be found in [6, 8, 9, 27, 36]. For large deformations, the number of required remeshing increase and this can be both costly and result in poor accuracy due to the solution transfer step.

Another approach to get around this problem was presented in [31], where the connectivity of the moving mesh is dynamically regenerated via a moving unstructured but conforming Voronoi tessellation of the domain. Building on the idea, a higher order ALE Discontinuous Galerkin (DG) method was designed by [19]. In these methods, a new mesh is generated at every time step, a correspondence is generated between the elements in the two time steps, and finally a ALE DG method is solved on the mesh, where the space-time information is gleaned from the correspondence information.

In the above method, there is a possibility of non-regular ("crazy") finite elements appearing in the mesh which is taken care of using the techniques in [19]. Due to no inherent need to remap the solutions in these methods, these methods are often referred to as *direct* ALE DG methods. Since there is no remapping of solutions needed, the methods do not suffer from the interpolation error, however, it requires a mesh to be generated at every time step and corresponding connection information to be determined which can be costly.

There is a third approach which was first outlined in [13, 16]. The approach was originally designed for mesh motion due to a moving interface or boundary, with applications in particular to fluid-structure interactions. The moving interface introduces deformation of the domain, which degrades the quality of the mesh. The focus in these papers is to preserve the adapt the mesh to the deforming domain by performing only local remeshing operations, such as vertex insertion, vertex collapse, connectivity changes and vertex displacements. The advantage of this method is that it maintains an acceptable mesh quality without needing to stop, remesh and resume the simulation. Another important aspect is that for all the cells which do not need to be modified, there is no interpolation error during the local cell remeshing. This helps in maintaining a high accuracy of solutions through the simulation. [1, 2] further showed that for ALE DG simulations with a elasticity-based mesh velocities, only mesh velocity corrections and connectivity changes (face swapping) is enough to maintain a good quality mesh. Finally, the local techniques ensure that there's no global synchronization bottleneck in the simulations, making such techniques suited large scale simulations.

In this report, we attempt to adopt the local remeshing techniques from [13, 16, 1, 2] for problems in compressible flows itself. In particular, we want to use the local adaptation techniques to preserve the quality of the mesh while carrying out an almost-Lagrangian motion of the mesh. The challenge in such an effort comes from the fact that the fluid motion is often much more irregular than moving interfaces. In particular, the fluid can demonstrate shear and vortical flow at the scale of single mesh cells and makes the choice of mesh velocity during simulation critical.

In this report, we propose a collection of methods to make such an approach possible for Euler equations in one and two dimensions. We propose an explicit single-step ALE DG scheme for hyperbolic conservation laws. The scheme considerably reduces the numerical dissipations introduced by the Riemann solvers. We show that the scheme also preserves the constant states for any mesh motion. We then study the effect of mesh quality on the accuracy of the simulations, and based on that, come up with a mesh quality indicator for the ALE DG method. Based on the considerations from the study on mesh quality, we design a local mesh velocity algorithm to compute the motion of the mesh. And finally, we propose a local mesh adaptation algorithm to control the quality of the mesh, and prevent the mesh from degradation.

## 1.1 Organization of the Part

**Chapter 2** We start by introducing the basic notations to describe hyperbolic systems of conservation laws and we discuss the weak formulation for hyperbolic conservation laws. We introduce some of the ideas in numerical method, and finally we introduce conservation laws in the ALE setting.

**Chapter 3** In this chapter, we introduce the ALE DG scheme. We define the mesh on which the scheme is to be defined, next, we describe the mesh evolution step in a time step, and the ALE DG Scheme in that time step. We prove that our scheme is conservative and preserves constant states under arbitrary mesh motion. We then describe the algorithm to compute vertex velocities in the mesh. We also study the dissipation of the ALE DG scheme, and study the accuracy of the solution with respect to the mesh quality. Finally, we describe the TVB limiter used in the simulations.

**Chapter 4** In this chapter, we describe different techniques to preserve the quality of the mesh in the simulation. We describe mesh velocity algorithms, following which, we describe some edge swapping algorithms. We also describe our method of transformation of solution on arbitrary meshes. And finally, we describe the complete ALE DG algorithm.

**Chapter 5** In this chapter, we use the method to study on one dimensional problems of Euler equations.

**Chapter 6** In this chapter, we use the method to study on two dimensional problems of Euler equations.

Chapter

# 2 Hyperbolic Conservation Laws

In this chapter, we introduce hyperbolic conservation laws in Section 2.1, the theory of weak formulation in Section 2.2 and a brief introduction to the theory of numerics in Section 2.3. Finally, for the purpose of the ALE DG method, we introduce the ALE formulation and state the hyperbolic conservation laws in the ALE form in Section 2.4.

## 2.1 Hyperbolic Conservation Laws

Consider an open bounded set $\Omega \subset \mathbb{R}^d$ with a Lipschitz boundary. Let $\boldsymbol{\nu}$ denote the outward unit normal along $\Omega$. Let $\mathbf{u}(\mathbf{x}, t)$, $\mathbf{u} \colon \mathbb{R}^d \times \mathbb{R}^+ \to \mathbb{R}^n$, be a function describing the density of various conserved quantities, let $\mathcal{F}(\mathbf{x}, t, \mathbf{u})$, $\mathcal{F} : \mathbb{R}^d \times \mathbb{R}^+ \times \mathbb{R}^n \to \mathbb{R}^{d \times n}$, $\mathcal{F} = (\mathbf{F}_1, \mathbf{F}_2, \dots, \mathbf{F}_d)^{\mathsf{T}}$, govern the rate of change of $\mathbf{u}(\mathbf{x}, t)$ within $\Omega$. One can then assume that the behavior of the quantity $\mathbf{u}$ in $\Omega$ can be described by the the integral form of the conservation laws

$$\frac{\mathrm{d}}{\mathrm{d}t} \int_{\Omega} \mathbf{u}(\mathbf{x}, t) \, \mathrm{d}\mathbf{x} + \int_{\partial\Omega} \mathcal{F}(\mathbf{x}, t, \mathbf{u}) \cdot \boldsymbol{\nu} \, \mathrm{d}s = 0 \tag{2.1.1}$$

**Equation 2.1.1:** Integral Form of Conservation Law

We can apply the Gauss divergence theorem to the second integral term to get

$$\frac{\mathrm{d}}{\mathrm{d}t} \int_{\Omega} \mathbf{u}(\mathbf{x}, t) \, \mathrm{d}\mathbf{x} + \int_{\Omega} \nabla_{\mathbf{x}} \cdot \mathcal{F}(\mathbf{x}, t, \mathbf{u}) \, \mathrm{d}\mathbf{x} = 0 \tag{2.1.2}$$

**Equation 2.1.2:** Integral Form of Conservation Law : Divergence Form

Next, by noting that the domain is independent of time, we take the time derivative inside the first integral and then substitute it with a partial derivative instead and finally, since the choice of initial domain $\Omega$ was arbitrary, we obtain the conservation law in a partial differential equation (PDE) form Eq. (2.1.3).

$$\frac{\partial}{\partial t}\mathbf{u}(\mathbf{x},t) + \nabla \cdot \mathcal{F}(\mathbf{x},t,\mathbf{u}) = 0 \tag{2.1.3}$$

**Equation 2.1.3:** PDE Form of Conservation Law

**Definition 2.1.1** (Hyperbolic Conservation Law)**.** *Let $D_{\mathbf{u}}\mathcal{F}$ be the Jacobian tensor of $\mathcal{F}$ with respect to $\mathbf{u}$. Then the system of conservation laws (Eq. (2.1.3)) is said to be* ***hyperbolic*** *if for all $\mathbf{u} \in \mathbb{R}^n$ and all $\boldsymbol{\nu} \in \mathbb{R}^d$, the matrix $\langle D_{\mathbf{u}}\mathcal{F}, \boldsymbol{\nu}\rangle$ has $n$ real eigenvalues and $n$ linearly independent eigenvectors. If, in addition, these eigenvalues are distinct, then the system is said to be* ***strictly hyperbolic****.*

Primary uses of conservation laws are in situations where the value of the conserved quantity is known at an initial time and we want to compute the behavior at some future time, these problems are called initial value problems. Sometimes, apart from the initial, one is also provided with the value at the boundary of the domain, in which case, the problem is then called an initial-boundary value problem. For most of the theoretical analysis in this part, we will restrict ourselves to initial value problem.

**Definition 2.1.2** (Initial Value Problem)**.** *The initial value problem for the conservation law Eq. (2.1.3) is a Cauchy problem given as*

$$\frac{\partial \mathbf{u}(\mathbf{x},t)}{\partial t} + \nabla_{\mathbf{x}} \cdot \mathcal{F}(\mathbf{x},t,\mathbf{u}) = 0 \qquad (\mathbf{x},t) \text{ in } \mathbb{R} \times [0,T] \tag{2.1.4a}$$

$$\mathbf{u}(\mathbf{x},0) = \mathbf{u}_0(\mathbf{x}) \tag{2.1.4b}$$

**Definition 2.1.3** (Initial-Boundary Value Problem)**.** *Assume that $\Omega \subset \mathbb{R}^d$ is a domain with a Lipschitz boundary $\Gamma \in \mathbb{R}^{d-1}$. The initial-boundary value problem is a Cauchy problem for given as*

$$\frac{\partial \mathbf{u}(\mathbf{x},t)}{\partial t} + \nabla_{\mathbf{x}} \cdot \mathcal{F}(\mathbf{x},t,\mathbf{u}) = 0 \qquad (\mathbf{x},t) \text{ in } \Omega \times [0,T] \tag{2.1.5a}$$

$$\mathbf{u}(\mathbf{x},0) = \mathbf{u}_0(\mathbf{x}) \qquad \mathbf{u}(\sigma,0) = \mathbf{u}_\Gamma(\sigma) \tag{2.1.5b}$$

As one can see, it is necessary for the function $\mathbf{u}(\mathbf{x},t)$ to be sufficiently smooth in order to meaningfully talk about derivatives and hence solutions to Eq. (2.1.4) and Eq. (2.1.5). Such a solution is known as a **classical solution** or a **strong solution**. However, in practice, it found that this condition is too restrictive. For example, even for smooth initial and boundary conditions , the solutions can develop discontinuities in finite time [14]. The equation then ceases to have solutions in the classical sense. In order to obtain a more useful theory, it is thus necessary to look at a more permissive frameworks.

## 2.2 Weak Formulation

Weak formulation is an important tool for analysis of mathematical equations since it allows porting of techniques from linear algebra to functional analysis. The formulation has been successful in allowing us to deal with discontinuous functions as solutions of conservation laws. The solutions of a partial differential equations which satisfy the weak formulation are known as the weak solutions.

However, unlike the classical solutions, weak solutions are not unique, and we need additional admissibility conditions to single out the physically relevant solutions. We use the framework of entropy solutions introduced by Kruzkov and later built on by others for this purpose. In this section, we introduce the concept of weak solutions and entropy solutions to the hyperbolic conservation law and discuss the well-posedness of the hyperbolic conservation laws in context of entropy solutions.

**Definition 2.2.1** (Weak Solution). *A measurable function $\mathbf{u}(\mathbf{x}, t)$ is said to be a* ***weak solution*** *to Eq. (2.1.3), if it satisfies the equation in the distributional sense. In the weak form, it implies that for all $\phi \in C_c^\infty(\mathbb{R}^d \times \mathbb{R}^+; \mathbb{R})$, $\mathbf{u}(\mathbf{x}, t)$ satisfies the weak formulation of the conservation law given by*

$$\int\limits_{\mathbb{R}^d \times \mathbb{R}} \left[ \mathbf{u}(\mathbf{x}, t) \frac{\partial \phi}{\partial \mathbf{t}} + \mathscr{F}(\mathbf{x}, t, \mathbf{u}) \nabla_{\mathbf{x}} \phi(\mathbf{x}, t) \right] \mathrm{d}\mathbf{x} \, \mathrm{d}t + \int\limits_{\mathbb{R}^d} \mathbf{u}(\mathbf{x}, 0) \phi(\mathbf{x}, 0) \, \mathrm{d}\mathbf{x} = 0 \qquad (2.2.1)$$

We derive Eq. (2.2.1) by computing the inner product of Eq. (2.1.3) with $\phi$ over space and time and then applying integration by parts to transfer the derivatives onto $\phi$. The advantage of considering the weak solutions is that the range of admissible solutions is now much larger, since the solutions are now allowed to be discontinuous. However, the discontinuities of such weak solutions are not arbitrary and instead must satisfy the Rankine–Hugoniot conditions.

**Definition 2.2.2** (Rankine–Hugoniot Condition). *Let $\Gamma(\mathbf{x}, t)$ be the surface of a discontinuity in the $\mathbf{x} - t$ plane of a weak solution $\mathbf{u}$ to Eq. (2.1.3). Let $\boldsymbol{\nu}(\mathbf{x}, t)$ be the unit normal vector to $\Gamma$. Let $\mathbf{u}_+(\mathbf{x}, t)$ and $\mathbf{u}_i(\mathbf{x}, t)$ denote the traces of $\mathbf{u}$ at $\Gamma^+$ and $\Gamma^-$ respectively. Then, the* ***Rankine–Hugoniot condition*** *is given as*

$$(\mathbf{u}^+ - \mathbf{u}^-) \frac{\partial \boldsymbol{\nu}}{\partial t} + (\mathbf{u}^+ - \mathbf{u}^-) \cdot \boldsymbol{\nu} = 0 \qquad (2.2.2)$$

As one can see, a classical solution satisfies a Rankine–Hugoniot condition and hence is also a weak solution. However, as mentioned before, the weak solutions are not unique, and we consider the framework of entropy solutions to single out the physically relevant solution.

**Definition 2.2.3** (Entropy - Entropy Flux Pair). *Given the conservation law Eq. (2.1.3), assume that there exists smooth functions $\boldsymbol{\eta}(\mathbf{x}, t, \mathbf{u})$ and $\mathscr{Q}(\mathbf{x}, t, \mathbf{u}) = (\mathbf{Q}_1, \mathbf{Q}_2, \dots, \mathbf{Q}_d)$ such that*

$$D_{\mathbf{u}} \mathbf{Q}_i = (D_{\mathbf{u}} \boldsymbol{\eta})^{\mathsf{T}} D_{\mathbf{u}} \mathbf{F}_i \qquad \forall \mathbf{u} \in \mathbb{R}^n \text{ and } i = 1, 2, \dots, d \qquad (2.2.3)$$

*Then the conservation law Eq. (2.1.3) is said to be endowed with* ***entropy*** *$\boldsymbol{\eta}$ and* ***entropy flux*** *$\mathscr{Q}$. The pair of $\boldsymbol{\eta}$ and $\mathscr{Q}$ together is called* ***entropy-entropy flux pair****. In case of scalar hyperbolic conservation laws, the regularity conditions on $\boldsymbol{\eta}$ (and thereby on $\mathscr{Q}$) can be relaxed to allow $\boldsymbol{\eta}$ (and thereby $\mathscr{Q}$) to only be locally Lipschitz continuous. In that case, Eq. (2.2.3) is written as*

$$\mathbf{Q}_i = \int\limits_{-\infty}^{\mathbf{u}} (D_{\mathbf{u}} \boldsymbol{\eta})^{\mathsf{T}} D_{\mathbf{u}} \mathbf{F}_i \, \mathrm{d}\mathbf{u} \qquad \forall \mathbf{u} \in \mathbb{R}^n \text{ and } i = 1, 2, \dots, d \qquad (2.2.4)$$

**Definition 2.2.4** (Entropy Solution). *A weak solution, $\mathbf{u}(\mathbf{x}, t)$, to the initial value problem Eq. (2.1.4) is said to be an* ***entropy solution*** *if it satisfies, in the distributional sense, the entropy* *admissibility condition*

$$\frac{\partial \boldsymbol{\eta}(\mathbf{x},t)}{\partial t} + \nabla_{\mathbf{x}} \cdot \boldsymbol{\mathcal{Q}}(\mathbf{x},t,\mathbf{u}) \leq 0 \tag{2.2.5}$$

*This then implies that the map* $t \mapsto \boldsymbol{\eta}(\mathbf{u}(\cdot,t))$ *is continuous on* $[0,T) \setminus \mathscr{F}$ *in* $L^{\infty}(\mathbb{R}^d)$ *weak*, where* $\mathscr{F}$ *is at most countable. Furthermore, for every non-negative Lipschitz test function* $\phi(x,t)$*, with compact support in* $\mathbb{R}^d \times [0,T)$ *and all* $\tau \in [0,T) \setminus \mathscr{F}$*, we have*

$$\int\limits_{\mathbb{R}^d \times \mathbb{R}} \left[\boldsymbol{\eta}(\mathbf{u}(\mathbf{x},t))\frac{\partial \phi}{\partial \mathbf{t}} + \boldsymbol{\mathcal{Q}}(\mathbf{x},t,\mathbf{u})\, \nabla_{\mathbf{x}}\, \phi(\mathbf{x},t)\ \mathrm{d}\mathbf{x}\,\mathrm{d}t\right] + \int\limits_{\mathbb{R}^d} \mathbf{u}(\mathbf{x},\tau)\phi(\mathbf{x},\tau)\ \mathrm{d}\mathbf{x} \geq 0 \tag{2.2.6}$$

In general, it is important to determine conditions under which the set $\mathscr{F}$ is actually empty, but this is difficult problem to solve. Even so, it is generally a great help to the analysis if at least $0 \notin \mathscr{F}$, in which case Eq. (2.2.6), with $\tau = 0$, becomes

$$\int\limits_{\mathbb{R}^d \times \mathbb{R}} \left[\boldsymbol{\eta}(\mathbf{u}(\mathbf{x},t))\frac{\partial \phi}{\partial \mathbf{t}} + \boldsymbol{\mathcal{Q}}(\mathbf{x},t,\mathbf{u})\, \nabla_{\mathbf{x}}\, \phi(\mathbf{x},t)\ \mathrm{d}\mathbf{x}\,\mathrm{d}t\right] + \int\limits_{\mathbb{R}^d} \mathbf{u}(\mathbf{x},0)\phi(\mathbf{x},0)\ \mathrm{d}\mathbf{x} \geq 0 \tag{2.2.7}$$

Accordingly, it is Eq. (2.2.7) rather than the slightly weaker condition Eq. (2.2.6) that is often postulated in the literature as the entropy admissibility criterion for the weak solution $\mathbf{u}$. We note that the admissible weak solutions characterized through Eq. (2.2.7) do not necessarily possess the semigroup property, i.e. $\mathbf{u}(x,t)$ admissible does not generally imply that $\mathbf{u}(x,t+\tau) = \mathbf{u}_{\tau}(x,t)$ is also admissible, for all $\tau \in [0,T)$.

Given an entropy–entropy flux pair, often, the corresponding entropy admissibility condition eliminates some, but not necessarily all, of the spurious weak solutions. One of the ways to overcome this deficiency is to require Eq. (2.2.5) to hold simultaneously for every convex entropy of the system. However, we note that this is viable only for those systems that are endowed with a rich family of entropies. This is particularly possible in the case of scalar hyperbolic conservation laws, wherein we have the family of **Kruzkov entropies**.

**Definition 2.2.5** (Kruzkov Entropy–Entropy Flux Pair). *Consider a scalar hyperbolic conservation law (*$\mathbf{u} \in \mathbb{R}$*) of the form of* Eq. (2.1.3)*. In Kruzkov entropy formulation, for a given constant* $\mathbf{c} \in \mathbb{R}$*, the corresponding entropy–entropy flux pair is given by*

$$\boldsymbol{\eta}_c(\mathbf{u}) = |\mathbf{u} - \mathbf{c}| \qquad \boldsymbol{\mathcal{Q}}_c(\mathbf{u}) = \mathrm{sgn}(\mathbf{u}-\mathbf{c})\Big(\boldsymbol{\mathcal{F}}(\mathbf{x},t,\mathbf{u}) - \boldsymbol{\mathcal{F}}(\mathbf{x},t,\mathbf{c})\Big) \tag{2.2.8}$$

**Definition 2.2.6** (Kruzkov Entropy Solution). *A weak solution* $\mathbf{u}(\mathbf{x},t)$ *to the initial value problem* Eq. (2.1.4) *is known as the Kruzkov entropy solution if* Eq. (2.2.7) *holds for every Kruzkov entropy-entropy flux pair and all non-negative Lipschitz continuous test functions* $\phi$ *on* $\mathbb{R}^d \times [0,\infty)$ *with compact support.*

**Theorem 2.2.1** (Well-Posedness of the Initial Value Problem for the Scalar Case). *Consider the initial value problem* Eq. (2.1.4) *for a scalar hyperbolic conservation law* $(\mathbf{u} \in \mathbb{R})$ *where* $\boldsymbol{\mathcal{F}}$ *is locally Lipschitz continuous and* $\mathbf{u}_0 \in L^{\infty}(\mathbb{R}^d)$*. Then, there exists a unique Kruzkov entropy solution to* Eq. (2.1.4) *and*

$$\mathbf{u}(\cdot,t) \in C\left([0,\infty); L^1_{loc}(\mathbb{R}^d)\right) \tag{2.2.9}$$

The result completes the well-posedness theory for the scalar case. In contrast, the situation for systems of hyperbolic conservation laws is quite different. Apart from some partial results for one-dimensional systems in [7, 10, 20], no well-posedness results are available for general multi-dimensional systems. In fact, in recent years, [26] have shown that the entropy solutions to some specific multi-dimensional systems are not unique.

This suggests that the notion of entropy solutions may not be adequate to establish the existence and uniqueness of solutions for a general system, and there might be a need for a different solution framework. The search for a correct framework is still an active area of research and we will not delve further into it, apart from noting the popularity of the framework of measure valued solutions introduced in [15] in the current efforts.

As noted above, the entropy conditions are not enough to establish well-posedness of systems of conservation laws. Nevertheless, they play an important role in analysis of systems of conservation laws by providing various useful global estimates.

**Theorem 2.2.2** (Global Stability Estimate). *Let $\mathbf{u}(\mathbf{x},t)$ be a weak solution of Eq. (2.1.4) on $\mathbb{R}^d \times [0,T)$ with $\mathbf{u}(\cdot,0) - \overline{\mathbf{u}} \in L^2(\mathbb{R}^d)$ for some constant $\overline{\mathbf{u}}$. Assume that $\mathbf{u}$ satisfies the entropy admissibility condition Eq. (2.2.5) relative to a uniformly convex entropy $\boldsymbol{\eta}$, normalized so that $\boldsymbol{\eta}(\overline{\mathbf{u}}) = 0$, $D\boldsymbol{\eta}(\overline{\mathbf{u}}) = 0$. Then $\mathbf{u}$ satisfies Eq. (2.2.7) if and only if*

$$\int_{\mathbb{R}^d} \boldsymbol{\eta}(\mathbf{u}(\mathbf{x},\tau))\,\mathrm{d}\mathbf{x} \leq \int_{\mathbb{R}^d} \boldsymbol{\eta}(\mathbf{u}(\mathbf{x},0))\,\mathrm{d}\mathbf{x} \qquad \forall 0 < \tau < T \tag{2.2.10}$$

*It thus follows that Eq. (2.2.6) holds for all $\tau \in [0,T)$ if and only if the integral $\int_{\mathbb{R}^d} \boldsymbol{\eta}(\mathbf{u}(\mathbf{x},t))\,\mathrm{d}x$ is non-increasing on $[0,T)$.*

**Theorem 2.2.3** (Stability of Classical Solutions). *Let $\overline{\mathcal{B}}_\rho$ be a closed ball in $\mathbb{R}^n$ centered at origin. Assume that the system of conservation laws in Eq. (2.1.3) is endowed with an entropy–entropy flux pair $(\boldsymbol{\eta}, \boldsymbol{\mathcal{Q}})$, where $D^2\boldsymbol{\eta}(\mathbf{u})$ is positive definite on $\overline{\mathcal{B}}_\rho$. Suppose that $\overline{\mathbf{u}}(\mathbf{x},t)$ is a classical solution of the initial value problem Eq. (2.1.4) on $[0,T)$, taking values in $\overline{\mathcal{B}}_\rho$. Let $\mathbf{u}$ be any weak solution of the initial value problem Eq. (2.1.4) on $[0,T)$, taking values in $\overline{\mathcal{B}}_\rho$, which satisfies the entropy admissibility condition Eq. (2.2.7). Then*

$$\int_{|\mathbf{x}|<r} |\mathbf{u}(\mathbf{x},t) - \overline{\mathbf{u}}(\mathbf{x},t)|^2\,\mathrm{d}\mathbf{x} \leq ae^{bt} \int_{|\mathbf{x}|<r+st} |\mathbf{u}(\mathbf{x},0) - \overline{\mathbf{u}}(\mathbf{x},0)|^2\,\mathrm{d}\mathbf{x} \tag{2.2.11}$$

*holds for any $r > 0$ and $t \in [0,T)$, with positive constants $s, a$, depending solely on bounds on $\boldsymbol{\mathcal{F}}, \boldsymbol{\eta}, \boldsymbol{\mathcal{Q}}$ and their derivatives on $\overline{\mathcal{B}}_\rho$, and $b$ that also depends on the Lipschitz constant of $\overline{\mathbf{u}}$. In particular, $\overline{\mathbf{u}}$ is the unique admissible weak solution of the initial value problem Eq. (2.1.4).*

## 2.3 Numerical Methods in Hyperbolic Conservation Laws

As we have seen in Section 2.2, discontinuous solutions of conservation laws are not unique in general and we routinely need additional admissibility conditions to single out physically relevant solutions. In numerics, the concerns translate to being able to develop numerical schemes which provide a physically relevant accurate approximation to the solutions.

The computation is straightforward in case of smooth solutions with the method of finite differences being good enough to compute physically correct approximations, even if they are not very accurate. The difficulty appears near the discontinuities, where the solution is supposed to satisfy the PDE only in the weak form. As such, the classical finite difference methods are supposed to be inappropriate near the discontinuities, as is indeed evinced by their poor performance in such scenarios.

Hence, there is a strong motivation to understand the structure of the solution for discontinuous data. Understanding the structure requires a good understanding of the theory as well as some physical intuition about the behavior of solutions. For linear hyperbolic systems, the **characteristics** play a major role in such understanding. For non-linear problems, the generalization which is most often used is the solution of a **Riemann problem**.

Solutions of Riemann problem have a relatively simple structure, and can often be calculated explicitly. In numerical methods, we typically encounter a situation, where we need to understand the local structure of solution around a point $\mathbf{x}_j$ in the grid, while being provided with data at the point $\mathbf{x}_j$ and neighboring points ($\mathbf{x}_{j-1}$ for example. A great deal of such information is obtained by solving a Riemann problem at point $\mathbf{x}_j$ and many numerical methods make use of these Riemann solutions for their design.

For scalar equations, one can prove often prove convergence results of the numerical methods and thus be assured about the correctness of the methods. Such results often rely on non-linear stability estimates for the numerical methods. In particular, for many methods, it is possible to show that the total variation of the solution is non-increasing with time. This is often enough to obtain some convergence results and also guarantees that the method does not give rise to any spurious oscillations. Such methods are known as Total Variation Diminishing (TVD) methods.

The non-linear stability estimates for these methods however are often not enough to prove convergence in case of systems of conservation laws, and hence, cannot provide similar guarantees as in the scalar case. Nevertheless, using this condition as a design requirement for new methods has been found to be extremely useful in practice and has led to the design of many good numerical methods. We end the section with a brief description of the **Riemann problem** for scalar hyperbolic conservation laws in one dimension.

### 2.3.1 Riemann Problem in One Dimensional for the Scalar Case

**Definition 2.3.1** (Riemann Problem in One Dimension). *Let $\mathbf{u} \in \mathbb{R}^n$ and let $\mathbf{F}(x,t,\mathbf{u}) \colon \mathbb{R} \times \mathbb{R}^+ \times \mathbb{R}^n \to \mathbb{R}^n$. We assume that flux Jacobian $\partial \mathbf{F}/\partial \mathbf{u}$ has real and distinct eigenvalues. The* ***Riemann problem*** *in one dimension is then defined as the initial value problem illustrated in Fig. 2.3.1*

$$\frac{\partial \mathbf{u}}{\partial x} + \frac{\partial \mathbf{F}(x,t,\mathbf{u})}{\partial t} = 0 \qquad (x,t) \in (\mathbb{R}, \mathbb{R}^+) \tag{2.3.1a}$$

$$\mathbf{u}(x,0) = \begin{cases} \mathbf{u}_l & \text{if } x < 0 \\ \mathbf{u}_r & \text{if } x > 0 \end{cases} \tag{2.3.1b}$$

We note that the initial condition is discontinuous if $\mathbf{u}_l \neq \mathbf{u}_r$. As shown in Fig. 2.3.2, the solution to the Riemann problem consists of $n$ waves emanating from the origin, one corresponding to each eigenvalue. Furthermore, [21] showed that the solutions of Eq. (2.3.1a) are self-similar, that is, the solutions are of the form $\mathbf{u}(x,t) = \mathbf{w}(x/t)$ and consist of $n+1$ constant states separated by $n$ waves. The $n+1$ states are connected by the following types of waves.

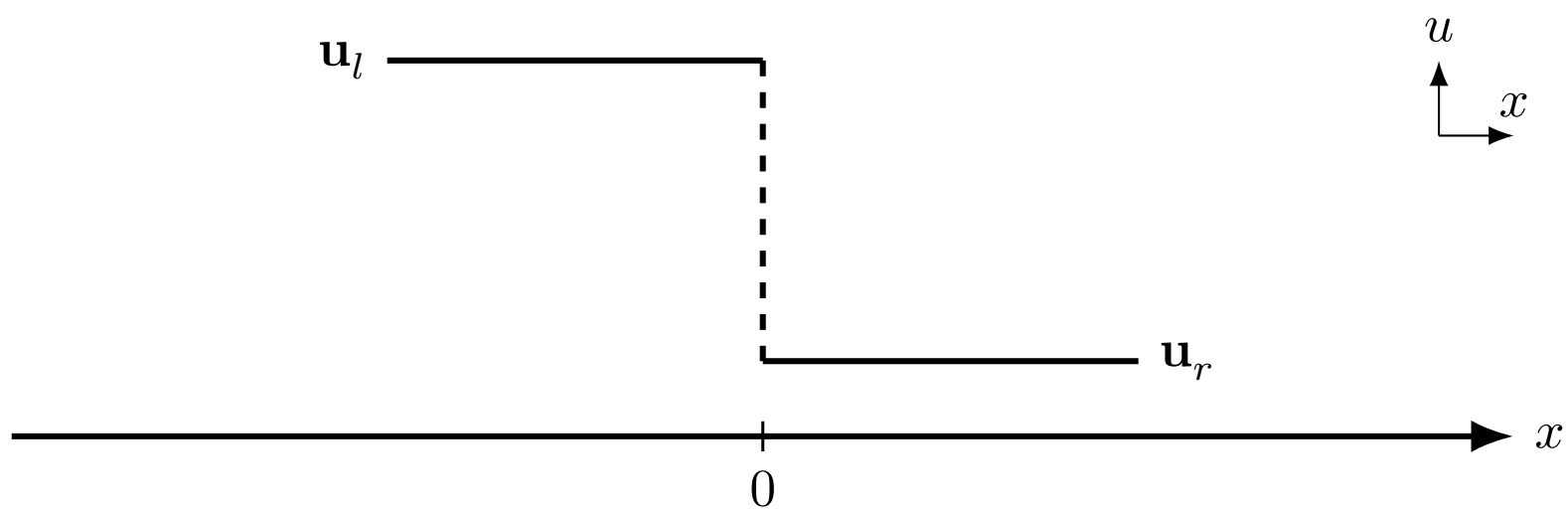

Figure 2.3.1: Illustration of the Initial Data for the Riemann Problem

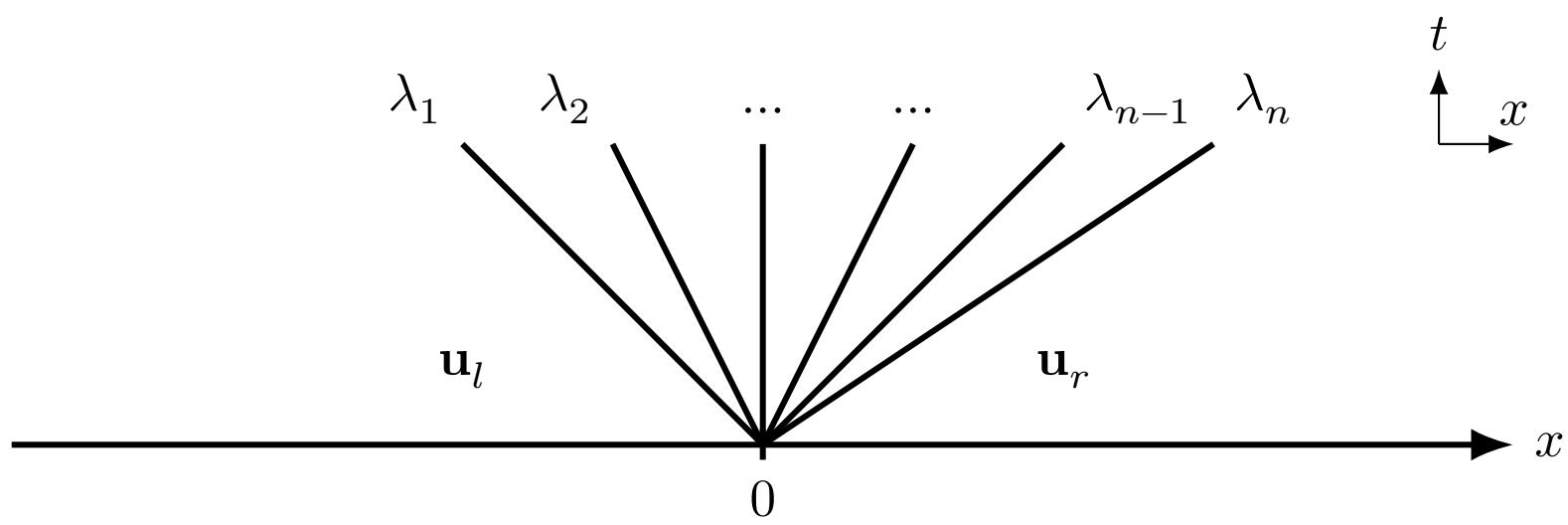

Figure 2.3.2: Structure of Solution of a Riemann Problem

1. **Shock Waves** The $\lambda_i$-wave is a shock wave, if it corresponds to a genuinely non-linear field and connects two states $\mathbf{u}_l$ and $\mathbf{u}_r$ through a single jump discontinuity. The discontinuity moves with a speed $s_i$ given by the Rankine–Hugoniot condition Eq. (2.2.2). Furthermore, the Lax entropy condition holds, i.e.

$$\lambda_i(\mathbf{u}_l) \geq s_i \geq \lambda_i(\mathbf{u}_r) \tag{2.3.2}$$

which can be deduced from the entropy condition Eq. (2.2.5) for a convex flux. This implies that the **characterics** $\mathrm{d}x/\,\mathrm{d}t = \lambda_i$ on both sides of the shock line $\mathrm{d}x/\,\mathrm{d}t = \lambda_i$ run into the shock wave as shown in Fig. 2.3.3.

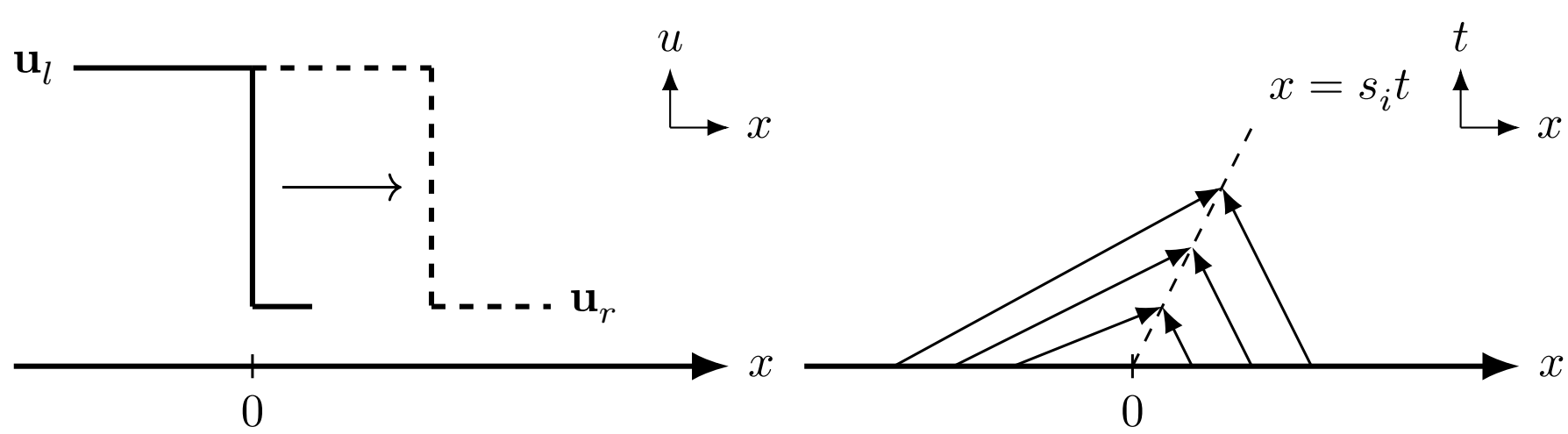

(a) Evolution of Solution (b) Characteristics

Figure 2.3.3: Shock Wave

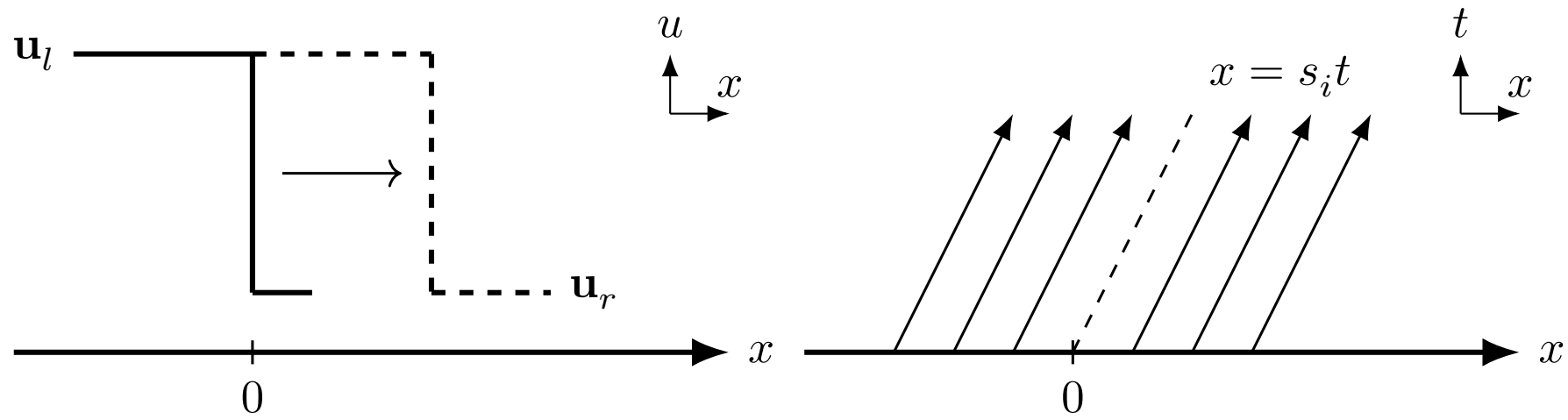

(a) Evolution of Solution (b) Characteristics

Figure 2.3.4: Contact Wave

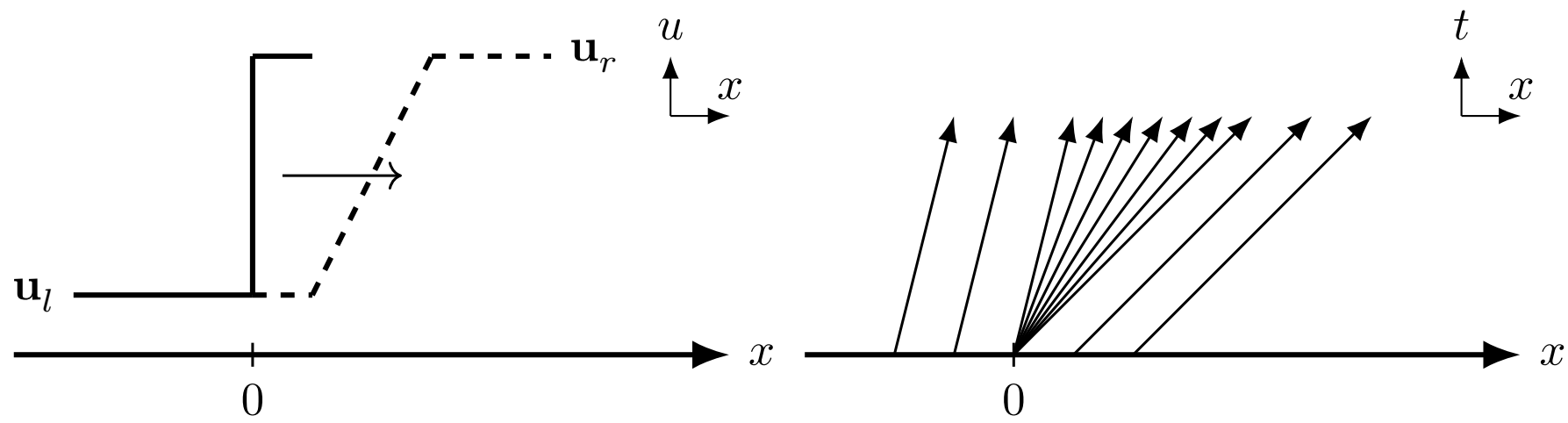

(a) Evolution of Solution (b) Characteristics

Figure 2.3.5: Rarefaction

2. **Contact Wave** The $\lambda_i$-wave is a contact wave, if it corresponds to a linearly degenerate field and connects two state $\mathbf{u}_l$ and $\mathbf{u}_r$ through a single jump discontinuity. The discontinuity moves with the speed $s_i$ given by the Rankine–Hugoniot condition Eq. (2.2.2). It satisfies the parallel characteristic condition

$$\lambda_i(\mathbf{u}_l) = s_i = \lambda_i(\mathbf{u}_r) \tag{2.3.3}$$

In other words, characteristic lines on either side of the contact line $\mathrm{d}x/\,\mathrm{d}t = s_i$ run parallel to it as shown in Fig. 2.3.4.

3. **Rarefaction** The $\lambda_i$-wave corresponds to a rarefaction, if it connects two states $\mathbf{u}_l$ and $\mathbf{u}_r$ through a smooth transition in a genuinely non-linear field. The characteristics corresponding to a rarefaction satisfy the property

$$\lambda_i(\mathbf{u}_l) \leq \lambda_i(\mathbf{u}_r) \tag{2.3.4}$$

and hence diverge from each other as shown in Fig. 2.3.5.

## 2.4 ALE Formulation for Hyperbolic Conservation Laws

Classically, the coordinate system used to describe hyperbolic conservation laws is assumed to be fixed in position respect to the laboratory. The corresponding frame of reference is known as a fixed frame or an Eulerian frame. The other frame we consider is the ALE frame where coordinate system can move arbitrarily with respect to the laboratory. In this section, we describe the transformation for physical quantities and their derivatives between the two kinematic frames. We then write down the hyperbolic conservation laws for an ALE frame.

**Definition 2.4.1** (ALE Frame)**.** *Let $(\mathbf{x}, t)$ describe the coordinate system for the Eulerian frame. We define the ALE frame by the coordinate system $(\boldsymbol{\chi}, t)$ where the relation between the ALE and the Eulerian coordinates is given by function $\mathbf{\Phi}$*

$$\mathbf{\Phi} : \Omega_{\boldsymbol{\chi}} \times [0, T) \to \Omega_{\mathbf{x}} \times [0, T) \tag{2.4.1a}$$

$$(\boldsymbol{\chi}, t) \mapsto \mathbf{\Phi}(\boldsymbol{\chi}, t) = (\mathbf{x}, t) \tag{2.4.1b}$$

*where $\mathbf{\Phi}$ is a homeomorphism and differentiable in $t$ almost everywhere over the interval $[0, T)$. We define the ALE Velocity by*

$$\mathbf{w} = \left.\frac{\partial \mathbf{x}}{\partial t}\right|_{\boldsymbol{\chi}} \tag{2.4.2}$$

**Theorem 2.4.1.** *The gradient for the transformation $\mathbf{\Phi}$ is*

$$\frac{\partial \mathbf{\Phi}}{\partial(\boldsymbol{\chi}, t)} = \begin{pmatrix} \frac{\partial \mathbf{x}}{\partial \boldsymbol{\chi}} & \mathbf{w} \\ \mathbf{0}^{\mathsf{T}} & 1 \end{pmatrix} \tag{2.4.3}$$

*Furthermore, the mapping $\mathbf{\Phi}$ is bijective and orientation-preserving if and only if*

$$\det\left(\frac{\partial \mathbf{x}}{\partial \boldsymbol{\chi}}\right) > 0 \qquad \forall t \in [0, T) \tag{2.4.4}$$

**Theorem 2.4.2** (Change of Variables)**.** *Let $\mathbf{u} \in \mathbb{R}^n$ be a physical quantity described by $\mathbf{u}(\mathbf{x}, t)$ and $\mathbf{u}(\boldsymbol{\chi}, t)$ in the Eulerian frame and ALE frame respectively. Then, we have*

$$\frac{\partial \mathbf{u}(\boldsymbol{\chi}, t)}{\partial \boldsymbol{\chi}} = \nabla_{\mathbf{x}}(\mathbf{x}, t)\mathbf{u}\, \frac{\partial \mathbf{x}}{\partial \boldsymbol{\chi}} \qquad \frac{\partial \mathbf{u}(\boldsymbol{\chi}, t)}{\partial t} = [\nabla_{\mathbf{x}}\mathbf{u}(\mathbf{x}, t)]\mathbf{w} + \frac{\partial \mathbf{u}(\mathbf{x}, t)}{\partial t} \tag{2.4.5}$$

**Theorem 2.4.3** (Reynold's Transport Theorem)**.** *Let $\Omega$ be a physical domain and let $\Omega_{\mathbf{x}}(t)$ be its description in Eulerian frame. Let $\mathbf{u} \in \mathbb{R}^n$ be a physical quantity described by $\mathbf{u}(\mathbf{x}, t)$ in the Eulerian frame. Then, we have*

$$\frac{\mathrm{d}}{\mathrm{d}t} \int\limits_{\Omega_{\mathbf{x}}(t)} \mathbf{u}(\mathbf{x}, t)\, \mathrm{d}\mathbf{x} = \int\limits_{\Omega_{\mathbf{x}}(t)} \left[\frac{\partial}{\partial t}\mathbf{u}(\mathbf{x}, t) + \nabla_{\mathbf{x}} \cdot [\mathbf{w} \otimes \mathbf{u}]\right] \mathrm{d}\mathbf{x} \tag{2.4.6}$$

**Theorem 2.4.4** (Hyperbolic Conservation Laws on Moving Domains)**.** *Let $\Omega$ be a physical domain and let $\Omega_{\mathbf{x}}(t)$ be its description in Eulerian frame. Let $\mathbf{u} \in \mathbb{R}^n$ be a physical quantity* *described by* $\mathbf{u}(\mathbf{x},t)$ *in the Eulerian frame. Then the integral form of the hyperbolic conservation law on the moving domains is given by*

$$\frac{\mathrm{d}}{\mathrm{d}t}\int\limits_{\Omega_{\mathbf{x}}(t)} \mathbf{u}(\mathbf{x},t)\,\mathrm{d}\mathbf{x} + \int\limits_{\Omega_{\mathbf{x}}(t)} \nabla_{\mathbf{x}}\cdot[\mathcal{F}(\mathbf{x},t,\mathbf{u}) - \mathbf{w}\otimes\mathbf{u}]\,\mathrm{d}\mathbf{x} = 0 \tag{2.4.7}$$

**Proof.** By the Reynold's transport theorem, Eq. (2.4.6), we have

$$\frac{\mathrm{d}}{\mathrm{d}t}\int\limits_{\Omega_{\mathbf{x}}(t)} \mathbf{u}(\mathbf{x},t)\,\mathrm{d}\mathbf{x} = \int\limits_{\Omega_{\mathbf{x}}(t)} \left[\frac{\partial}{\partial t}\mathbf{u}(\mathbf{x},t) + \nabla_{\mathbf{x}}\cdot[\mathbf{w}\otimes\mathbf{u}]\right]\mathrm{d}\mathbf{x}$$

Then, by using Eq. (2.1.3), we have

$$\frac{\mathrm{d}}{\mathrm{d}t}\int\limits_{\Omega_{\mathbf{x}}(t)} \mathbf{u}(\mathbf{x},t)\,\mathrm{d}\mathbf{x} = \int\limits_{\Omega_{\mathbf{x}}(t)} \left[-\nabla_{\mathbf{x}}\cdot\mathcal{F}(\mathbf{x},t,\mathbf{u}) + \nabla_{\mathbf{x}}\cdot[\mathbf{w}\otimes\mathbf{u}]\right]\mathrm{d}\mathbf{x}$$

which then gives us Eq. (2.4.7). ■

**Theorem 2.4.5** (ALE Formulation of Hyperbolic Conservation Laws)**.** *Let* $\Omega \subset \mathbb{R}^n$ *be a physical domain and let* $\Omega_{\chi}$ *be its description the ALE frame respectively. Let* $\mathbf{u} \in \mathbb{R}^n$ *be a physical quantity defined on* $\Omega$ *. Then, the ALE form of the hyperbolic conservation law is given by*

$$\frac{\mathrm{d}}{\mathrm{d}t}\int\limits_{\Omega_{\chi}} \mathbf{u}(\chi,t)\,\mathrm{d}\chi + \int\limits_{\Omega_{\chi}} \left[\nabla_{\chi}\cdot\mathcal{F}(\chi,t,\mathbf{u}) - [\nabla_{\chi}\mathbf{u}(\chi,t)]\mathbf{w}\right]\left|\frac{\partial\chi}{\partial\mathbf{x}}\right|\mathrm{d}\chi = 0 \tag{2.4.8}$$

$$\frac{\partial}{\partial t}\mathbf{u}(\chi,t) + \left[\nabla_{\chi}\cdot\mathcal{F}(\chi,t,\mathbf{u}) - [\nabla_{\chi}\mathbf{u}(\chi,t)]\mathbf{w}\right]\left|\frac{\partial\chi}{\partial\mathbf{x}}\right| = 0 \tag{2.4.9}$$

**Proof.** Since $\Omega_{\chi}$ is not a function of time, we can take the time derivative under the integral sign and use change of variables for time derivative, Eq. (2.4.5), to write

$$\frac{\mathrm{d}}{\mathrm{d}t}\int\limits_{\Omega_{\chi}} \mathbf{u}(\chi,t)\,\mathrm{d}\chi = \int\limits_{\Omega_{\chi}} \left[\left[\nabla_{\mathbf{x}}\mathbf{u}(\mathbf{x},t)\right]\mathbf{w} + \frac{\partial\mathbf{u}(\mathbf{x},t)}{\partial t}\right]\mathrm{d}\chi$$

Next, we can use Eq. (2.1.3) to get

$$\frac{\mathrm{d}}{\mathrm{d}t}\int\limits_{\Omega_{\chi}} \mathbf{u}(\chi,t)\,\mathrm{d}\chi = \int\limits_{\Omega_{\chi}} \left[[\nabla_{\mathbf{x}}\mathbf{u}(\mathbf{x},t)]\mathbf{w} - \nabla_{\mathbf{x}}\cdot\mathcal{F}(\mathbf{x},t,\mathbf{u})\right]\mathrm{d}\chi$$

Finally, we do a change of coordinates to convert everything in terms of $\chi$,

$$\frac{\mathrm{d}}{\mathrm{d}t}\int\limits_{\Omega_{\chi}} \mathbf{u}(\chi,t)\,\mathrm{d}\chi = \int\limits_{\Omega_{\chi}} \left[[\nabla_{\chi}\mathbf{u}(\chi,t)]\mathbf{w} - \nabla_{\chi}\cdot\mathcal{F}(\chi,t,\mathbf{u})\right]\left|\frac{\partial\chi}{\partial\mathbf{x}}\right|\mathrm{d}\chi$$

which after rearranging the terms gives us Eq. (2.4.8) Next, by again using the fact the $\Omega_{\chi}$ is independent of time, we can take the time derivative inside the integral sign for the first term.

And then by the arbitrariness of $\Omega_{\mathbf{x}}$, we get Eq. (2.4.9). ■

# Chapter 3 An ALE DG Scheme

In this chapter, we design an ALE DG method to solve the conservation law over a Lipschitz domain $\Omega(t) \subset \mathbb{R}^d$ in a single time step $[t_n, t_{n+1}]$

$$\frac{\partial}{\partial t}\mathbf{u}(\mathbf{x},t) + \nabla_{\mathbf{x}} \cdot \mathcal{F}(\mathbf{x},t,\mathbf{u}) = 0 \qquad (\mathbf{x},t) \in \Omega(t) \times [t_n, t_{n+1}) \qquad (3.0.1)$$

The Lipschitz domain $\Omega(t)$ at any time $t$ is partitioned into a simplex mesh $\mathcal{T}(t)$ at time $t$. Given the solution $\mathbf{u}(\mathbf{x},t)$ at $t_n$, the ALE DG method provides us with the velocity of the vertices. With the velocities at vertices, the velocities in the interior of the mesh is now computed using barycentric interpolation of the velocities of the vertices.

The ALE DG scheme then provides us with the time step $\Delta t_n = t_{n+1} - t_n$ that ensures that the topology of the mesh does not change between $t_n$ and $t_{n+1}$. We assume that the mesh velocity is constant between the time steps $t_n$ and $t_{n+1}$ and is equal to the mesh velocity computed at time $t_n$. This determines the space time discretization of the space in the time step $[t_n, t_{n+1}]$. We then solve the conservation law over the resulting space time configuration in the given time step.

Towards that end, we first describe the mesh and its motion in Section 3.1. Then, we describe the finite element space over the moving mesh in .

## 3.1 The Mesh for the ALE DG Scheme

The ALE DG method is designed to work on simplicial meshes. In this section, we introduce the simplicial mesh and describe its motion for the ALE DG method. We start by the definition of a simplex itself, the different elements in a simplex and the simplicial mesh.

**Definition 3.1.1** (Simplex, Faces, Edges and Vertices)**.** *Let $X$ be a collection of $d+1$ affinely independent points.*

1. ***$d$-dimensional Simplex :*** *A **$d$-dimensional simplex** $\tau$ is the convex hull of a set $X$ of $d+1$ affinely independent points. In particular, $0$-simplex is a point, $1$-simplex is an edge, $2$-simplex is a triangle, and $3$-simplex is a tetrahedron.*
2. ***Elements :*** *A simplex is an element of $\tau$ if it is a convex hull of a non-empty subset of $X$.*
3. ***Vertices :*** *The points in the set $X$ are known as the vertices of the simplex. One can also think of vertices as elements of the simplex $\tau$ with dimension $d$.*
4. ***Edges:*** *A simplex is an edge of $\tau$ if it is the convex hull of a non-empty subset of $X$ and has a dimension $1$.*

5. ***Faces:*** *A simplex is a face of $\tau$ if it is the convex hull of a non-empty subset of $X$ and has a dimension $k-1$.*

**Definition 3.1.2** (Simplicial Complex)**.** *A set $\mathcal{T}$ of finitely many simplexes is called a **simplicial complex** if*

- *$\mathcal{T}$ contains every element of every simplex in $\mathcal{T}$*
- *For any two simplexes $\sigma, \tau \in \mathcal{T}$, $\sigma \cap \tau$ is either empty or an element of both $\sigma$ and $\tau$.*

*The $d$ dimensional simplexes in $\mathcal{T}$ are called **cells** and the set of cells is denoted by $\mathcal{K}$. All the $d-1$ dimensional elements of the simplexes in $\mathcal{K}$ are defined as **faces** and the set of faces is denoted by $\mathcal{F}$. The $0$ dimensional elements of the simplex in $\mathcal{K}$ are defined as the **vertices** of the complex and the set of vertices is known as $V$.*

**Definition 3.1.3** (Orientation of A Cell in a Complex)**.** *Let $K_i$ be a cell in a simplicial complex which is the convex hull of the subset $V_{K_i}$ of vertices $V$. Associate a permutation $\sigma = (v_0, v_1, \dots, v_d)$, $v_i \in V_{K_i}$ of the vertices and call the resultant cell as $K_{i,\sigma}$. Then one can define an orientation of the simplex $K_{i,\sigma}$ as the determinant of the matrix $A_\sigma$, where $A_\sigma$ is defined as*

$$A_\sigma = \begin{pmatrix} v_1 - v_0 & v_2 - v_1 & \cdots & v_d - v_{d-1} \end{pmatrix} \tag{3.1.1}$$

*By definition, the orientation can either be positive or negative. In future, we will assume that all cells have an associated permutation without explicitly specifying it in the subscript and just write $K_i$ instead of $K_{i,\sigma}$.*

**Definition 3.1.4** (Simplicial Mesh)**.** *Let $V$ be the set of points in $\mathbb{R}^d$. A $d$-dimensional simplicial mesh of $V$ is a simplicial complex $\mathcal{T}$ such that*

1. *$V$ is the set of vertices in $\mathcal{T}$.*
2. *The union of all the simplexes in $\mathcal{T}$ is the convex hull of $V$.*
3. *Orientation of all cells $K_i \in \mathcal{K} \subset \mathcal{T}$ is positive.*

## 3.2 Evolution of The Mesh during a Time Step

Consider an ALE DG simulation starting at time $t_n$ with the time step determined to be $[t_n, t_{n+1}]$. Assume that a simplicial mesh $\mathcal{T}(t_n)$ exists at time $t_n$. An example of such a mesh in two dimensions is shown in **??**.

Let $V_{\mathcal{T}}(t)$ be the set of vertices of the mesh. Given the solution $\mathbf{u}(\mathbf{x}, t)$ at $t_n$, the ALE DG method provides us with the velocity of the vertices in $V_{\mathcal{T}}(t_n)$. We denote these velocities by $\mathbf{w}_{\mathbf{V}_i}$. Given the velocities at vertices, the velocities in the interior of the mesh is now computed using barycentric interpolation of the velocities of the vertices

$$\mathbf{w}_{K_i}(\mathbf{x}) = \sum_{j=1}^{d+1} \alpha_j \mathbf{w}_{v_{k_j}} \tag{3.2.1}$$

where $\boldsymbol{\alpha} = (\alpha_0, \alpha_1, \dots, \alpha_{d+1})$ is the barycentric co-ordinate of $\mathbf{x}$ in the cell $K_i$. The ALE DG scheme provides us with the time step $\Delta t_n = t_{n+1} - t_n$ t and vertex velocities such that ensures that the topology of the mesh does not change between $t_n$ and $t_{n+1}$. We assume that the mesh velocity is constant between the time steps $t_n$ and $t_{n+1}$ and is equal to the mesh velocity computed at time $t_n$. This ensures that the mesh motion is affine and the mesh remains a simplicial mesh. It is also important to ensure that the orientation of the cells is also preserved during the evolution of the mesh. This can be achieved by placing the following conditions on the evolution of the mesh.

**Theorem 3.2.1** (Preservation of Orientation)**.** *Let $\mathcal{K}(t_n)$ be the mesh at time $t_n$ and let $\mathcal{K}(t_{n+1})$ be the mesh at time $t_{n+1}$. Let $K_i^n$ be the cells in $\mathcal{K}(t_n)$ and let $K_i^{n+1}$ be the cells in $\mathcal{K}(t_{n+1})$. Let $A_{K_i}^n$ and $A_{K_i}^{n+1}$ be the orientation matrices of the cells. Further assume that*

$$\|\mathbf{A}_{K^n}^{-1}\mathbf{A}_{K^{n+1}}\|_{L(\mathbb{R}^d,\mathbb{R}^d)} \leq \frac{t_{n+1}-t}{t-t_n} \qquad \|\mathbf{A}_{K^{n+1}}^{-1}\mathbf{A}_{K^n}\|_{L(\mathbb{R}^d,\mathbb{R}^d)} \leq \frac{t-t_n}{t_{n+1}-t} \tag{3.2.2}$$

## 3.3 Finite Element Space for ALE DG Scheme

*Notation* 1 (Restriction of a Function). Let $\Omega \subset \mathbb{R}^d$ be a Lipschitz domain. Given a function $v$ defined on $\mathbb{R}^d$, we define

$$v|_\Omega \stackrel{\text{def}}{=} v(\mathbf{x}) \cdot \mathbb{1}_\Omega \tag{3.3.1}$$

**Definition 3.3.1** (Test Function Space of Degree $k$)**.** *Let $\mathcal{T}$ be a simplicial mesh. Then the **test function space of degree** $k$, $V_{\mathcal{T},k}$, is defined as*

$$V_{\mathcal{T}} \stackrel{\text{def}}{=} \left\{v \in L^2(\Omega) : v|_K \in \mathbb{P}^k \,\forall K \in \mathcal{T}\right\} \tag{3.3.2}$$

*where $\mathbb{P}^k$ is space of polynomials of degree at most $k$.*

**Theorem 3.3.1** (Transport Function for Test Functions)**.** *Let $u\colon \Omega \times [t_n, t_{n+1}] \to \mathbb{R}$ be sufficiently smooth function. Then, for all $t \in [t_n, t_{n+1}]$, for all $K(t) \in \mathcal{T}(t)$ and for all $v \in V_{\mathcal{T}}$, we have ,*

$$\frac{\mathrm{d}}{\mathrm{d}t}(u, v)_{K(t)} = \left(\frac{\partial u}{\partial t}, v\right) + (\nabla\cdot(\mathbf{w} \otimes u), v)_{K(t)} \tag{3.3.3}$$

## 3.4 Semi-Discrete Formulation

We first define the ALE flux which we will need later in our computations.

**Definition 3.4.1** (ALE Flux)**.** *Define the ALE flux $\mathcal{G}$ as*

$$\mathcal{G}(\mathbf{x}, t, \mathbf{u}, \mathbf{w}) = \mathcal{F}(\mathbf{x}, t, \mathbf{u}) - \mathbf{u} \otimes \mathbf{w} \tag{3.4.1}$$

At any time $t$, the approximate solution $\mathbf{u}(\mathbf{x}, t)$ in some grid cell $K(t)$ is completely specified by the modes $(\mathbf{u}_m(t))$, $m = 1, 2, \dots, M$ of the solution as

$$\mathbf{u}(\mathbf{x}, t) = \sum_{m=1}^{M} \mathbf{u}_m(t)\phi_m(\mathbf{x}, t) \tag{3.4.2}$$

We can transform the solution to the reference cell using the transformation from previous section, when we can write

$$\mathbf{u}(\mathbf{x},t) = \hat{\mathbf{u}}(\boldsymbol{\xi},\tau) = \sum_{m=1}^{M} \hat{\mathbf{u}}_m(\tau)\hat{\phi}_m(\boldsymbol{\xi}) \tag{3.4.3}$$

Next, we can use the orthonormality of the basis functions to write

$$\mathbf{u}_m(t) = \hat{\mathbf{u}}_m(\tau) = \int_{K(t)} \mathbf{u}(\mathbf{x},t)\phi_m(\mathbf{x},t)\mathrm{d}\mathbf{x} = \int_{\widehat{K}} \hat{\mathbf{u}}(\boldsymbol{\xi},\tau)\hat{\phi}_m(\boldsymbol{\xi},\tau)\left|\boldsymbol{\mathcal{J}}_{K(t)}\right|\mathrm{d}\boldsymbol{\xi} \tag{3.4.4}$$

Whence, we can compute the evolution of the moments as

$$\begin{aligned}\frac{\mathrm{d}}{\mathrm{d}t}\int_{K(t)} \mathbf{u}(\mathbf{x},t)\phi_m(\mathbf{x},t)\mathrm{d}\mathbf{x} &= \frac{\mathrm{d}}{\mathrm{d}t}\int_{\widehat{K}} \hat{\mathbf{u}}(\boldsymbol{\xi},\tau)\phi_m(\boldsymbol{\xi})\left|\boldsymbol{\mathcal{J}}_{K(t)}\right|\mathrm{d}\boldsymbol{\xi} \\ &= \int_{\widehat{K}} \frac{\partial}{\partial\tau}\left[\hat{\mathbf{u}}(\boldsymbol{\xi},\tau)\phi_m(\boldsymbol{\xi})\left|\boldsymbol{\mathcal{J}}_{K(t)}\right|\right]\mathrm{d}\boldsymbol{\xi} \\ &= \int_{\widehat{K}} \phi_m(\boldsymbol{\xi})\left[\left|\boldsymbol{\mathcal{J}}_{K(t)}\right|\frac{\partial}{\partial\tau}\hat{\mathbf{u}}(\boldsymbol{\xi},\tau) + \hat{\mathbf{u}}(\boldsymbol{\xi},\tau)\frac{\partial}{\partial\tau}\left|\boldsymbol{\mathcal{J}}_{K(t)}\right|\right]\mathrm{d}\boldsymbol{\xi}\end{aligned}$$

$$\frac{\partial}{\partial\tau}\hat{\mathbf{u}}(\boldsymbol{\xi},\tau) = \frac{\partial}{\partial t}\mathbf{u}(\mathbf{x},t) + \nabla_{\mathbf{x}}\mathbf{u}\cdot\mathbf{w} \tag{3.4.5a}$$

$$\frac{\partial}{\partial\tau}\left|\boldsymbol{\mathcal{J}}_{K(t)}\right| = \left|\boldsymbol{\mathcal{J}}_{K(t)}\right|\nabla_{\mathbf{x}}\cdot\mathbf{w} \tag{3.4.5b}$$

$$\begin{aligned}\frac{\mathrm{d}}{\mathrm{d}t}\int_{K(t)} \mathbf{u}(\mathbf{x},t)\phi_m(\mathbf{x},t)\mathrm{d}\mathbf{x} &= \int_{K(t)} \phi_m(\mathbf{x},t)\left[\left(\frac{\partial}{\partial t}\mathbf{u}(\mathbf{x},t) + \nabla_{\mathbf{x}}\mathbf{u}\cdot\mathbf{w}\right) + \mathbf{u}(\mathbf{x},t)\nabla_{\mathbf{x}}\cdot\mathbf{w}\right]\mathrm{d}\mathbf{x} \\ &= \int_{K(t)} \phi_m(\mathbf{x},t)\nabla_{\mathbf{x}}\cdot\left[\mathbf{u}\otimes\mathbf{w} - \boldsymbol{\mathcal{F}}(\mathbf{x},t,\mathbf{u})\right]\mathrm{d}\mathbf{x}\end{aligned}$$

Hence, the equation for the evolution of moment is given by

$$\frac{\mathrm{d}}{\mathrm{d}t}\int_{K(t)} \mathbf{u}(\mathbf{x},t)\phi_m(\mathbf{x},t)\mathrm{d}\mathbf{x} = -\int_{K(t)} \nabla_{\mathbf{x}}\cdot\left[\boldsymbol{\mathcal{F}}(\mathbf{x},t,\mathbf{u}) - \mathbf{u}\otimes\mathbf{w}\right]\phi_m(\mathbf{x},t)\mathrm{d}\mathbf{x} \tag{3.4.6}$$

Apply integration by parts to get

$$\begin{aligned}\frac{\mathrm{d}}{\mathrm{d}t}\int_{K(t)} \mathbf{u}(\mathbf{x},t)\phi_m(\mathbf{x},t)\mathrm{d}\mathbf{x} = &\int_{K(t)} \left[\boldsymbol{\mathcal{F}}(\mathbf{x},t,\mathbf{u}) - \mathbf{u}\otimes\mathbf{w}\right]\cdot\nabla_{\mathbf{x}}\phi_m(\mathbf{x},t)\mathrm{d}\mathbf{x} \\ &- \int_{\partial K(t)} \left(\boldsymbol{\mathcal{F}}(\mathbf{x},t,\mathbf{u}) - \mathbf{u}\otimes\mathbf{w}\right)\cdot\hat{\boldsymbol{\nu}}\phi_m(\mathbf{s},t)\mathrm{d}\mathbf{s}\end{aligned} \tag{3.4.7}$$

We now use the ALE version $\mathcal{G}$ of the analytical flux $\mathcal{F}$ to get

$$\frac{\mathrm{d}}{\mathrm{d}t}\int_{K(t)} \mathbf{u}(\mathbf{x},t)\phi_m(\mathbf{x},t)\mathrm{d}\mathbf{x} = \int_{K(t)} \mathcal{G}(\mathbf{x},t,\mathbf{u},\mathbf{w})\cdot\nabla_{\mathbf{x}}\phi_m(\mathbf{x},t)\,\mathrm{d}\mathbf{x} - \int_{\partial K(t)} \mathcal{G}(\mathbf{u},\mathbf{w})\cdot\hat{\boldsymbol{\nu}}\phi_m(\mathbf{s},t)\mathrm{d}\mathbf{s} \tag{3.4.8}$$

The solution is discontinuous at the faces, and hence, we must instead use a numerical flux to evaluate the solution at the boundary. Let $\mathcal{H}(\mathbf{x},t,\mathbf{u}^-,\mathbf{u}^+,\mathbf{w};\hat{\boldsymbol{\nu}})$ be the numerical flux function for $\mathcal{G}(\mathbf{x},t,\mathbf{u},\mathbf{w})$ which is consistent in the sense that

$$\mathcal{H}(\mathbf{x},t,\mathbf{u},\mathbf{u},\mathbf{w};\hat{\boldsymbol{\nu}}) = \mathcal{G}(\mathbf{x},t,\mathbf{u},\mathbf{w})\cdot\hat{\boldsymbol{\nu}}, \qquad \forall\mathbf{u},\mathbf{w},\hat{\boldsymbol{\nu}} \tag{3.4.9}$$

Then we can rewrite the equation in the semi-discrete form

$$\frac{\mathrm{d}}{\mathrm{d}t}\int_{K(t)} \mathbf{u}(\mathbf{x},t)\phi_m(\mathbf{x},t)\mathrm{d}\mathbf{x} = \int_{K(t)} \mathcal{G}(\mathbf{x},t,\mathbf{u},\mathbf{w})\cdot\nabla_{\mathbf{x}}\phi_m(\mathbf{x},t)\mathrm{d}\mathbf{x} - \int_{\partial K(t)} \mathcal{H}(\mathbf{x},t,\mathbf{u}^-,\mathbf{u}^+,\mathbf{w};\hat{\boldsymbol{\nu}})\phi_m(\mathbf{s},t)\mathrm{d}\mathbf{s} \tag{3.4.10}$$

## 3.5 Properties of ALE DG Scheme

In this section, we list some of the properties of the ALE DG scheme. In particular, we show that the ALE DG scheme is conservative and preserves constant states over arbitrary mesh velocities.

### 3.5.1 Conservation

The semi-discrete scheme can be integrated over time as

$$\begin{aligned}\int_{K(t_{n+1})} &\mathbf{u}(\mathbf{x},t_{n+1})\phi_m(\mathbf{x},t_{n+1})\,\mathrm{d}\mathbf{x} - \int_{K(t_n)} \mathbf{u}(\mathbf{x},t_n)\phi_m(\mathbf{x},t_n)\,\mathrm{d}\mathbf{x} \\ &= \int_{t}^{t_{n+1}}\int_{K(t)} \mathcal{G}(\mathbf{x},t,\mathbf{u},\mathbf{w})\cdot\nabla_{\mathbf{x}}\phi_m(\mathbf{x},t)\,\mathrm{d}\mathbf{x}\,\mathrm{d}t - \int_{t}^{t_{n+1}}\int_{\partial K(t)} \mathcal{H}(\mathbf{x},t,\mathbf{u}^-,\mathbf{u}^+,\mathbf{w};\hat{\boldsymbol{\nu}})\phi_m(\mathbf{s},t)\,\mathrm{d}\mathbf{s}\,\mathrm{d}t\end{aligned} \tag{3.5.1}$$

For the case $m = 0$, $\phi_m(\mathbf{x}) = \phi_0$ is a constant w.r.t. $\mathbf{x}$. Hence we can write

$$\begin{aligned}\int_{K(t_{n+1})} \mathbf{u}(\mathbf{x},t_{n+1})\,\mathrm{d}\mathbf{x} &- \int_{K(t_n)} \mathbf{u}(\mathbf{x},t_n)\,\mathrm{d}\mathbf{x} \\ &- \int_{t}^{t_{n+1}}\int_{\partial K(t)} \mathcal{H}(\mathbf{x},t,\mathbf{u}^-,\mathbf{u}^+,\mathbf{w};\hat{\boldsymbol{\nu}})\,\mathrm{d}\mathbf{s}\,\mathrm{d}t\end{aligned} \tag{3.5.2}$$

Summing over all the cells, we get

$$\sum_{K_i} \int_{K(t_{n+1})} \mathbf{u}(\mathbf{x}, t_{n+1})\,\mathrm{d}\mathbf{x} - \sum_{K_i} \int_{K(t_n)} \mathbf{u}(\mathbf{x}, t_n)\,\mathrm{d}\mathbf{x} = -\sum_{K_i} \int_{t}^{t_{n+1}} \int_{\partial K(t)} \mathcal{H}(\mathbf{x}, t, \mathbf{u}^-, \mathbf{u}^+, \mathbf{w}; \hat{\boldsymbol{\nu}})\,\mathrm{d}\mathbf{s}\,\mathrm{d}t \tag{3.5.3}$$

Next, we use the fact that

$$\mathcal{H}(\mathbf{x}, t, \mathbf{u}^-, \mathbf{u}^+, \mathbf{w}; \hat{\boldsymbol{\nu}}) = -\mathcal{H}(\mathbf{x}, t, \mathbf{u}^+, \mathbf{u}^-, \mathbf{w}; -\hat{\boldsymbol{\nu}}) \tag{3.5.4}$$

and that for $\partial K_i$ such that $\partial K_i = K_i \cap K_j$, we have

$$\mathcal{H}_{K_i}(\mathbf{x}, t, \mathbf{u}^-, \mathbf{u}^+, \mathbf{w}; \hat{\boldsymbol{\nu}}) = -\mathcal{H}_{K_j}(\mathbf{x}, t, \mathbf{u}^-, \mathbf{u}^+, \mathbf{w}; -\hat{\boldsymbol{\nu}}) \tag{3.5.5}$$

$$\sum_{K_i} \int_{K(t_{n+1})} \mathbf{u}(\mathbf{x}, t_{n+1})\,\mathrm{d}\mathbf{x} - \sum_{K_i} \int_{K(t_n)} \mathbf{u}(\mathbf{x}, t_n)\,\mathrm{d}\mathbf{x} = -\sum_{K_i} \int_{t}^{t_{n+1}} \int_{\partial K(t)} \mathcal{H}(\mathbf{x}, t, \mathbf{u}^-, \mathbf{u}^+, \mathbf{w}; \hat{\boldsymbol{\nu}})\,\mathrm{d}\mathbf{s}\,\mathrm{d}t \tag{3.5.6}$$

### 3.5.2 Preservation of Constant States

Let the solution be a constant, i.e. $\mathbf{u}^n = \mathbf{c}$. Then any consistent predictor will also be constant $\tilde{\mathbf{u}}(\mathbf{x}, t) = \mathbf{c}$. The semi-discrete scheme can be integrated over time as

$$\int_{K^{n+1}} \mathbf{u}(\mathbf{x}, t_{n+1})\phi_m(\mathbf{x}, t_{n+1})\,\mathrm{d}\mathbf{x} - \int_{K^n} \mathbf{c}\phi_m(\mathbf{x}, t_n)\,\mathrm{d}\mathbf{x} = \int_{t_n}^{t_{n+1}} \int_{K(t)} \mathcal{G}(\mathbf{x}, t, \mathbf{c}, \mathbf{w}) \cdot \nabla_{\mathbf{x}}\phi_m(\mathbf{x}, t)\,\mathrm{d}\mathbf{x}\,\mathrm{d}t - \int_{t_n}^{t_{n+1}} \int_{\partial K(t)} \mathcal{H}(\mathbf{x}, t, \mathbf{c}, \mathbf{c}, \mathbf{w}; \hat{\boldsymbol{\nu}})\phi_m(\mathbf{s}, t)\,\mathrm{d}\mathbf{s}\,\mathrm{d}t \tag{3.5.7}$$

We use the consistency of numerical flux to write

$$\int_{K^{n+1}} \mathbf{u}(\mathbf{x}, t_{n+1})\phi_m(\mathbf{x}, t_{n+1})\,\mathrm{d}\mathbf{x} - \int_{K^n} \mathbf{c}\phi_m(\mathbf{x}, t_n)\,\mathrm{d}\mathbf{x} = \int_{t_n}^{t_{n+1}} \int_{K(t)} \mathcal{G}(\mathbf{x}, t, \mathbf{c}, \mathbf{w}) \cdot \nabla_{\mathbf{x}}\phi_m(\mathbf{x}, t)\,\mathrm{d}\mathbf{x}\,\mathrm{d}t - \int_{t_n}^{t_{n+1}} \int_{\partial K(t)} \phi_m(\mathbf{s}, t)\mathcal{G}(\mathbf{x}, t, \mathbf{c}, \mathbf{w}) \cdot \boldsymbol{\nu}\,\mathrm{d}\mathbf{s}\,\mathrm{d}t \tag{3.5.8}$$

Using Gauss Divergence Theorem, we can write

$$\int_{K^{n+1}} \mathbf{u}(\mathbf{x}, t_{n+1})\phi_m(\mathbf{x}, t_{n+1})\,\mathrm{d}\mathbf{x} - \int_{K^n} \mathbf{c}\phi_m(\mathbf{x}, t_n)\,\mathrm{d}\mathbf{x}$$
$$= \int_{t_n}^{t_{n+1}} \int_{K(t)} \mathcal{G}(\mathbf{x}, t, \mathbf{c}, \mathbf{w}) \cdot \nabla_{\mathbf{x}}\phi_m(\mathbf{x}, t)\,\mathrm{d}\mathbf{x}\,\mathrm{d}t - \int_{t_n}^{t_{n+1}} \int_{K(t)} \nabla_{\mathbf{x}} \cdot \Big[\mathcal{G}(\mathbf{x}, t, \mathbf{c}, \mathbf{w})\phi_m(\mathbf{x}, t)\Big]\,\mathrm{d}\mathbf{x}\,\mathrm{d}t \tag{3.5.9}$$

which can now be written as

$$\int_{K^{n+1}} \mathbf{u}(\mathbf{x}, t_{n+1})\phi_m(\mathbf{x}, t_{n+1})\,\mathrm{d}\mathbf{x} - \int_{K^n} \mathbf{c}\phi_m(\mathbf{x}, t_n)\,\mathrm{d}\mathbf{x}$$
$$= \int_{t_n}^{t_{n+1}} \int_{K(t)} \mathcal{G}(\mathbf{x}, t, \mathbf{c}, \mathbf{w})\cdot\nabla_{\mathbf{x}}\phi_m(\mathbf{x}, t)\,\mathrm{d}\mathbf{x}\,\mathrm{d}t - \int_{t_n}^{t_{n+1}} \int_{K(t)} \phi_m(\mathbf{x})\nabla_{\mathbf{x}}\cdot\mathcal{G}(\mathbf{x}, t, \mathbf{c}, \mathbf{w}) + \mathcal{G}(\mathbf{x}, t, \mathbf{c}, \mathbf{w})\cdot\nabla_{\mathbf{x}}\phi_m(\mathbf{x}, t)$$

On simplification, we get

$$\int_{K^{n+1}} \mathbf{u}(\mathbf{x}, t_{n+1})\phi_m(\mathbf{x}, t_{n+1})\,\mathrm{d}\mathbf{x} - \int_{K^n} \mathbf{c}\phi_m(\mathbf{x}, t_n)\,\mathrm{d}\mathbf{x} = -\int_{t_n}^{t_{n+1}} \int_{K(t)} \phi_m(\mathbf{x}, t)\nabla_{\mathbf{x}} \cdot \mathcal{G}(\mathbf{x}, t, \mathbf{c}, \mathbf{w})\,\mathrm{d}\mathbf{x}\,\mathrm{d}t \tag{3.5.11}$$

We can now expand the flux and write

$$\int_{K^{n+1}} \mathbf{u}(\mathbf{x}, t_{n+1})\phi_m(\mathbf{x}, t_{n+1})\,\mathrm{d}\mathbf{x} - \int_{K^n} \mathbf{c}\phi_m(\mathbf{x}, t_n)\,\mathrm{d}\mathbf{x} = -\int_{t_n}^{t_{n+1}} \int_{K(t)} \phi_m(\mathbf{x}, t)\nabla_{\mathbf{x}} \cdot \Big[\mathcal{F}(\mathbf{x}, t, \mathbf{c}) - \mathbf{c} \otimes \mathbf{w}\Big]\,\mathrm{d}\mathbf{x} \tag{3.5.12}$$

If we have $\nabla_{\mathbf{x}} \cdot \mathcal{F}(\mathbf{x}, t, \mathbf{c}) = 0$, then

$$\int_{K^{n+1}} \mathbf{u}(\mathbf{x}, t_{n+1})\phi_m(\mathbf{x}, t_{n+1})\,\mathrm{d}\mathbf{x} = \int_{K^n} \mathbf{c}\phi_m(\mathbf{x}, t_n)\,\mathrm{d}\mathbf{x} + \mathbf{c}\left[\int_{t_n}^{t_{n+1}} \int_{K(t)} \phi_m(\mathbf{x}, t)\nabla_{\mathbf{x}} \cdot \mathbf{w}\,\mathrm{d}\mathbf{x}\,\mathrm{d}t\right] \tag{3.5.13}$$

We now consider two cases:

$\boxed{\text{Case} : m = 0}$ $\phi_m(\mathbf{x}, t) = d$, and hence we have

$$\int_{K^{n+1}} \mathbf{u}(\mathbf{x}, t_{n+1})d\,\mathrm{d}\mathbf{x} = \int_{K^n} \mathbf{c}d\,\mathrm{d}\mathbf{x} + \mathbf{c}\left[\int_{t_n}^{t_{n+1}} \int_{K(t)} d\nabla_{\mathbf{x}} \cdot \mathbf{w}\,\mathrm{d}\mathbf{x}\,\mathrm{d}t\right] \tag{3.5.14}$$

and we get

$$\int_{K^{n+1}} \mathbf{u}(\mathbf{x}, t_{n+1})\, \mathrm{d}\mathbf{x} = \mathbf{c} \left[ \int_{K^n} \mathrm{d}\mathbf{x} + \int_{t_n}^{t_{n+1}} \int_{K(t)} \nabla_{\mathbf{x}} \cdot \mathbf{w}\, \mathrm{d}\mathbf{x}\, \mathrm{d}t \right] = \left|K^{n+1}\right| \mathbf{c} \tag{3.5.15}$$

Case : $m \neq 0$ Using the fact that $\nabla_{\mathbf{x}} \cdot \mathbf{w} =$ constant and $\int_K \phi_m(\mathbf{x}, t)\, \mathrm{d}\mathbf{x} = 0$, we get

$$\int_{K^{n+1}} \mathbf{u}(\mathbf{x}, t_{n+1}) \phi_m(\mathbf{x}, t_{n+1})\, \mathrm{d}\mathbf{x} = 0 \tag{3.5.16}$$

The two cases together imply $\mathbf{u}_h^{n+1} = \mathbf{c}$ .

## 3.6 Fully Discrete Form

We now integrate over time to get

$$\begin{aligned} &\int_{K^{n+1}} \mathbf{u}(\mathbf{x}, t_{n+1}) \phi_m(\mathbf{x}, t_{n+1}) \mathrm{d}\mathbf{x} - \int_{K^n} \mathbf{u}(\mathbf{x}, t_n) \phi_m(\mathbf{x}, t_n) \mathrm{d}\mathbf{x} \\ &= \int_{t_n}^{t_{n+1}} \int_{K(t)} \mathcal{G}(\mathbf{x}, t, \mathbf{u}, \mathbf{w}) \cdot \nabla_{\mathbf{x}} \phi_m(\mathbf{x}) \mathrm{d}\mathbf{x} \mathrm{d}t - \int_{t_n}^{t_{n+1}} \int_{\partial K(t)} \mathcal{H}(\mathbf{x}, t, \mathbf{u}^-, \mathbf{u}^+, \mathbf{w}; \hat{\boldsymbol{\nu}}) \phi_m(\mathbf{s}) \mathrm{d}s \mathrm{d}t \end{aligned} \tag{3.6.1}$$

Since we are using a triangular mesh, we can divide the boundary term as a sum over the faces and then, we get

$$\begin{aligned} &\int_{K^{n+1}} \mathbf{u}(\mathbf{x}, t_{n+1}) \phi_m(\mathbf{x}, t_{n+1}) \mathrm{d}\mathbf{x} - \int_{K^n} \mathbf{u}(\mathbf{x}, t_n) \phi_m(\mathbf{x}, t_n) \mathrm{d}\mathbf{x} \\ &= \int_{t_n}^{t_{n+1}} \int_{K(t)} \mathcal{G}(\mathbf{x}, t, \mathbf{u}, \mathbf{w}) \cdot \nabla_{\mathbf{x}} \phi_m(\mathbf{x}) \mathrm{d}\mathbf{x} \mathrm{d}t - \int_{t_n}^{t_{n+1}} \sum_{F(t)} \int_{F(t)} \mathcal{H}(\mathbf{x}, t, \mathbf{u}^-, \mathbf{u}^+, \mathbf{w}; \hat{\boldsymbol{\nu}}) \phi_m(\mathbf{s}) \mathrm{d}s \mathrm{d}t \end{aligned} \tag{3.6.2}$$

Using the orthonormality of basis function and the polynomial nature of solution, we can write

$$\begin{aligned} &|\mathcal{J}_{K^{n+1}}|\, \mathbf{u}_m(t_{n+1}) - |\mathcal{J}_{K^n}|\, \mathbf{u}_m(t_n) \\ &= \int_{t_n}^{t_{n+1}} \int_{K(t)} \mathcal{G}(\mathbf{x}, t, \mathbf{u}, \mathbf{w}) \cdot \nabla_{\mathbf{x}} \phi_m(\mathbf{x}) \mathrm{d}\mathbf{x} \mathrm{d}t - \int_{t_n}^{t_{n+1}} \sum_{F(t)} \int_{F(t)} \mathcal{H}(\mathbf{x}, t, \mathbf{u}^-, \mathbf{u}^+, \mathbf{w}; \hat{\boldsymbol{\nu}}) \phi_m(\mathbf{s}) \mathrm{d}s \mathrm{d}t \end{aligned} \tag{3.6.3}$$

This is an implicit scheme which we can convert to an explicit scheme by using a predictor function $\tilde{\mathbf{u}}(x,t)$

$$
\begin{aligned}
&|\mathcal{J}_{K^{n+1}}|\,\mathbf{u}_m(t_{n+1}) - |\mathcal{J}_{K^n}|\,\mathbf{u}_m(t_n) \\
&= \int_{t_n}^{t_{n+1}} \int_{K(t)} \mathcal{G}(\mathbf{x},t,\tilde{\mathbf{u}},\mathbf{w}) \cdot \nabla_{\mathbf{x}}\phi_m(\mathbf{x})\mathrm{d}\mathbf{x}\mathrm{d}t - \int_{t_n}^{t_{n+1}} \sum_{F(t)} \int_{F(t)} \mathcal{H}(\mathbf{x},t,\tilde{\mathbf{u}}^-,\tilde{\mathbf{u}}^+,\mathbf{w};\hat{\boldsymbol{\nu}})\phi_m(\mathbf{s})\mathrm{d}\mathbf{s}\mathrm{d}t
\end{aligned}
\tag{3.6.4}
$$

In order to obtain a fully discrete scheme, we next convert the integrals in space to reference coordinates

$$
\begin{aligned}
|\mathcal{J}_{K^{n+1}}|\,\mathbf{u}_m(t_{n+1}) - |\mathcal{J}_{K^n}|\,\mathbf{u}_m(t_n) = &\int_{t_n}^{t_{n+1}} \left|\mathcal{J}_{K(t)}\right| \int_{\widehat{K}} \mathcal{G}(\mathbf{x},t,\tilde{\mathbf{u}},\mathbf{w}) \cdot \left[\nabla_{\boldsymbol{\xi}}\phi_m(\boldsymbol{\xi}) \cdot \mathcal{J}_{K(t)}^{-1}\right]\mathrm{d}\boldsymbol{\xi}\mathrm{d}t \\
&- \int_{t_n}^{t_{n+1}} \sum_{F(t)} \left|\mathcal{J}_{F(t)}\right| \int_{\hat{I}} \mathcal{H}(\mathbf{x},t,\tilde{\mathbf{u}}^-,\tilde{\mathbf{u}}^+,\mathbf{w};\hat{\boldsymbol{\nu}})\phi_m(\boldsymbol{\sigma})\mathrm{d}\boldsymbol{\sigma}\mathrm{d}t
\end{aligned}
\tag{3.6.5}
$$

Simplification of Jacobian gives us

$$
\begin{aligned}
|\mathcal{J}_{K^{n+1}}|\,\mathbf{u}_m(t_{n+1}) - |\mathcal{J}_{K^n}|\,\mathbf{u}_m(t_n) = &\int_{t_n}^{t_{n+1}} \int_{\widehat{K}} \mathcal{G}(\mathbf{x},t,\tilde{\mathbf{u}},\mathbf{w}) \cdot \left[\nabla_{\boldsymbol{\xi}}\phi_m(\boldsymbol{\xi}) \cdot \widehat{\mathcal{J}}_{K(t)}^{-1}\right]\mathrm{d}\boldsymbol{\xi}\mathrm{d}t \\
&- \int_{t_n}^{t_{n+1}} \sum_{F(t)} \left|\mathcal{J}_{F(t)}\right| \int_{\hat{I}} \mathcal{H}(\mathbf{x},t,\tilde{\mathbf{u}}^-,\tilde{\mathbf{u}}^+,\mathbf{w};\hat{\boldsymbol{\nu}})\phi_m(\boldsymbol{\sigma})\mathrm{d}\boldsymbol{\sigma}\mathrm{d}t
\end{aligned}
\tag{3.6.6}
$$

Next, we plug in the quadrature rules on the right hand side to write

$$
\begin{aligned}
&|\mathcal{J}_{K^{n+1}}|\,\mathbf{u}_m(t_{n+1}) - |\mathcal{J}_{K^n}|\,\mathbf{u}_m(t_n) \\
&\qquad = |\mathcal{J}_{\Delta t}| \sum_{g=1}^{n_g} \left[ w_g \sum_{c=1}^{n_c} w_c\, \mathcal{G}(\mathbf{x}_c^g, t^g, \mathcal{U}_c^g, \mathbf{w}_c) \cdot \left[\nabla_{\boldsymbol{\xi}_c}\phi_m(\boldsymbol{\xi}_c) \cdot \widehat{\mathcal{J}}_{K^g}^{-1}\right]\right] \\
&\quad - |\mathcal{J}_{\Delta t}| \sum_{g=1}^{n_g} \left[ w_g \sum_{F^g} |\mathcal{J}_{F^g}| \left[\sum_{f=1}^{n_f} w_f\, \mathcal{H}(\mathbf{x}^g, t^g, \mathcal{U}_f^{-,g}, \mathcal{U}_f^{+,g}; \mathbf{w}_f, \hat{\boldsymbol{\nu}})\phi_m(\boldsymbol{\sigma})\right]\right]
\end{aligned}
\tag{3.6.7}
$$

We can divide by $|\mathcal{J}_{K^{n+1}}|$ throughout and get the final form of the scheme

$$\mathbf{u}_m(t_{n+1}) = \frac{|\mathcal{J}_{K^n}|}{|\mathcal{J}_{K^{n+1}}|}\mathbf{u}_m(t_n) + \frac{|\mathcal{J}_{\Delta t}|}{|\mathcal{J}_{K^{n+1}}|}\sum_{g=1}^{n_g}\left[w_g\sum_{c=1}^{n_c} w_c\,\mathcal{G}(\mathbf{x}_c^g, t^g, \mathcal{U}_c^g, \mathbf{w}_c)\cdot\left[\nabla_{\boldsymbol{\xi}_c}\phi_m(\boldsymbol{\xi}_c)\cdot\widehat{\mathcal{J}}_{K^g}^{-1}\right]\right]$$

$$-\frac{|\mathcal{J}_{\Delta t}|}{|\mathcal{J}_{K^{n+1}}|}\sum_{g=1}^{n_g}\left[w_g\sum_{F^g}|\mathcal{J}_{F^g}|\left[\sum_{f=1}^{n_f} w_f\,\mathcal{H}(\mathbf{x}^g, t^g, \mathcal{U}_f^{-,g}, \mathcal{U}_f^{+,g}; \mathbf{w}_f, \hat{\boldsymbol{\nu}})\phi_m(\boldsymbol{\sigma})\right]\right]. \quad (3.6.8)$$

We can finally write it in a detailed form as

$$\mathbf{u}_m(t_{n+1}) = \frac{|\mathcal{J}_{K^n}|}{|\mathcal{J}_{K^{n+1}}|}\mathbf{u}_m(t_n) + \sum_{g=1}^{n_g}\sum_{c=1}^{n_c}\frac{|\mathcal{J}_{\Delta t}|}{|\mathcal{J}_{K^{n+1}}|}\left[w_g\, w_c\,\mathcal{G}(\mathbf{x}_c^g, t^g, \mathcal{U}_c^g, \mathbf{w}_c)\cdot\left[\nabla_{\boldsymbol{\xi}_c}\phi_m(\boldsymbol{\xi}_c)\cdot\widehat{\mathcal{J}}_{K^g}^{-1}\right]\right]$$

$$-\sum_{g=1}^{n_g}\sum_{F^g}\sum_{f=1}^{n_f}\frac{|\mathcal{J}_{\Delta t}|\,|\mathcal{J}_{F^g}|}{|\mathcal{J}_{K^{n+1}}|}\left[w_g\, w_f\,\phi_m(\boldsymbol{\sigma})\mathcal{H}(\mathbf{x}^g, t^g, \mathcal{U}_f^{-,g}, \mathcal{U}_f^{+,g}; \mathbf{w}_f, \hat{\boldsymbol{\nu}})\right]. \quad (3.6.9)$$

## 3.7 Computing Vertex Velocities

The DG solution is discontinuous at the interfaces of cells and therefore at the vertices. This implies that there is no unique fluid velocity available from the scheme at the vertices. Instead, we must compute an approximation for the fluid velocity which is appropriate in the context.

The mesh velocity must be close to the local fluid velocity in order to have a Lagrangian character to the scheme. Some researchers, especially in the context of Lagrangian methods, solve a Riemann problem at the cell face to determine the face velocity. Since we use an ALE formulation, we do not require the exact fluid velocity which is anyway not available to use since we only have a predicted solution. Following the exact trajectory of the fluid would also lead to curved trajectories for the grid point, which is an unnecessary complication. In our work, we make two different choices for the mesh velocities:

1. **Average Mesh Velocities** The first choice is to take the average of the velocities at the barycentre of the neighboring cells.

$$\mathbf{w}_v = \frac{1}{N_v}\sum_{i=1}^{N_v}\mathbf{w}_{N_v} \quad (3.7.1)$$

2. **Linearized Riemann Problem** The second choice is to solve a linearized Riemann problem at the vertex $V$ at time $t_n$. For Euler equations in one dimension, the linearized Riemann solver gives us

$$\mathbf{w}_{j+\frac{1}{2}}^n = \frac{\rho_j^n c_j^n v_j^n + \rho_{j+1}^n c_{j+1}^n v_{j+1}^n}{\rho_j^n c_j^n + \rho_{j+1}^n c_{j+1}^n} + \frac{p_j^n - p_{j+1}^n}{\rho_j^n c_j^n + \rho_{j+1}^n c_{j+1}^n} \quad (3.7.2)$$

where $\mathbf{w}_{j+\frac{1}{2}}$ is the mesh velocity at $j+\frac{1}{2}$ face, $\rho$ is density, $c$ is the speed of sound, $v$ is the velocity and $p$ is the pressure.

### 3.7.1 Computation of Velocities at the Boundary

The above mentioned methods to compute mesh velocity algorithms are appropriate for the interior vertices. On the boundary vertices, however, we need to ensure that the mesh velocities accurately reflect the fluid motion without negatively affecting the mesh motion.

1. Open Boundary

   For open boundary, the situation is comparatively straight forward. We select all the neighbors in the average mesh velocity formula Eq. (3.7.1).

2. Reflective Boundary

   In case of reflective boundary, we compute the velocity by selecting all the neighbors in the average mesh velocity formula Eq. (3.7.1) and then reflect the velocity along the normal of the boundary. Let $\mathbf{w}_0$ be the velocity computed by the average formula and let $\boldsymbol{\nu}$ be the normal to the boundary at the vertex, then the resultant mesh velocity is given by

$$\mathbf{w} = \mathbf{w}_0 - 2(\mathbf{w} \cdot \boldsymbol{\nu})\hat{\boldsymbol{\nu}} \tag{3.7.3}$$

3. Periodic Boundary

   For periodic boundaries, the vertex on one boundary is a clone of the vertex on the periodic boundary. This requires that the velocity computed for the paired vertices should be equal. In order to achieve that, we generate the list of cells which are neighbors of either of the vertices and then use the average mesh velocity formula. Let $NK_{n_k}$ be the neighbors of the vertex $K_k$ and let $NK_{n_l}$ be the neighbors of the vertex $K_{k_p}$ where $K_{k_p}$ is the paired vertex. Let $NK = NK_{n_k} \cup NK_{n_l}$, then we have

$$\mathbf{w} = \frac{1}{|NV|} \sum_{i \in |NV|} \mathbf{w}_i \tag{3.7.4}$$

4. Closed Boundary

   At closed boundaries, the fluid can either go into or come out of the fluid boundaries. Ideally, this would require insertion or deletion of points in the mesh. In absence of such technology, we move the mesh along with the fluid just like in open boundary case.

## 3.8 Dissipation in ALE DG Scheme

This section has been adapted from [35]. In this section, we study the dissipation of the ALE DG Scheme for a scalar transport equation in one dimension

$$\frac{\partial u}{\partial t} + a\frac{\partial u}{\partial x} = 0 \tag{3.8.1}$$

We assume that the mesh velocity is a constant $w$ over the whole domain. The semi-discrete formulation of the ALE DG Scheme in this particular case can be written as

$$\frac{\mathrm{d}}{\mathrm{d}t}\int_{K_j(t)} u(x,t)\phi(x,t)\mathrm{d}x = \int_{K_j(t)} g(u,w)\frac{\mathrm{d}\phi(x,t)}{\mathrm{d}x}\mathrm{d}x - h_{j+\frac{1}{2}}\phi^-_{j+\frac{1}{2}} + h_{j-\frac{1}{2}}\phi^+_{j-\frac{1}{2}} \tag{3.8.2}$$

where,

$$g(u,w) = (a-w)u \tag{3.8.3}$$

and $h(u^-,u^+)$ is a consistent numerical flux. An alternative form of the above equation can be obtained by applying integration by parts to the first term on the right hand side, to get

$$\frac{\mathrm{d}}{\mathrm{d}t}\int_{K_j(t)} u(x,t)\phi(x,t)\mathrm{d}x = -\int_{K_j(t)} \frac{\mathrm{d}g(u,w)}{\mathrm{d}x}\phi(x,t)\mathrm{d}x + \phi^-_{j+\frac{1}{2}}\left(g^-_{j+\frac{1}{2}} - h_{j+\frac{1}{2}}\right) - \phi^+_{j-\frac{1}{2}}\left(g^+_{j-\frac{1}{2}} - h_{j-\frac{1}{2}}\right) \tag{3.8.4}$$

We use an upwind numerical flux

$$h_{j+\frac{1}{2}} = g^-_{j+\frac{1}{2}} \qquad \forall j \tag{3.8.5}$$

Substituting in the above equation,we get

$$\frac{\mathrm{d}}{\mathrm{d}t}\int_{K_j(t)} u(x,t)\phi(x,t)\mathrm{d}x = -\int_{K_j(t)} \frac{\mathrm{d}g(u,w)}{\mathrm{d}x}\phi(x,t)\mathrm{d}x - \phi^+_{j-\frac{1}{2}}\left(g^+_{j-\frac{1}{2}} - g^-_{j-\frac{1}{2}}\right) \tag{3.8.6}$$

Next, we substitute the expression for $g(u,w)$ in the equation

$$\frac{\mathrm{d}}{\mathrm{d}t}\int_{K_j(t)} u(x,t)\phi(x,t)\mathrm{d}x = -\int_{K_j(t)} (a-w)\frac{\mathrm{d}u}{\mathrm{d}x}\phi(x,t)\mathrm{d}x - (a-w)\phi^+_{j-\frac{1}{2}}\left(u^+_{j-\frac{1}{2}} - u^-_{j-\frac{1}{2}}\right) \tag{3.8.7}$$

Next, we write the solution in form of the local basis and use the DG formulation by specifying $\phi_m$ as the test function and assuming that the solution is written as

$$u_j(x,t) = \sum_{m=0}^{M} u_{j,m}(t)\phi_{j,m}(x,t) \tag{3.8.8}$$

$$\frac{\mathrm{d}u_j(x,t)}{\mathrm{d}x} = \sum_{m=0}^{M} u_{j,m}(t)\frac{\mathrm{d}\phi_{j,m}(x,t)}{\mathrm{d}x} \tag{3.8.9}$$

we can write

$$\frac{\mathrm{d}}{\mathrm{d}t}\int_{K_j(t)} u(x,t)\phi_{j,l}(x,t)\mathrm{d}x = -\int_{K_j(t)} (a-w)\frac{\mathrm{d}u}{\mathrm{d}x}\phi_{j,l}(x,t)\mathrm{d}x - (a-w)\phi^+_{j-\frac{1}{2},l}\left(u^+_{j-\frac{1}{2}} - u^-_{j-\frac{1}{2}}\right)$$

$$\begin{aligned}\frac{\mathrm{d}u_{j,l}}{\mathrm{d}t} &= -\int\limits_{K_j(t)} (a-w)\left(\sum_{m=0}^{M} u_{j,m}(t)\frac{\mathrm{d}\phi_{j,m}(x,t)}{\mathrm{d}x}\right)\phi_{j,l}(x,t)\mathrm{d}x - (a-w)\phi^{+}_{j-\frac{1}{2},l}\left(u^{+}_{j-\frac{1}{2}} - u^{-}_{j-\frac{1}{2}}\right)\\ &= -(a-w)\sum_{m=0}^{M} u_{j,m}(t)\int\limits_{K_j(t)}\left(\frac{\mathrm{d}\phi_{j,m}(x,t)}{\mathrm{d}x}\right)\phi_{j,l}(x,t)\mathrm{d}x - (a-w)\phi^{+}_{j-\frac{1}{2},l}\left(u^{+}_{j-\frac{1}{2}} - u^{-}_{j-\frac{1}{2}}\right)\end{aligned}$$

Define

$$d_{r,s} = \int\limits_{K_j(t)}\left(\frac{\mathrm{d}\phi_{j,r}(x,t)}{\mathrm{d}x}\right)\phi_{j,s}(x,t) \qquad p_{r,s} = \phi_{j-1,r}\phi_{j,s}\big|_{x_{j-\frac{1}{2}}} \qquad q_{r,s} = \phi_{j,r}\phi_{j,s}\big|_{x_{j-\frac{1}{2}}} \tag{3.8.10}$$

Note that due to the our choice of basis functions, $d_{r,s}, p_{r,s}, q_{r,s}$ are independent of $t, j$.

$$\frac{\mathrm{d}u_{j,l}}{\mathrm{d}t} = -(a-w)\sum_{m=0}^{M}\left[\left[d_{m,l} + q_{m,l}\right]u_{j,m} + p_{m,l}u_{j-1,m}\right] \tag{3.8.11}$$

Define

$$D = (d_{r,s}) \qquad P = (p_{r,s}) \qquad Q = (q_{r,s}) \tag{3.8.12}$$

$$W = (D+Q) \qquad V = P \tag{3.8.13}$$

Next, assume that we have $N$ cells with periodic boundary conditions. Define

$$U = \left[u_{j,l}\right] \qquad A = \begin{bmatrix} U & 0 & \cdots & 0 & L\\ L & U & & & \\ & \ddots & \ddots & & \\ & & \ddots & \ddots & \\ & & & L & U\end{bmatrix} \tag{3.8.14}$$

whence, we get the matrix equation

$$\frac{\mathrm{d}U}{\mathrm{d}t} = -(a-w)AU \tag{3.8.15}$$

This means that the dissipation is directly proportional to $|a-w|$.

## 3.9 Dependence of DG Error on Mesh Equality

With evolution of mesh, the mesh can become distorted, and the errors in the simulation can rise with degrading mesh quality. It is then important to have a quantitative measure to determine the quality of the mesh. There are some straightforward measures which have been discussed in the path, for example, obtuse triangles are bad for errors. Most of the analysis found in literature deals with finite volume methods with [30] having an extensive set of mesh quality indicators. In this section, we propose a mesh quality indicator motivated by the interpolation errors of Sobolev spaces on the polynomial spaces over simplexes. In this report, we focus on two dimensional case.

Let $\mathit{\Pi} \in L(W^{k+1,p}(K), W^{m,q}(K))$, then we have due to [11]

$$\|\mathbf{u} - \mathit{\Pi}_k \mathbf{u}\|_{W^{m,q}(K,\mathbb{V})} \leq C \frac{R_K^{k+1}}{r_K^m} |K|^{\frac{1}{q}-\frac{1}{p}} \|\mathbf{u}\|_{W^{k+1,p}(K,\mathbb{V})} \tag{3.9.1}$$

where $r_K$ is the inradius, and $R_K$ is the circumradius of the triangle $K$. We assume that the solution exists in a Sobolev space $H^{k+1}$, and let $\mathit{\Pi}_k \in L(H^{k+1}(K), \mathbb{P}^m(K))$, then we substitute $p = q = 2$ in the above equation and get

$$\|\mathbf{u} - \mathit{\Pi}_k \mathbf{u}\|_{H^m(K,\mathbb{V})} \leq C \frac{R_K^{k+1}}{r_K^m} \|\mathbf{u}\|_{H^{k+1}(K,\mathbb{V})} \tag{3.9.2}$$

$$\|\mathbf{u} - \mathit{\Pi}_k \mathbf{u}\|_{H^1(K,\mathbb{V})} \leq C \frac{R_K^{k+1}}{r_k} \|\mathbf{u}\|_{H^{k+1}(K,\mathbb{V})} \tag{3.9.3}$$

$$\|\nabla_{\mathbf{x}} \mathbf{u} - \mathit{\Pi}_k \nabla_{\mathbf{x}} \mathbf{u}\|_{L^2(K,\mathbb{V})} \leq C \frac{R_K^{k+1}}{r_k} \|\mathbf{u}\|_{H^{k+1}(K,\mathbb{V})} \tag{3.9.4}$$

$$\|\mathbf{u} - \mathit{\Pi}_k \mathbf{u}\|_{L^2(K,\mathbb{V})} \leq C R_K^{k+1} \|\mathbf{u}\|_{H^{k+1}(K,\mathbb{V})} \tag{3.9.5}$$

It is helpful to replace $R_K$ using the right-most term given below

$$C_K \stackrel{\text{def}}{=} \frac{s^2}{|K|} \qquad s = \sum_i \frac{s_i}{2} \qquad s_i \text{ are side lengths} \tag{3.9.6}$$

$$\frac{1}{r_K} = \frac{s}{|K|} = C\frac{C_K}{s} \qquad R_K = \frac{\prod_i s_i}{4\,|K|} \stackrel{\text{AM-GM}}{\leq} C\frac{s^3}{|K|} = C C_K\, s \tag{3.9.7}$$

Doing so, we get

$$\|\mathbf{u} - \mathit{\Pi}_k \mathbf{u}\|_{H^1(K,\mathbb{V})} \leq C C_K^{k+2} s^k \|\mathbf{u}\|_{H^{k+1}(K,\mathbb{V})} \tag{3.9.8}$$

$$\|\nabla_{\mathbf{x}} \mathbf{u} - \mathit{\Pi}_k \nabla_{\mathbf{x}} \mathbf{u}\|_{L^2(K,\mathbb{V})} \leq C C_K^{k+2} s^k \|\mathbf{u}\|_{H^{k+1}(K,\mathbb{V})} \tag{3.9.9}$$

$$\|\mathbf{u} - \mathit{\Pi}_k \mathbf{u}\|_{L^2(K,\mathbb{V})} \leq C_K^{k+1} s^{k+1} \|\mathbf{u}\|_{H^{k+1}(K,\mathbb{V})} \tag{3.9.10}$$

For DG schemes, the accuracy of gradients is as important as accuracy of solutions. In the above inequalities, the gradient inequality has a stricter conditions. So, motivated by the above inequalities, we propose $C_k$ as a mesh quality parameter.

**Definition 3.9.1** ($Q_K$ : RMS Quality of a Simplex)**.** *The quality of a simplex element $K$ with*

*faces $F_i$ is defined as*

$$Q_K = \frac{S_{rms}^2}{3\sqrt{3}\,|K|} - 1 \tag{3.9.11}$$

| Triangle | Mesh Quality |
|---|---|
| 60-60-60 | 0 |
| 45-45-90 | |

Table 3.9.1: Mesh Quality and Interpolation Error

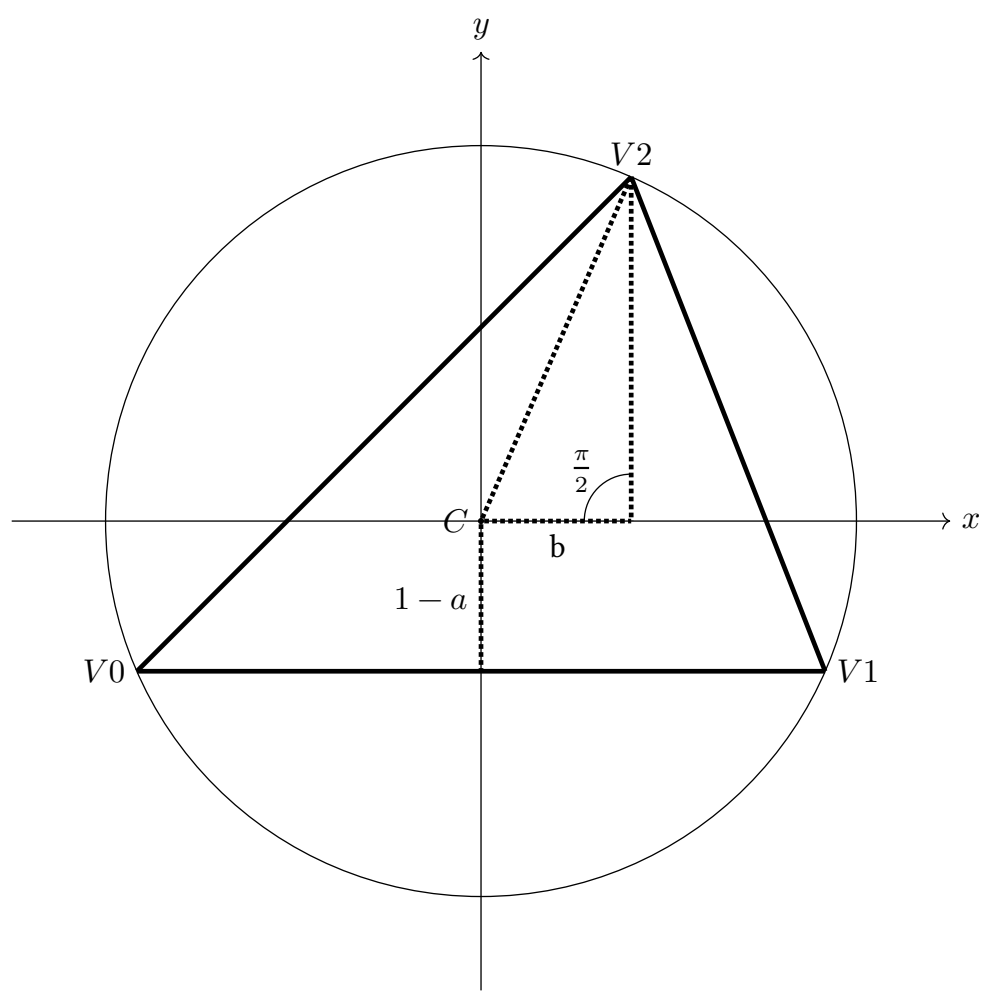

Figure 3.9.1: An Arbitrary Triangle

In the figures below, we evaluate the performance mesh quality indicator in comparison with actual simulations. In the comparison, we can see that the indicator is much stricter than the empirical results, which is a good thing.

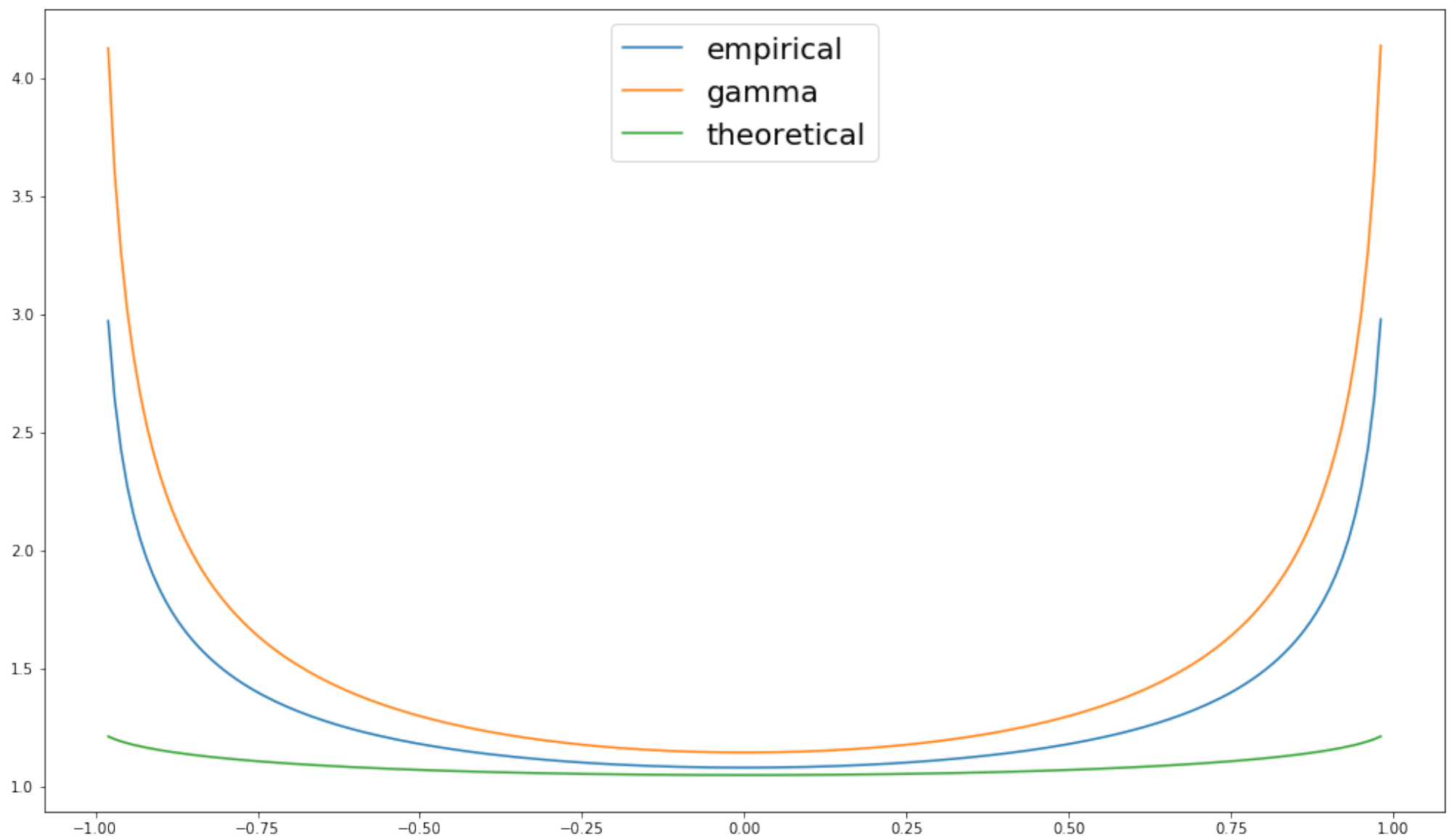
empirical
gamma
theoretical
4.0
3.5
3.0
2.5
2.0
1.5
1.0
−1.00
−0.75
−0.50
−0.25
0.00
0.25
0.50
0.75
1.00

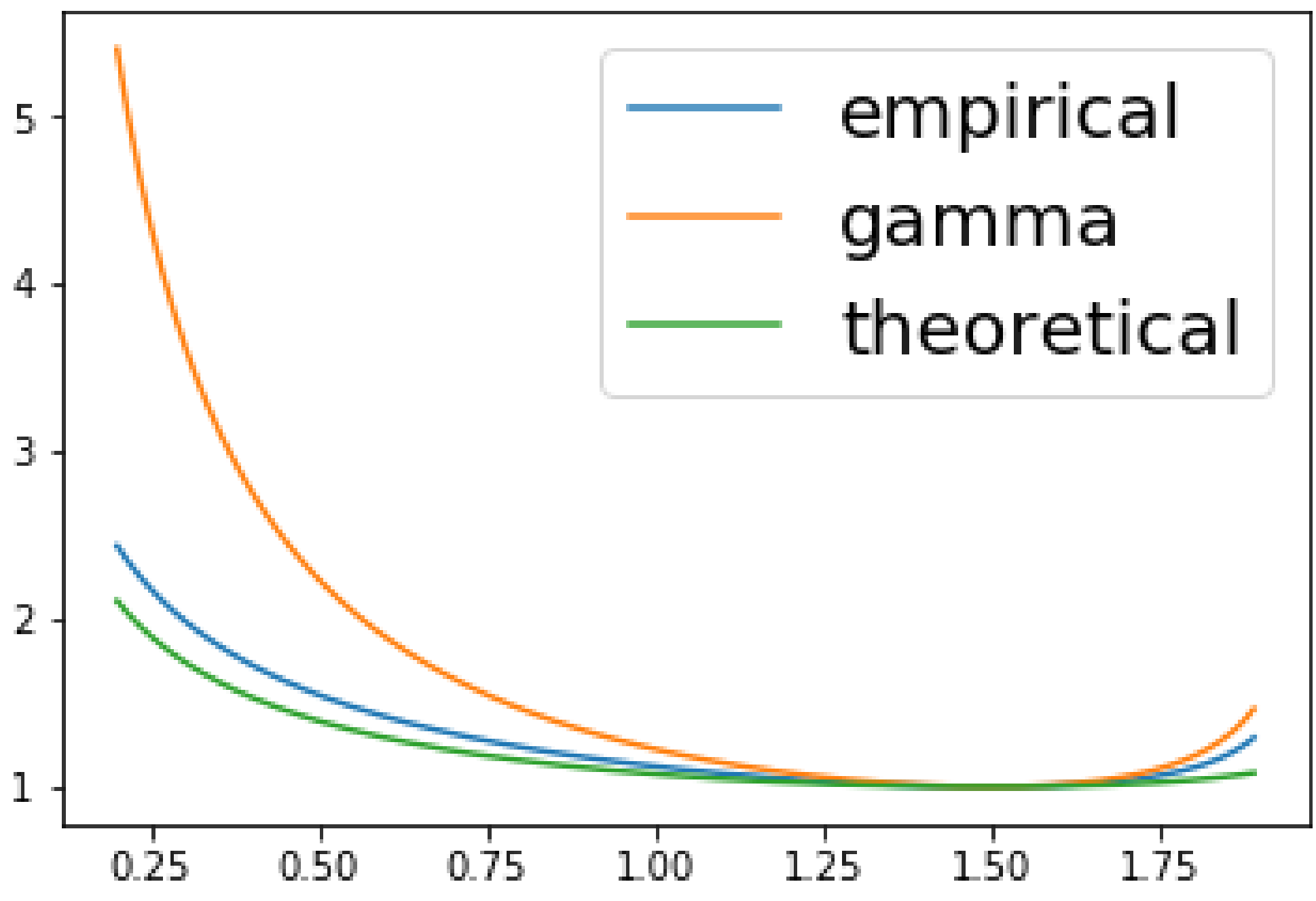
empirical
gamma
theoretical
5
4
3
2
1
0.25
0.50
0.75
1.00
1.25
1.50
1.75

## 3.10 Total Variation Bounded Limiter in Two Dimensions

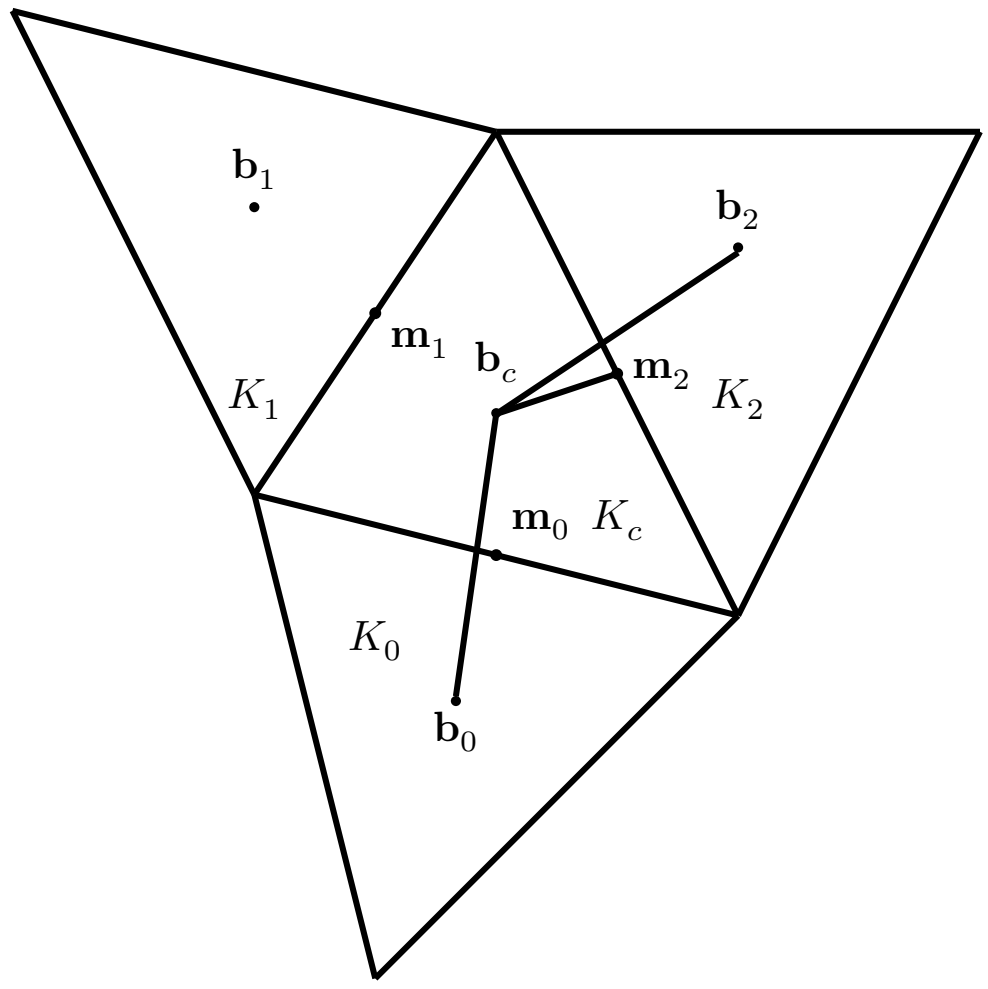

Figure 3.10.1: TVB Limiter in Two Dimensions

The TVB limiter presented here was first presented in [12]. In this discussion, we assume that $u_i$ is a piecewise linear scalar function over $K_i$. We write this function in terms of function $\psi(\mathbf{x})$ as

$$u_{K_0}(\mathbf{x}, t) = \sum_{i=1}^{3} u(\mathbf{m}_i)\psi_i(\mathbf{x}) \tag{3.10.1}$$

The function $\psi_i(\mathbf{x})$ is a linear function defined as

$$\psi_i(\mathbf{x}) = \begin{cases} 1 & \text{if } \mathbf{x} = \mathbf{m}_i \\ 0 & \text{if } \mathbf{x} = \mathbf{m}_j, j \neq i \end{cases} \tag{3.10.2}$$

The above conditions uniquely determine $\psi_i(\mathbf{x})$. The averages over the cells are given by

$$\overline{u}_i = \frac{1}{|K_i|} \int_{K_i} u \,\mathrm{d}x = u(\mathbf{b}_i) \tag{3.10.3}$$

Define $\tilde{u}(\mathbf{m}_i, K_c)$ as

$$\tilde{u}(\mathbf{m}_i, K_c) \stackrel{\text{def}}{=} u(\mathbf{m}_i) - \overline{u}_c \tag{3.10.4}$$

Let $\mathbf{b}_c$ be the barycentre of triangle $K_c$. Let $\mathbf{b}_i$, $i = 0, 1, 2$ be the barycentre of triangles $K_i$, $i = 0, 1, 2$. Let $\mathbf{m}_i$ be the midpoint of the face $F_i$. For each $\mathbf{m}_i$, define $\mathbf{b}_{i,k}$ and $\mathbf{b}_{i,l}$ such $\mathbf{m}_i$ lies in the directed regular angle between $\mathbf{b}_{i,k}$ and $\mathbf{b}_{i,l}$, that is the directed angles have the relation $\measuredangle \mathbf{m}_i\mathbf{b}_c\mathbf{b}_{i,l} \leq \measuredangle \mathbf{b}_{i,k}\mathbf{b}_c\mathbf{b}_{i,l} \leq \pi$. Then there exists $\alpha_{i,j} \geq 0$, $i = 0, 1, 2$, $j = 0, 1$ such

that

$$\mathbf{m}_0 - \mathbf{b}_c = \alpha_{0,0}(\mathbf{b}_{0,k} - \mathbf{b}_c) + \alpha_{0,1}(\mathbf{b}_{1,l} - \mathbf{b}_c) \tag{3.10.5a}$$

$$\mathbf{m}_1 - \mathbf{b}_c = \alpha_{1,0}(\mathbf{b}_{1,k} - \mathbf{b}_c) + \alpha_{1,1}(\mathbf{b}_{2,l} - \mathbf{b}_c) \tag{3.10.5b}$$

$$\mathbf{m}_2 - \mathbf{b}_c = \alpha_{2,0}(\mathbf{b}_{2,k} - \mathbf{b}_c) + \alpha_{2,1}(\mathbf{b}_{0,l} - \mathbf{b}_c) \tag{3.10.5c}$$

whence we define

$$\Delta\overline{u}(\mathbf{m}_0, K_c) \overset{\text{def}}{=} \alpha_{0,0}(\overline{u}_{K_0} - \overline{u}_{K_c}) + \alpha_{0,1}(\overline{u}_{K_1} - \overline{u}_{K_c}) \tag{3.10.6a}$$

$$\Delta\overline{u}(\mathbf{m}_1, K_c) \overset{\text{def}}{=} \alpha_{1,0}(\overline{u}_{K_1} - \overline{u}_{K_c}) + \alpha_{1,1}(\overline{u}_{K_2} - \overline{u}_{K_c}) \tag{3.10.6b}$$

$$\Delta\overline{u}(\mathbf{m}_2, K_c) \overset{\text{def}}{=} \alpha_{2,0}(\overline{u}_{K_2} - \overline{u}_{K_c}) + \alpha_{2,1}(\overline{u}_{K_0} - \overline{u}_{K_c}) \tag{3.10.6c}$$

We now compute the minmod slope

$$\Delta_i \overset{\text{def}}{=} \overline{m}\Big(\tilde{u}(\mathbf{m}_i, K_c), \nu\Delta\overline{u}(\mathbf{m}_i, K_c)\Big) \tag{3.10.7}$$

where

$$\overline{m}(a_1, a_2, \dots, a_m) \overset{\text{def}}{=} \begin{cases} a_1, & \text{if } |a_1| \le M\Delta x^2 \\ m(a_1, a_2, \dots, a_m), & \text{otherwise} \end{cases} \tag{3.10.8}$$

$$m(a_1, a_2, \dots, a_m) \overset{\text{def}}{=} \begin{cases} s \min_i |a_i| & \text{if } s = \operatorname{sgn}(a_1) = \dots = \operatorname{sgn}(a_m) \\ 0 & \text{otherwise} \end{cases} \tag{3.10.9}$$

Then, we can write the limiting function $\Lambda\Pi$ for the following two different cases as shown

Case 1: $\sum_{i=0}^{2} \Delta_i = 0$ :

$$\Lambda\Pi u(\mathbf{x}) \overset{\text{def}}{=} \overline{u}_{K_c} + \sum_{i=0}^{2} \Delta_i \psi_i(\mathbf{x}) \tag{3.10.10a}$$

Case 2: $\sum_{i=0}^{2} \Delta_i \neq 0$ :

$$\text{pos} \overset{\text{def}}{=} \sum_{i=0}^{2} \max(0, \Delta_i) \qquad \text{neg} \overset{\text{def}}{=} \sum_{i=0}^{2} \max(0, -\Delta_i) \tag{3.10.10b}$$

and set

$$\theta^+ \overset{\text{def}}{=} \min\left(1, \frac{\text{neg}}{\text{pos}}\right) \qquad \theta^- \overset{\text{def}}{=} \min\left(1, \frac{\text{pos}}{\text{neg}}\right) \tag{3.10.10c}$$

$$\hat{\Delta}_i = \theta^+ \max(0, \Delta_i) - \theta^- \max(0, -\Delta_i) \tag{3.10.10d}$$

Then we define

$$\Lambda\Pi u(\mathbf{x}) = \overline{u}_{K_c} + \sum_{i=0}^{2} \hat{\Delta}_i \psi_i(\mathbf{x}) \tag{3.10.10e}$$

### 3.10.1 TVB Limiter for Systems in Two Dimensions

If $u$ is a solution to a hyperbolic system of conservation laws then, instead of limiting $\mathbf{u}$ directly, we limit in the characteristic variables. Define the Jacobians $\mathcal{J}_i$ as

$$\mathcal{J}_i = \frac{\partial}{\partial \mathbf{u}} \mathcal{F}\left(\mathbf{u}_{K_0}\right) \cdot \frac{\mathbf{m}_i - \mathbf{b}_0}{|\mathbf{m}_i - \mathbf{b}_0|} \tag{3.10.11}$$

Find matrices $\mathcal{R}_i$ and $\mathcal{R}_i^{-1}$ such that

$$\mathcal{R}_i^{-1} \mathcal{J}_i \mathcal{R}_i = \Lambda \tag{3.10.12}$$

Transform the quantities $\Delta\overline{u}(\mathbf{m}_i, K_0)$ in Eq. (3.10.6) by left-multiplying by $\mathcal{R}_i^{-1}$. Apply limiter component by component on the transformed quantities. Transform back the result by right-multiplying by $\mathcal{R}_i$.

### 3.10.2 Basis Functions for TVB Limiter

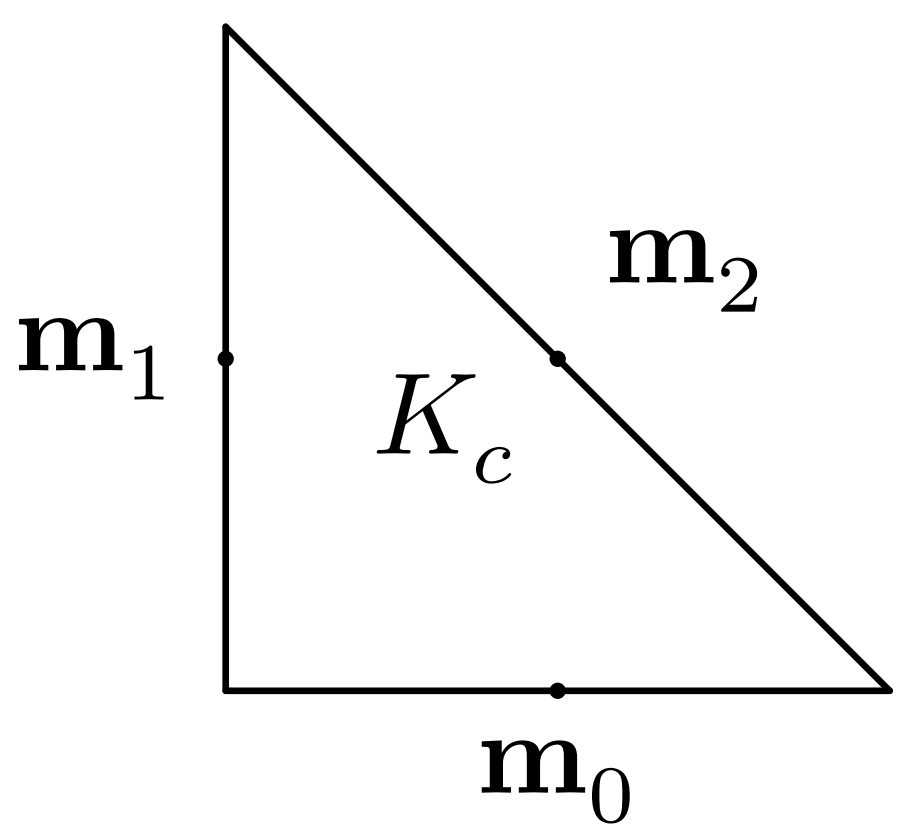

Figure 3.10.2: Basis Functions for TVB Limiter in Two Dimensions

As can be seen by Eq. (3.10.10), the projected solution is obtained in terms of $\psi_i(\mathbf{x})$. For the purpose of the TVB limiter, it is necessary to write the solution in terms of special functions. Let $\mathbf{u}(x, t)$ be the solution on the cell $K_c$. Then the solution needs to be represented in the form

$$u_{K_0}(\mathbf{x}, t) = \sum_{i=1}^{3} u(\mathbf{m}_i) \psi_i(\mathbf{x}) \tag{3.10.13}$$

where the function $\psi_i(\mathbf{x})$ is a linear function defined as

$$\psi_i(\mathbf{x}) = \begin{cases} 1 & \text{if } \mathbf{x} = \mathbf{m}_i \\ 0 & \text{if } \mathbf{x} = \mathbf{m}_j, j \neq i \end{cases} \tag{3.10.14}$$

For a linear function, this condition uniquely determines the function. Assume that the function is of the form

$$\psi_i(\mathbf{m}_1) = \mathbf{P}^1(\mathbf{x}) \cdot \mathbf{a}_i = a_{i,0} + x a_{i,1} + y a_{i,2} \tag{3.10.15}$$

where

$$\mathbf{P}^1(\mathbf{x}) = \begin{bmatrix} 1 & x & y \end{bmatrix} \qquad \mathbf{a}_i = \begin{bmatrix} a_{i,0} \\ a_{i,1} \\ a_{i,2} \end{bmatrix} \tag{3.10.16}$$

Then, we have

$$\psi_i(\mathbf{m}_1) = a_{i,0} + x_1 a_{i,1} + y_1 a_{i,2} = 1 \tag{3.10.17a}$$
$$\psi_i(\mathbf{m}_2) = a_{i,0} + x_2 a_{i,1} + y_2 a_{i,2} = 0 \tag{3.10.17b}$$
$$\psi_i(\mathbf{m}_3) = a_{i,0} + x_3 a_{i,1} + y_3 a_{i,2} = 0 \tag{3.10.17c}$$

This gives us the linear equation

$$\begin{bmatrix} 1 & x_1 & y_1 \\ 1 & x_2 & y_2 \\ 1 & x_3 & y_3 \end{bmatrix} \begin{bmatrix} a_{i,0} \\ a_{i,1} \\ a_{i,2} \end{bmatrix} = \begin{bmatrix} \delta_{0,i} \\ \delta_{1,i} \\ \delta_{2,i} \end{bmatrix} \tag{3.10.18}$$

which gives us

$$\mathbf{a}_0 = \frac{1}{2|K|} \begin{bmatrix} x_3 y_2 - x_2 y_3 \\ -y_2 + y_3 \\ x_2 - x_3 \end{bmatrix} \quad \mathbf{a}_1 = \frac{1}{2|K|} \begin{bmatrix} -x_3 y_1 + x_1 y_3 \\ y_1 - y_3 \\ -x_1 + x_3 \end{bmatrix} \quad \mathbf{a}_2 = \frac{1}{2|K|} \begin{bmatrix} x_2 y_1 - x_1 y_2 \\ -y_1 + y_2 \\ x_1 - x_2 \end{bmatrix} \tag{3.10.19}$$

For the reference cell, this is gives us

$$\mathbf{a}_0 = \frac{1}{2|K|} \begin{bmatrix} 1 \\ 0 \\ -2 \end{bmatrix} \qquad \mathbf{a}_1 = \frac{1}{2|K|} \begin{bmatrix} 1 \\ -2 \\ 0 \end{bmatrix} \qquad \mathbf{a}_2 = \frac{1}{2|K|} \begin{bmatrix} -1 \\ 2 \\ 2 \end{bmatrix} \tag{3.10.20}$$

### 3.10.3 Algorithm for the TVB Limiter

**Step 1:** Compute the barycentre $\mathbf{b}_c$ of the current cell $K_c$.

**Step 2:** Compute the midpoints $\mathbf{m}_i$ of the faces $F_i$, $i = 0, 1, 2$.

This function is provided by the grid.

**Step 3:** Determine the neighboring triangles and mark ghost cells as such.

There are three cases one has to consider when determining neighboring triangles.

a) No Boundary

When the cell is not a boundary cell, then the computation of neighboring cells is straightforward and given by the grid.

b) Periodic Boundary

When the cell is a periodic boundary cell, then the computation of neighboring cells is straightforward and given by the grid.

c) Open Boundary

When the cell is a open boundary cell, then there is a possibility that up to two of the neighbors do not exist each corresponding to a face of the cell. In this case, the neighbor is considered to be a ghost cell which is a reflection of the cell across the corresponding face. In the figure below for example, the triangles $V_3 - V_0 - V_1$ is congruent to $V_3 - V_2 - V_1$.

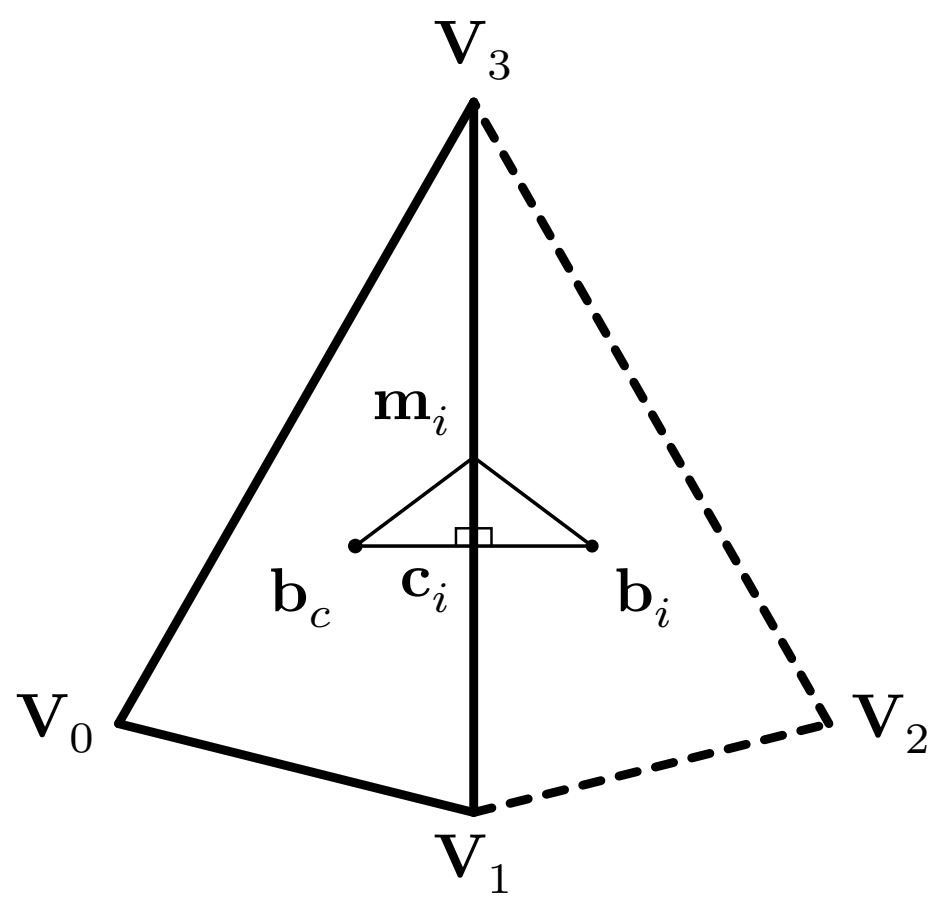

Figure 3.10.3: Ghost Cells for Open Boundary in Case of $i = 0$

**Step 4:** Compute the barycentre $\mathbf{b}_i$ of the neighboring cells $K_i$, $i = 0, 1, 2$.

a) For normal cells, the barycentre is provided by the grid.

b) For ghost cells on face $F_i$, we have the following

i. Compute the unit outward normal $\boldsymbol{\nu}_i$ on $F_i$ w.r.t $K_c$.

ii. Compute the $\boldsymbol{\nu}_i$ on $F_i$ w.r.t $K_c$.

iii. $\mathbf{b}_i = \mathbf{b}_c + 2[(\mathbf{m}_i - \mathbf{b}_c) \cdot \boldsymbol{\nu}_i]\boldsymbol{\nu}_i$

**Step 5:** Compute $\alpha_{i,j}$ which satisfies Eq. (3.10.5) as follows

a)

**Step 6:** Compute the solution $u$ at the barycentre of the triangles as follows

a) For normal cells, the method is provided by the finite element abstraction.

b) For ghost cells,

$$\overline{u}_{K_i} = \overline{u}_{K_c} \qquad i = 0, 1, 2 \tag{3.10.21}$$

**Step 7:** Compute $\Delta\overline{u}(\mathbf{m}_i, K_c)$ using the formula in Eq. (3.10.6)

**Step 8:** Compute the minmod slope $\Delta i$ as

$$\Delta_i \stackrel{\text{def}}{=} \overline{m}\Big(\tilde{u}(\mathbf{m}_i, K_c), \nu\Delta\overline{u}(\mathbf{m}_i, K_c)\Big) \tag{3.10.22}$$

where

$$\overline{m}(a_1, a_2, \dots, a_m) \stackrel{\text{def}}{=} \begin{cases} a_1, & \text{if } |a_1| \leq M\Delta x^2 \\ m(a_1, a_2, \dots, a_m), & \text{otherwise} \end{cases} \tag{3.10.23}$$

$$m(a_1, a_2, \dots, a_m) \stackrel{\text{def}}{=} \begin{cases} s\min_i |a_i| & \text{if } s = \operatorname{sgn}(a_1) = \dots = \operatorname{sgn}(a_m) \\ 0 & \text{otherwise} \end{cases} \tag{3.10.24}$$

**Step 9:** If $|\tilde{u}(\mathbf{m}_i, K_c) - \Delta_i| < \epsilon$, do nothing. Else, we can write the limiting function $\Lambda\Pi$ for the following two different cases as shown in Eq. (3.10.10).

**Step 10:** Next, we transform the solution in terms of basis function using inverse transform.

Chapter

# 4 Techniques to Prevent Mesh Degradation

In Chapter 3, we describe the numerical scheme which can be used to compute the solutions of conservation law on moving meshes. In the scheme, the mesh moves with velocities close to the velocity of the fluid. As we have seen in Section 3.9, the error of the ALE DG scheme depends on the mesh quality. And with the mesh moving with velocities close to the fluid velocity, it is inevitable that the mesh will be distorted near vortices and shocks and possibly other fluid features. This will lead to poor simulation results. One can see the distortion of the mesh when computing the isentropic vortex for the Euler equations using the ALE DG scheme. Apart from the distortion of the mesh, the mesh can also become degenerate, at which point, the simulation will stop completely.

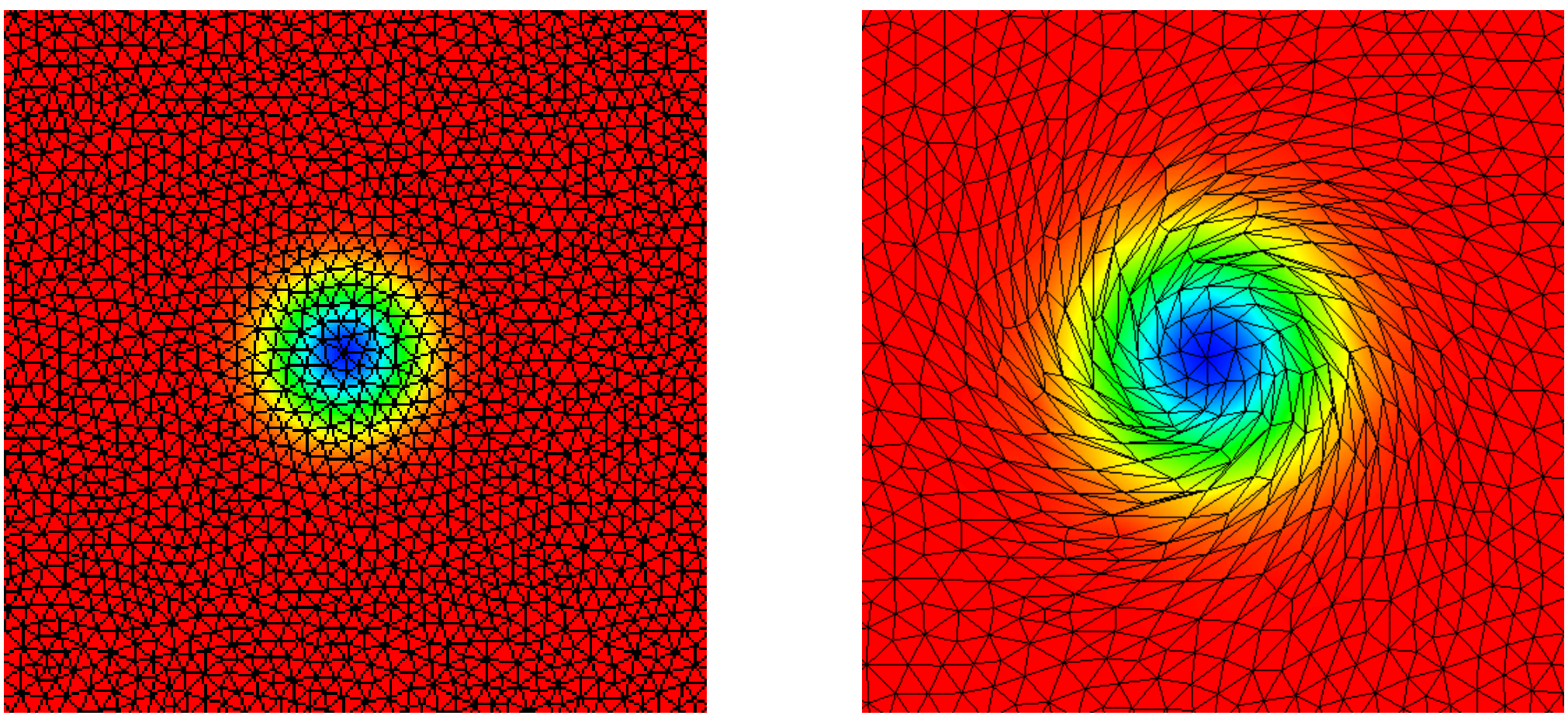

(a) $t = 0$ (b) $t = 20$

Figure 4.0.1: Degradation of Mesh Over Time $t = 20$

Therefore, it is important to ensure that the quality of the mesh stays above acceptable levels throughout the fluid simulation. In this chapter, we explore an array of techniques which allow us to achieve the goal. We prefer local techniques over global techniques in order to ensure that the techniques are suitable for simulation on high performance architectures.

## 4.1 Smoothing Algorithms for Mesh Velocity

One of the basic techniques to deal with mesh distortion is smoothing of mesh velocities. The technique prevents mesh quality from degrading severely in a single time step and hence, allowing mesh adaptation techniques to repair the mesh before the mesh quality degrades unacceptably. There are two different kinds of mesh smoothing algorithm which we explore in this report

1. Laplacian Smoothing Algorithm
2. Variable Diffusivity Laplacian Smoothing Algorithm

### 4.1.1 Laplacian Smoothing Algorithm

The objective of Laplacian smoothing is to minimize the distance between the position of a vertex and the position of the centroid of the polygon formed by neighboring vertices.

Let $\{\mathbf{w}_i^n\}_{i=1}^N$ be the velocity of the vertices of the mesh at time $t_n$ and $\mathbf{x}_i^{n+1}$ be the resultant position of the vertex $i$ obtained from. Let $\mathbf{X}_i^{n+1}$ be the average position of the vertexes connected to the $i$-th vertex. We would like the $\mathbf{x}_i^{n+1}$ to be close to $\mathbf{X}_i^{n+1}$. Hence, we can choose the mesh velocity to be

$$\mathbf{w}_i^n = \alpha \mathbf{w}_i^n + (1-\alpha)\frac{\mathbf{X}_i^{n+1} - \mathbf{x}_i^{n+1}}{\Delta t} \qquad \alpha \in [0,1] \tag{4.1.1}$$

The $\alpha$ factor moves $\mathbf{x}_i^{n+1}$ closer to $\mathbf{X}_i^{n+1}$ while the number of steps decrease the distance. We can apply a few iterations of the above smoothing in order to get a better velocity. The algorithm is described in Algorithm 1.

**Algorithm 1:** Algorithm Laplacian Smoothing

**parameter:** $\alpha$, `nsmooth`
**forall vertex** $i$ **do**
  $\mathbf{w}_i^n = \mathbf{w}_i^n$
**end**
**for** $s = 1$ **to** `nsmooth` **do**
  **forall vertex** $i$ **do**
    $\mathbf{x}_i^{n+1} = \mathbf{x}_i^n + \mathbf{w}_i^n \Delta t$
  **end**
  **forall vertex** $i$ **do**
    $\mathbf{X}_i^{n+1} = \frac{1}{N} \sum_{\text{neighbor}} \mathbf{x}_i^{n+1}$
  **end**
  **forall vertex** $i$ **do**
    $\mathbf{w}_i^n = \alpha \mathbf{w}_i^n + (1-\alpha)\frac{1}{\Delta t}\left[\mathbf{X}_i^{n+1} - \mathbf{x}_i^n\right]$
  **end**
**end**

In practice, we find that Laplacian smoothing is quite expensive due to the need of $\alpha$ and `nsmooth` to be very high. Otherwise one observes that the vertices collapse near the center of the vortices. [28] observed that this is due to non-uniform contributions to the velocity from the neighboring vertices irrespective of the distance between the vertices. They proposed that the contribution should also depend on the distance of one vertex from another. In light of this observation, they proposed variable diffusivity Laplacian smoothing algorithm.

### 4.1.2 Variable Diffusivity Laplacian Smoothing Algorithm

In variable diffusivity Laplacian smoothing algorithm, we solve

$$-\nabla_{\mathbf{x}}\left[\epsilon(\delta)\nabla_{\mathbf{x}}\mathbf{w}_i^n\right] + \mathbf{w}_i^n = \mathbf{w}_{i,0}^n \tag{4.1.2}$$

in order to derive smoothened velocities, where $\mathbf{w}_{i,0}^n$ is the velocity derived from the fluid velocity. $\epsilon(\delta)$ is given by

$$\epsilon(r) = \epsilon_0 + (1-\epsilon_0)\max\left(0, \min\left(1, \frac{\delta-\delta_l}{\delta_u-\delta_l}\right)\right) \qquad 0<\delta_l<\delta_u<1 \tag{4.1.3a}$$

$$\delta = \frac{r}{\max l_{ij}} \tag{4.1.3b}$$

We use a finite-difference based method to solve the equation as follows. Integrating over the neighboring cells, we get

$$\int \nabla_{\mathbf{x}}\cdot\left[\epsilon(\mathbf{x})\nabla_{\mathbf{x}}\mathbf{w}\right] + \int \mathbf{w} = \int \mathbf{v} \tag{4.1.4}$$

Apply divergence theorem and we get

$$\int_{\partial K}\left[\epsilon(\mathbf{x})\nabla_{\mathbf{x}}\mathbf{w}\right]\boldsymbol{\nu} + \int_K \mathbf{w} = \int_K \mathbf{v} \tag{4.1.5}$$

Splitting over the different cells, we get

$$\sum\int_{\partial K}\left[\epsilon(\mathbf{x})\nabla_{\mathbf{x}}\mathbf{w}\right]\boldsymbol{\nu} + \int \mathbf{w} = \int \mathbf{v} \tag{4.1.6}$$

$$\sum\int_{\partial K}\left[\epsilon(\mathbf{x})A_F\right] + \int \mathbf{w} = \int \mathbf{v} \tag{4.1.7}$$

$$\nabla_x\mathbf{w}\cdot\nu = 0.5\left[\frac{\mathbf{w}_m-\mathbf{c}}{\|\mathbf{m}-\mathbf{c}\|}\widehat{\mathbf{m}-\mathbf{c}}\cdot\nu + \frac{\mathbf{w}_{c2}-\mathbf{m}}{\|\mathbf{c}_2-\mathbf{m}\|}\widehat{\mathbf{c}_2-\mathbf{m}}\cdot\nu\right] \tag{4.1.8}$$

$$\mathbf{w}_0 = \frac{1}{\sum|K|}\left[\int \mathbf{v} - \sum\int_{\partial K}\left[\epsilon(\mathbf{x})A_F\right]\right] \tag{4.1.9}$$

$$\mathbf{w}_0 = \frac{1}{\sum|K|}\left[\int \mathbf{v} - \sum|K|\left[\epsilon(\mathbf{x}_g)A_F\right]\right] \tag{4.1.10}$$

The variable diffusivity Laplacian smoothing is less prone to collapse the vertices in a vortex and is less computationally expensive. However, it is less robust over long time simulations and we use Laplacian smoothing as a fallback in cases where the mesh quality becomes lower than is ideal.

## 4.2 Edge Swapping Techniques

Smoothing of velocity can maintain the mesh quality, however, it can destroy the almost-Lagrangian nature of the mesh motion. In order to maintain almost-Lagrangian, it is necessary to locally change the mesh topology. One of the methods to do so is the edge swapping methods.

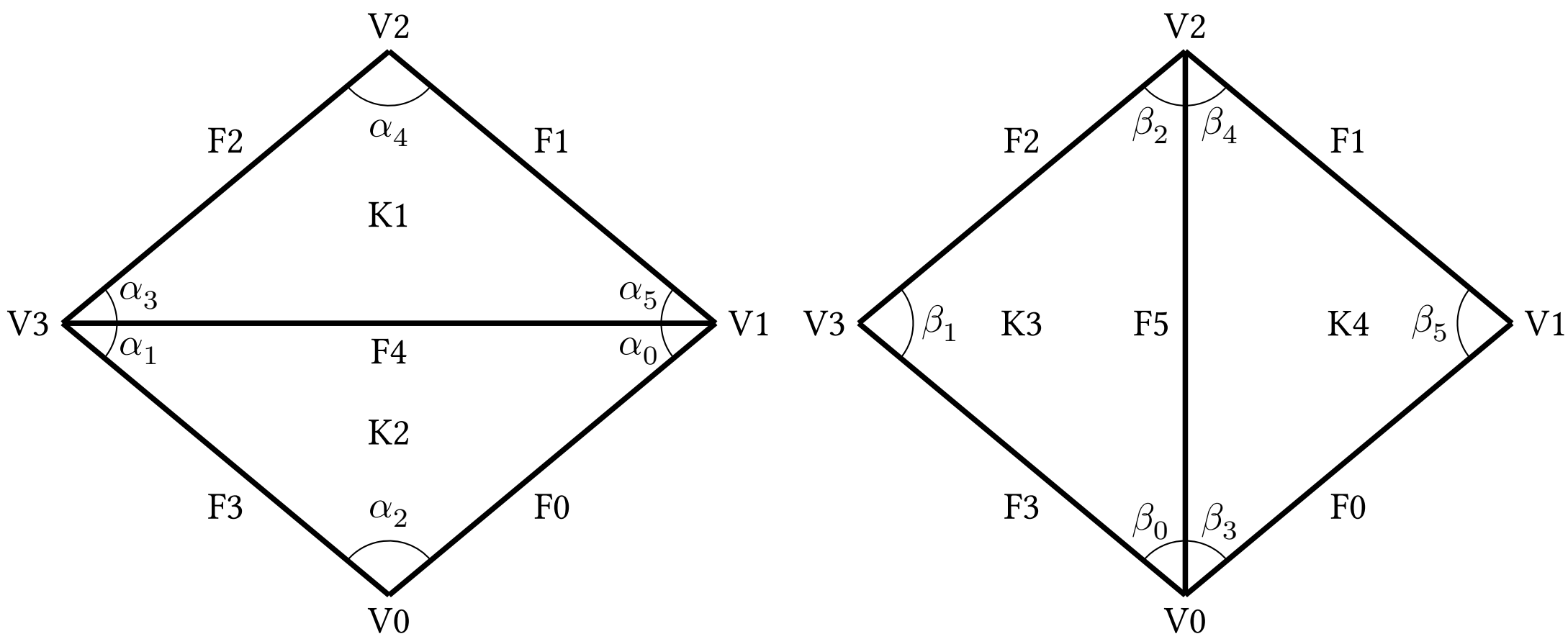

(a) Before (b) After

Figure 4.2.1: Face Swap

In an edge swap method, we change the connectivity of an edge in order to improve the quality of the mesh. As can be seen in Fig. 4.2.1, the mesh quality can improve in certain cases due to edge swaps. In order to figure out the situations in which the edge swaps are favorable, we use the mesh quality indicators introduced in **??**. However, direct usage the mesh quality indicator for swapping can often be suboptimal. This is because, often, a naive face swap algorithm is not commutative.

## 4.3 Solution Transfer for Arbitrary Remeshing

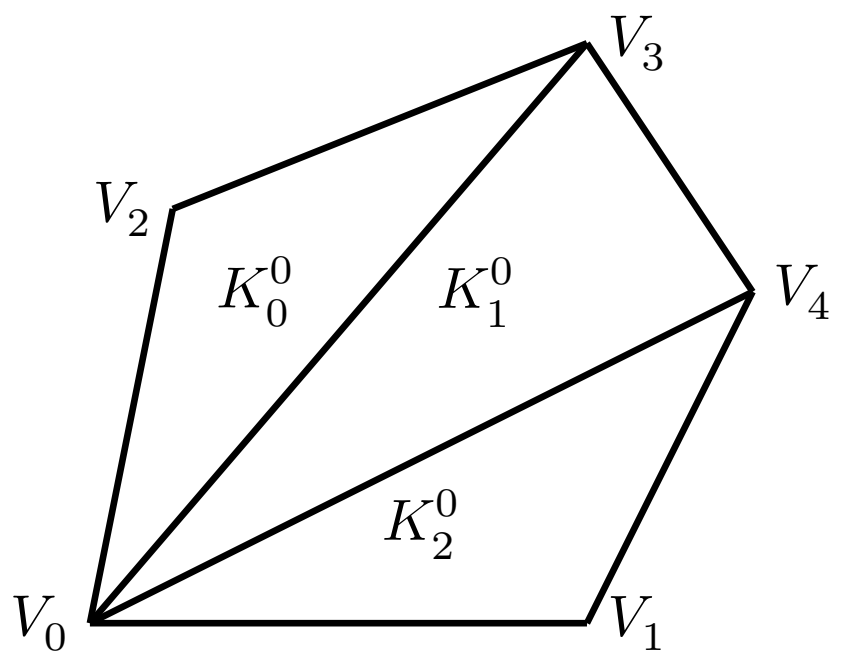

(a) Before

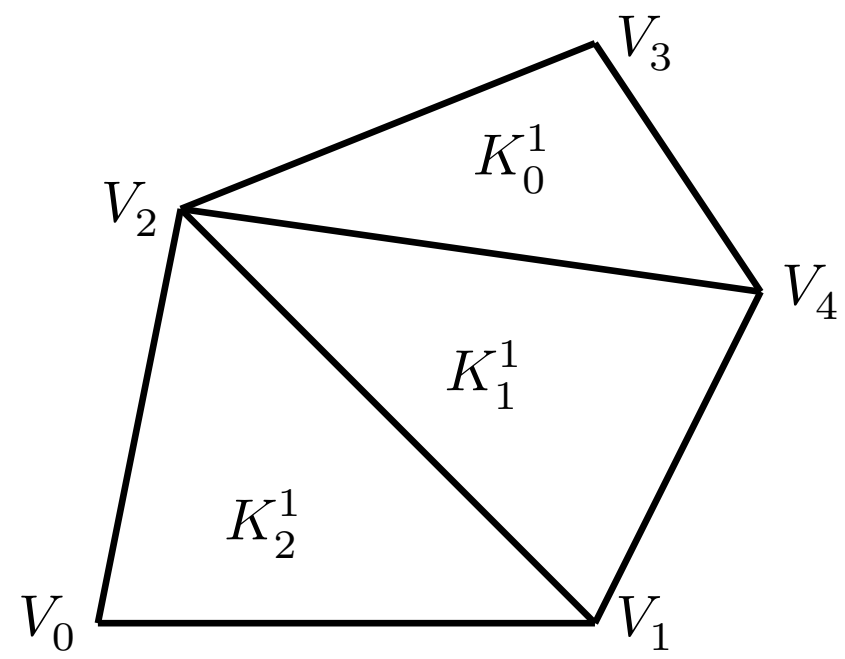

(b) After

Figure 4.3.1: Local Remesh

Consider that we have a section $\Omega_K$ of the mesh which is locally remeshed as shown in Fig. 4.3.1. Let $K_i^0$ be the list of triangles before remeshing and let $K_i^1$ be the list of triangles obtained after remeshing. An triangle $R_i \subset \Omega_K$ is an interpolation region if $R_i \subseteq K_j^0$ and $R_i \subseteq K_k^1$, for some $j, k$. Given a remeshing as shown above, one can decompose the domain into disjoint interpolation regions $\{R_i\}_{i=1}^M$ such that $R_i \cap R_j = \emptyset$ for $i \neq j$ and $\cup R_i = \Omega_K$. In general, one can always obtain such a region using a constrained Delaunay triangulation. However, for some special cases, like face swapping and face deletion, it is trivial to obtain such a region as will be shown later. For the case shown in Fig. 4.3.1, we show one such decomposition in Fig. 4.3.2.

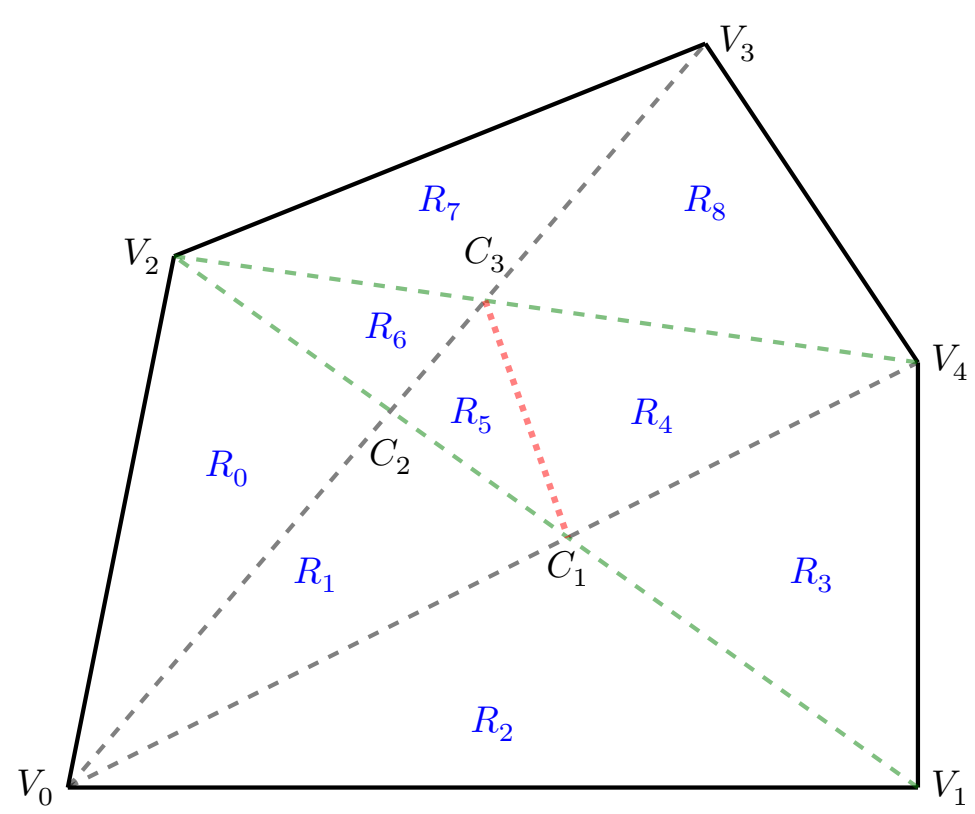

Figure 4.3.2: Interpolation Regions

Having obtained such a decomposition, the solution transfer can be done as follows. Let $\mathbf{u}$ be the solution over the domain $\Omega_K$. $\mathbf{u}_i^0$ be the solution in triangle $K_i^0$ for all $i$. We want to obtain the solution in triangle $\mathbf{u}_i^1$. We have

$$\mathbf{u}_{i,m}^1 = \int\limits_{K_i^1} \mathbf{u}\phi_{i,m}^1 \mathrm{d}\mathbf{x} \tag{4.3.1}$$

Now, we know that $K_i^1 = \cup R_j$ and in each of $R_j$, $\mathbf{u}$ is a polynomial, and hence we can write

$$\mathbf{u}_{i,m}^1 = \sum_{R_j} \int\limits_{R_j} \mathbf{u}\phi_{i,m}^1 \mathrm{d}\mathbf{x} \tag{4.3.2}$$

The integrals can now be evaluated exactly using quadrature rules. The remeshing algorithm now needs to provide only a description of how the regions are contained the triangles.

## 4.4 ALE DG Algorithm

Bringing all the techniques mentioned till now together, we have the following ALE DG algorithm:

1. Compute the fluid velocities at the vertices.
2. Compute the raw mesh velocities at the vertices using algorithms in section 3.7.1.
3. Compute the time step.
4. Compute the corrected mesh velocities using algorithms in 4.1
5. Compute the solution and the mesh at the time step $t_{n+1}$ according to the scheme, in the following steps.
   - Compute the solution at time step $t_{n+1}$.
   - Compute the position of vertices at the time step $t_{n+1}$.
   - Recompute the faces at the time step $t_{n+1}$.
   - Recompute the cells at time $t_{n+1}$ and compute the cell quality as well. Push all the cells with a cell quality below the threshold as described in REF onto a stack which is for named as "Poor Cell Stack".
6. Run the grid and solution adaptation algorithm as follows:
   a) Read a cell from poor cell stack.
   b) See if carrying out edge swap improves the quality of the mesh.
      - If yes, swap the cells, transfer the solution and push the neighbors onto the "Poor Cell Stack".
      - Otherwise, remove the cell of the poor cell stack.
   c) Continue till Poor Cell Stack is empty.
7. Go to step 1

Chapter

# 5 One Dimensional Applications

The numerical tests are performed with polynomials of degree one, two and three, together with the linear Taylor expansion, two stage CERK and four stage CERK, respectively, for the computation of the predictor. For the quadrature in time, we use the mid-point rule, two and three point Gauss-Legendre quadrature, respectively. The time step is chosen using the CFL condition,

$$\Delta t_n = \frac{\text{cfl}}{2k+1} \Delta t_n^{(1)}$$

where $\Delta t_n^{(1)}$ is given by equation (**??**), and the factor $(2k+1)$ comes from linear stability analysis [**cockburn1989**]. In most of the computations we use $\text{cfl} = 0.9$ unless stated otherwise. We observe that the results using average or linearized Riemann velocity are quite similar. We use the average velocity for most of the results and show the comparison between the two velocities for some results.

## 5.1 Order of Accuracy Test Case

We study the convergence rate of the schemes by applying them to a problem with a known smooth solution. The initial condition is taken as

$$\rho(x, 0) = 1 + \exp(-10x^2) \tag{5.1.1}$$

$$u(x, 0) = 1 \tag{5.1.2}$$

$$p(x, 0) = 1 \tag{5.1.3}$$

whose exact solution is $\rho(x,t) = \rho(x-t, 0)$, $u(x,t) = 1$, $p(x,t) = 1$. The initial domain is $[-5, +5]$ and the final time is $t = 1$ units. The results are presented using Rusanov and HLLC numerical fluxes. The $L^2$ norm of the error in density are shown in table (5.1.2), (5.1.3) for the static mesh and in table (5.1.4), (5.1.5) for the moving mesh. In each case, we see that the error behaves as $O(h^{k+1})$ which is the optimal rate we can expect for smooth solutions. In table (5.1.1), we show that the ALE DG methods preserves its higher order in presence of a limiter.

| $N$ | $k=1$ | | $k=2$ | |
|---|---|---|---|---|
| | Error | Rate | Error | Rate |
| 100 | 2.053E-02 | - | 2.277E-03 | - |
| 200 | 4.312E-03 | 2.251 | 3.425E-04 | 2.732 |
| 400 | 1.031E-03 | 2.064 | 4.565E-05 | 2.907 |
| 800 | 2.550E-04 | 2.015 | 5.812E-06 | 2.973 |
| 1600 | 6.356E-05 | 2.004 | 7.315E-07 | 2.990 |

Table 5.1.1: Order of accuracy study on moving mesh using Rusanov flux using Higher Order Limiter [**Zhong:2013:SWE:2397205.2397446**]

| $N$ | $k=1$ | | $k=2$ | | $k=3$ | |
|---|---|---|---|---|---|---|
| | Error | Rate | Error | Rate | Error | Rate |
| 100 | 4.370E-02 | - | 3.498E-03 | - | 3.883E-04 | - |
| 200 | 6.611E-03 | 2.725 | 4.766E-04 | 2.876 | 1.620E-05 | 4.583 |
| 400 | 1.332E-03 | 2.518 | 6.415E-05 | 2.885 | 9.376E-07 | 4.347 |
| 800 | 3.151E-04 | 2.372 | 8.246E-06 | 2.910 | 5.763E-08 | 4.239 |
| 1600 | 7.846E-05 | 2.280 | 1.031E-06 | 2.932 | 3.595E-09 | 4.180 |

Table 5.1.2: Order of accuracy study on static mesh using Rusanov flux

| $N$ | $k=1$ | | $k=2$ | | $k=3$ | |
|---|---|---|---|---|---|---|
| | Error | Rate | Error | Rate | Error | Rate |
| 100 | 4.582E-02 | | 3.952E-03 | | 3.464E-04 | |
| 200 | 9.611E-03 | 2.253 | 4.048E-04 | 3.287 | 2.058E-05 | 4.073 |
| 400 | 2.052E-03 | 2.240 | 4.640E-05 | 3.206 | 1.287E-06 | 4.036 |
| 800 | 4.803E-04 | 2.192 | 5.623E-06 | 3.152 | 8.061E-08 | 4.023 |
| 1600 | 1.184E-04 | 2.149 | 6.929E-07 | 3.119 | 5.050E-09 | 4.016 |

Table 5.1.3: Order of accuracy study on static mesh using HLLC flux

| $N$ | $k=1$ | | $k=2$ | | $k=3$ | |
|---|---|---|---|---|---|---|
| | Error | Rate | Error | Rate | Error | Rate |
| 100 | 2.331E-02 | - | 3.979E-03 | - | 8.633E-04 | - |
| 200 | 6.139E-03 | 1.9250 | 4.058E-04 | 3.294 | 1.185E-05 | 6.186 |
| 400 | 1.406E-03 | 2.0258 | 5.250E-05 | 3.122 | 7.079E-07 | 5.126 |
| 800 | 3.375E-04 | 2.0366 | 6.626E-06 | 3.077 | 4.340E-08 | 4.760 |
| 1600 | 8.278E-05 | 2.0344 | 8.304E-07 | 3.057 | 2.689E-09 | 4.573 |

Table 5.1.4: Order of accuracy study on moving mesh using Rusanov flux

| $N$ | $k=1$ | | $k=2$ | | $k=3$ | |
|---|---|---|---|---|---|---|
| | Error | Rate | Error | Rate | Error | Rate |
| 100 | 1.590E-02 | | 1.626E-03 | | 1.962E-04 | |
| 200 | 4.042E-03 | 1.977 | 2.072E-04 | 2.972 | 1.269E-05 | 3.950 |
| 400 | 1.014E-03 | 1.985 | 2.605E-05 | 2.982 | 7.983E-07 | 3.971 |
| 800 | 2.538E-04 | 1.990 | 3.261E-06 | 2.988 | 4.997E-08 | 3.980 |
| 1600 | 6.349E-05 | 1.992 | 4.077E-07 | 2.991 | 3.124E-09 | 3.985 |

Table 5.1.5: Order of accuracy study on moving mesh using HLLC flux

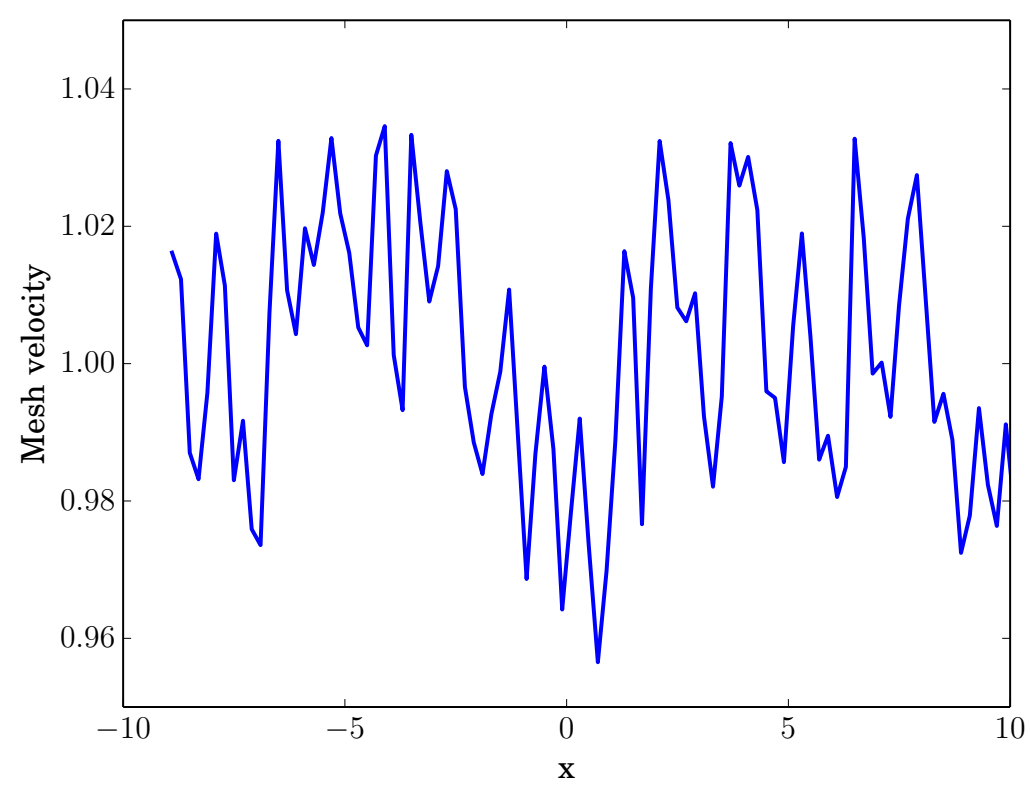

Figure 5.1.1: Example of randomized velocity distribution for smooth test case

| $N$ | $k=1$ | | $k=2$ | | $k=3$ | |
|---|---|---|---|---|---|---|
| | Error | Rate | Error | Rate | Error | Rate |
| 100 | 1.735E-02 | - | 1.798E-03 | - | 2.351E-04 | - |
| 200 | 4.179E-03 | 2.051 | 2.848E-04 | 2.676 | 1.416E-05 | 4.069 |
| 400 | 1.054E-03 | 2.035 | 4.301E-05 | 2.703 | 8.578E-07 | 4.041 |
| 800 | 2.615E-04 | 1.943 | 6.012E-06 | 2.838 | 5.476E-08 | 3.958 |
| 1600 | 7.279E-05 | 1.852 | 8.000E-07 | 2.909 | 3.505E-09 | 3.966 |

Table 5.1.6: Order of accuracy study on moving mesh using HLLC flux with randomly perturbed mesh velocity

| $N$ | $k=1$ | | $k=2$ | |
|---|---|---|---|---|
| | Error | Rate | Error | Rate |
| 100 | 8.535E-03 | - | 1.033E-03 | - |
| 200 | 1.958E-03 | 2.124 | 1.221E-04 | 3.08 |
| 400 | 4.721E-04 | 2.052 | 1.581E-05 | 2.95 |
| 800 | 1.238E-04 | 1.931 | 2.14E-06 | 2.89 |
| 1600 | 3.563E-05 | 1.796 | 2.63E-07 | 3.02 |

Table 5.2.1: Order of accuracy study on fixed mesh using Roe flux with Non-Constant Velocity Smooth Test Case

## 5.2 Smooth Test Case with Non-Constant Velocity

We also test the accuracy of our schemes on a isentropic problem with smooth solutions. The test case The initial conditions are given by

$$\rho(x,0) = 1 + 0.9999995\sin(\pi x) \qquad u(x,0) = 0 \qquad p(x,0) = \rho^{\gamma}(x,0) \tag{5.2.1}$$

with $\gamma = 3$ and periodic boundary conditions. For this kind of special isentropic problem, the Euler equations are equivalent to the two Burgers equations in terms of their two Riemann invariants which can then be used to derive the analytical solution. The errors are then computed with respect to the given analytical solution. In contrast to the previous test case, the velocity and pressure are not constant which makes this a more challenging test case. We run the simulation with a WENO-type limiter from [**Zhang:2010:PHO:1864819.1865066**] and positivity limiter enabled. As we can see from Tables 5.2.1, 5.2.2, the rate of convergence is maintained for the moving mesh method with the moving mesh methods exhibiting much lower errors.

| $N$ | $k=1$ | | $k=2$ | |
|---|---|---|---|---|
| | Error | Rate | Error | Rate |
| 100 | 4.235E-03 | - | 2.238E-04 | - |
| 200 | 1.058E-03 | 2.001 | 3.255E-05 | 2.87 |
| 400 | 2.586E-04 | 2.035 | 4.301E-05 | 3.133 |
| 800 | 5.804E-05 | 2.155 | 5.762E-06 | 2.901 |
| 1600 | 1.271E-05 | 2.192 | 7.401E-07 | 2.96 |

Table 5.2.2: Order of accuracy study on moving mesh using Roe flux with Non-Constant Velocity Smooth Test Case

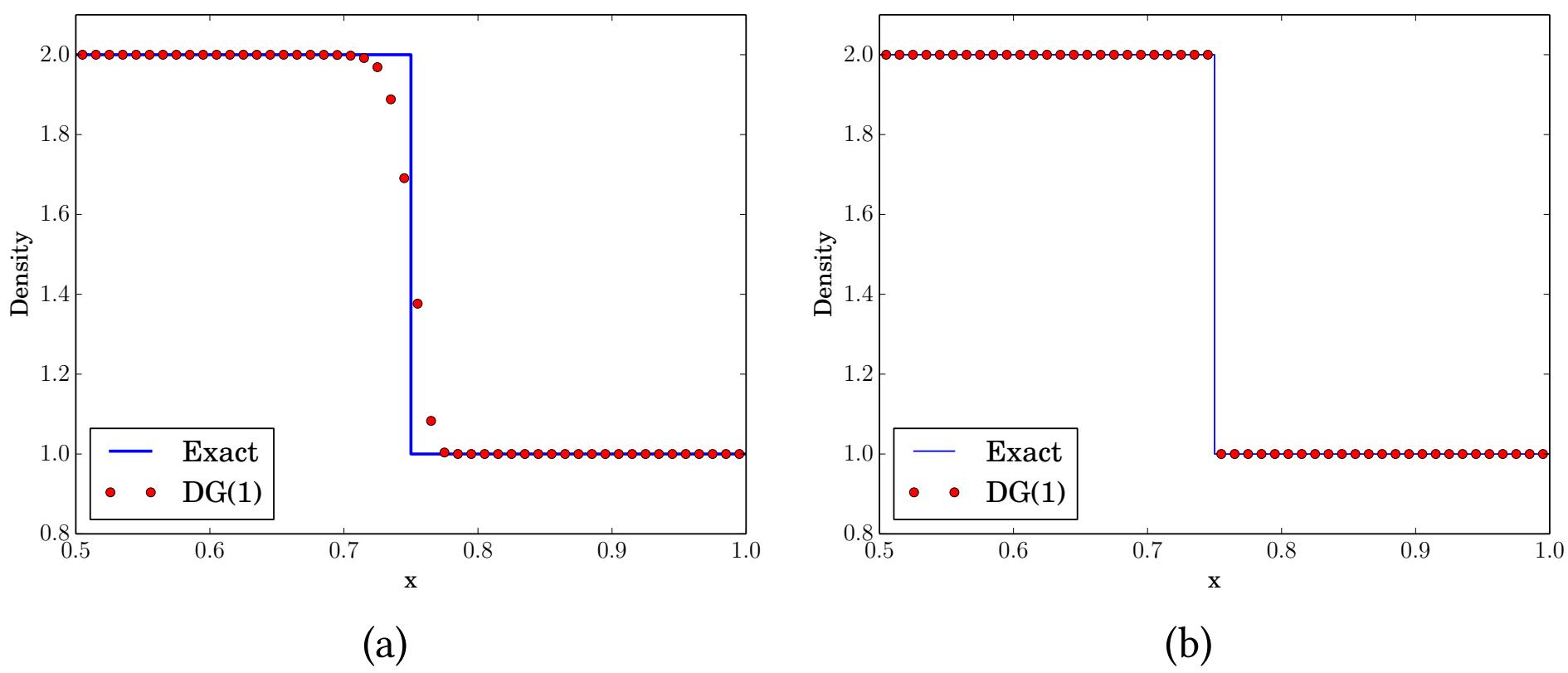

Figure 5.3.1: Single contact wave using Roe flux and 100 cells: (a) static mesh, (b) moving mesh

## 5.3 Single contact wave

In this example, we choose a Riemann problem which gives rise to a single contact wave in the solution that propagates with a constant speed. The initial condition is given by

$$(\rho, v, p) = \begin{cases} (2.0, 1.0, 1.0) & \text{if } x < 0.5 \\ (1.0, 1.0, 1.0) & \text{if } x > 0.5 \end{cases}$$

and the contact wave moves with a constant speed of 1.0. The solution on static and moving meshes are shown in figure (5.3.1) at time $t = 0.5$ using Roe flux. The moving mesh is able to exactly resolve the contact wave while the static mesh scheme adds considerable numerical dissipation that smears the discontinuity over many cells. The accurate resolution of contact waves is a key advantage of such moving mesh methods, which are capable of giving very good resolution of the contact discontinuity even on coarse meshes.

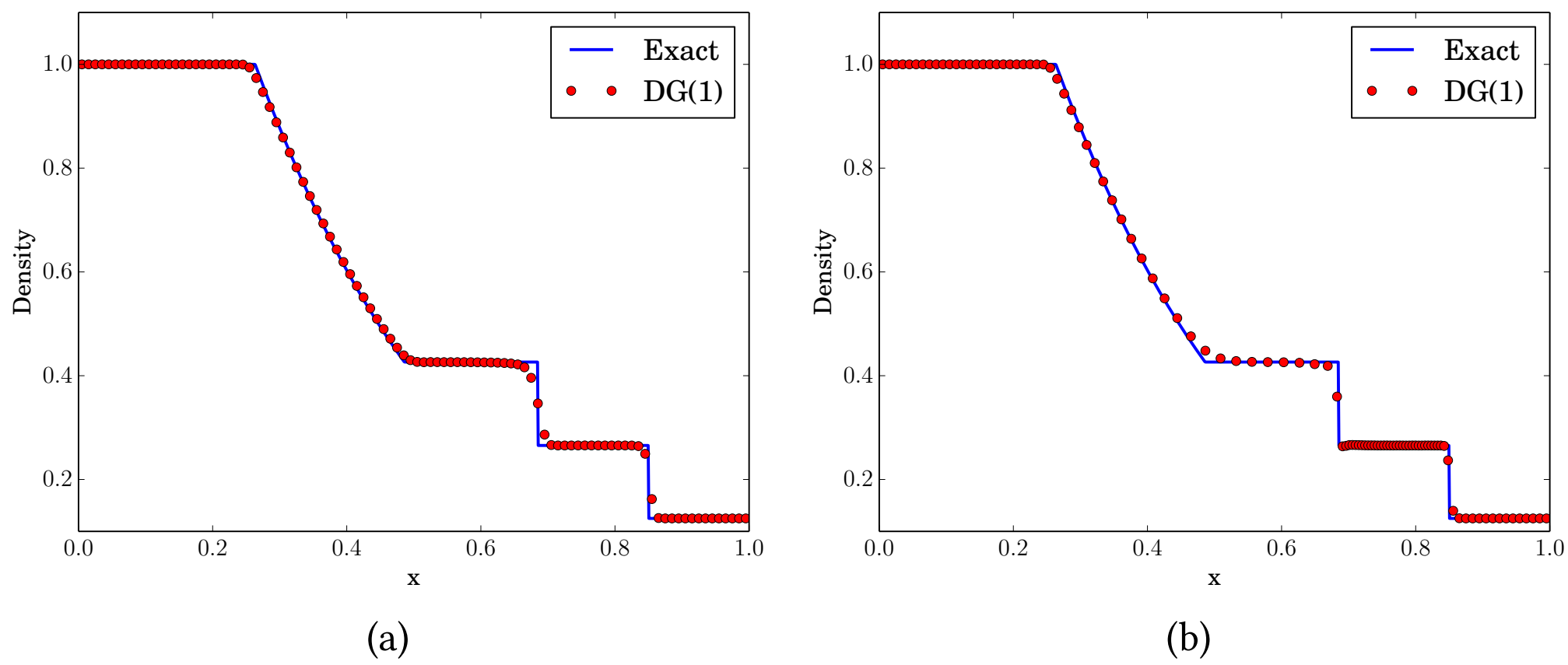

Figure 5.4.1: Sod problem using Roe flux, 100 cells and TVD limiter: (a) static mesh (b) moving mesh

## 5.4 Sod problem

The initial condition for the Sod test case is given by [**sod1978**]

$$(\rho, v, p) = \begin{cases} (1.0, 0.0, 1.0) & \text{if } x < 0.5 \\ (0.125, 0.0, 0.1) & \text{if } x > 0.5 \end{cases}$$

and the solution is computed upto a final time of $T = 0.2$ with the computational domain being $[0, 1]$. Since the fluid velocity is zero at the boundary, the computational domain does not change with time for the chosen final time. The exact solution consists of a rarefaction fan, a contact wave and a shock wave. In figure (5.4.1), we show the results obtained using Roe flux with 100 cells and TVD limiter on static and moving mesh. The contact wave is considerably well resolved on the moving mesh as compared to the static mesh due to reduced numerical dissipation on moving meshes.

To study the Galilean invariance or the dependence of the solution on the choice of coordinate frame, we add a boost velocity of $V = 10$ or $V = 100$ to the coordinate frame, while implies the initial fluid velocity is $v(x, 0) = V$ and the other quantities remain as before. Figure (5.4.4a) shows that the accuracy of the static mesh results degrades with increase in velocity of the coordinate frame, particularly the contact discontinuity is highly smeared. The results given in figure (5.4.4b) clearly show the independence of the results on the moving mesh with respect to the coordinate frame velocity. The allowed time step from CFL condition decreases with increase in coordinate frame speed for the static mesh case, while in case of moving mesh, it remains invariant. This means that in case of static mesh, we have to perform more time steps to reach the same final time as shown in table (5.4.1), which increases the computational time. Thus the moving mesh scheme has the additional advantage of allowing a larger time step compared to the fixed mesh scheme.

Finally, we compute the solutions using quadratic and cubic polynomials and the results are shown in figure (5.4.5). The solutions look similar to the case of linear polynomials and have the same sharp resolution of discontinuities.

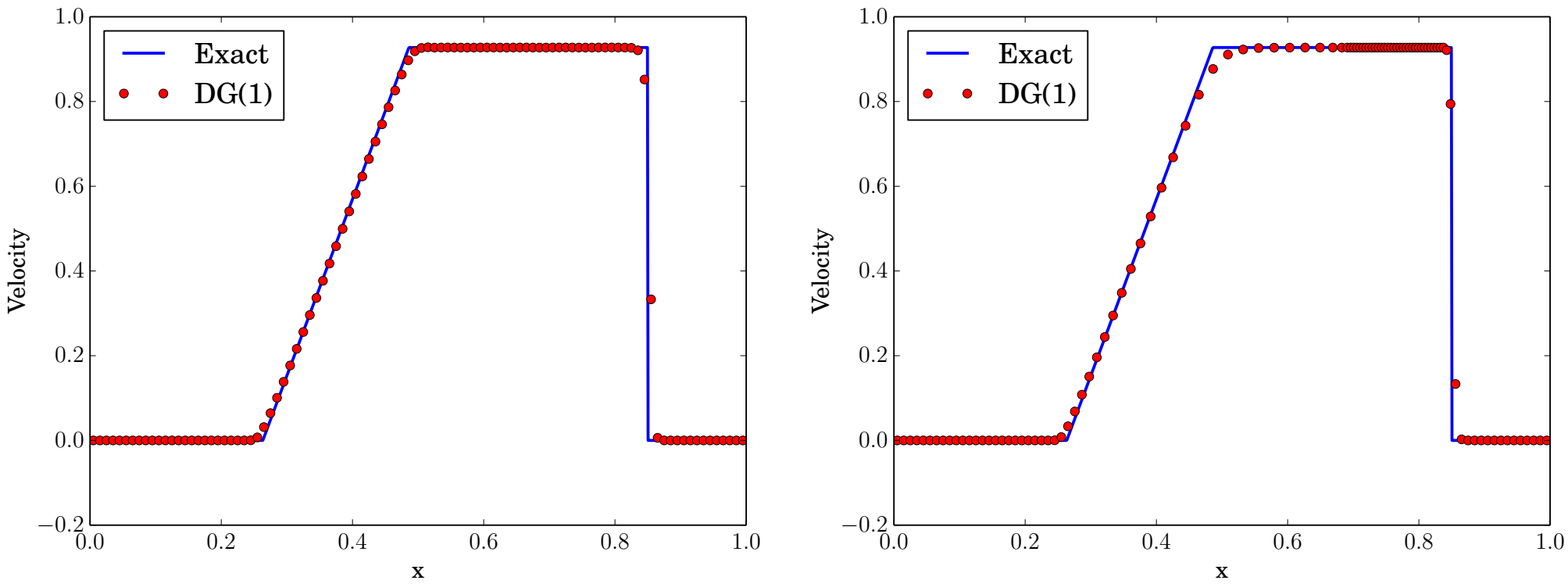

(a) (b)

Figure 5.4.2: Sod problem using Roe flux, 100 cells and TVD limiter: (a) static mesh (b) moving mesh

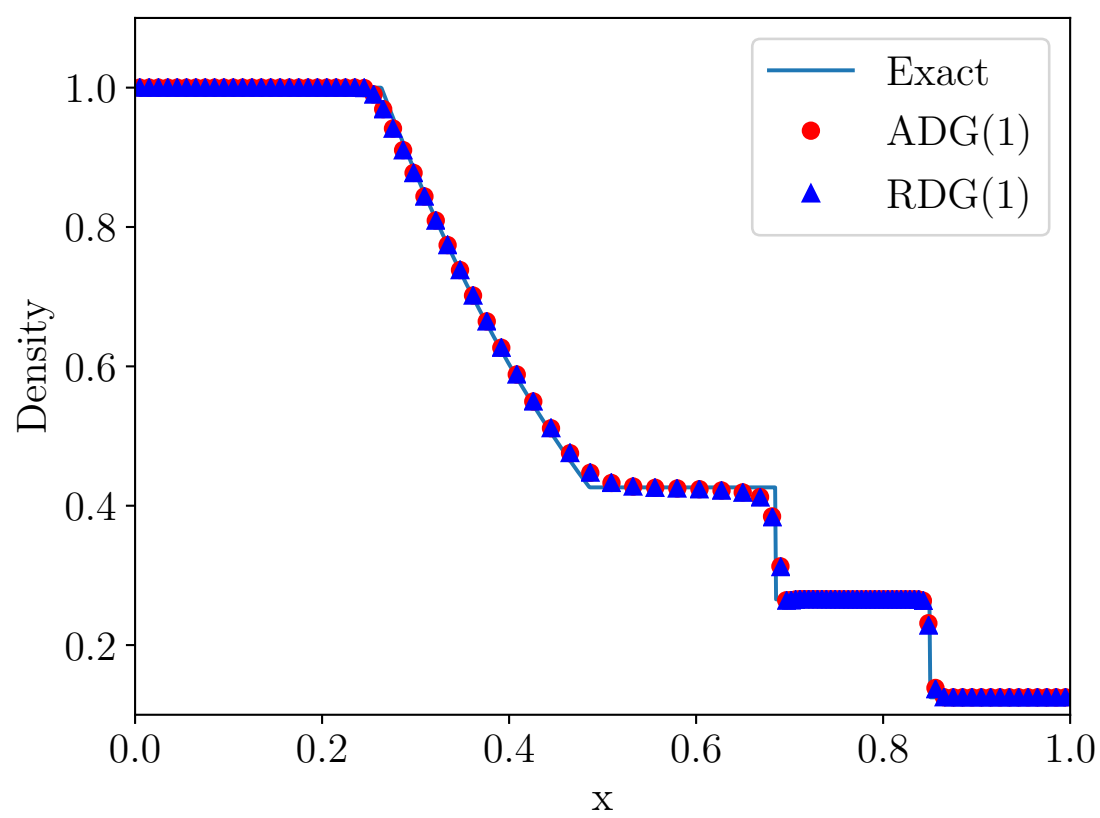

Figure 5.4.3: Sod problem using Roe flux, 100 cells and TVD limiter. ADG : Average Velocity, RDG : Linearized Riemann Velocity

| $V$ | 0 | 10 | 100 |
|---|---|---|---|
| static mesh | 144 | 810 | 6807 |
| moving mesh | 176 | 176 | 176 |

Table 5.4.1: Number of iterations required to reach time $t = 0.2$ for Sod test for different boost velocity of the coordinate frame

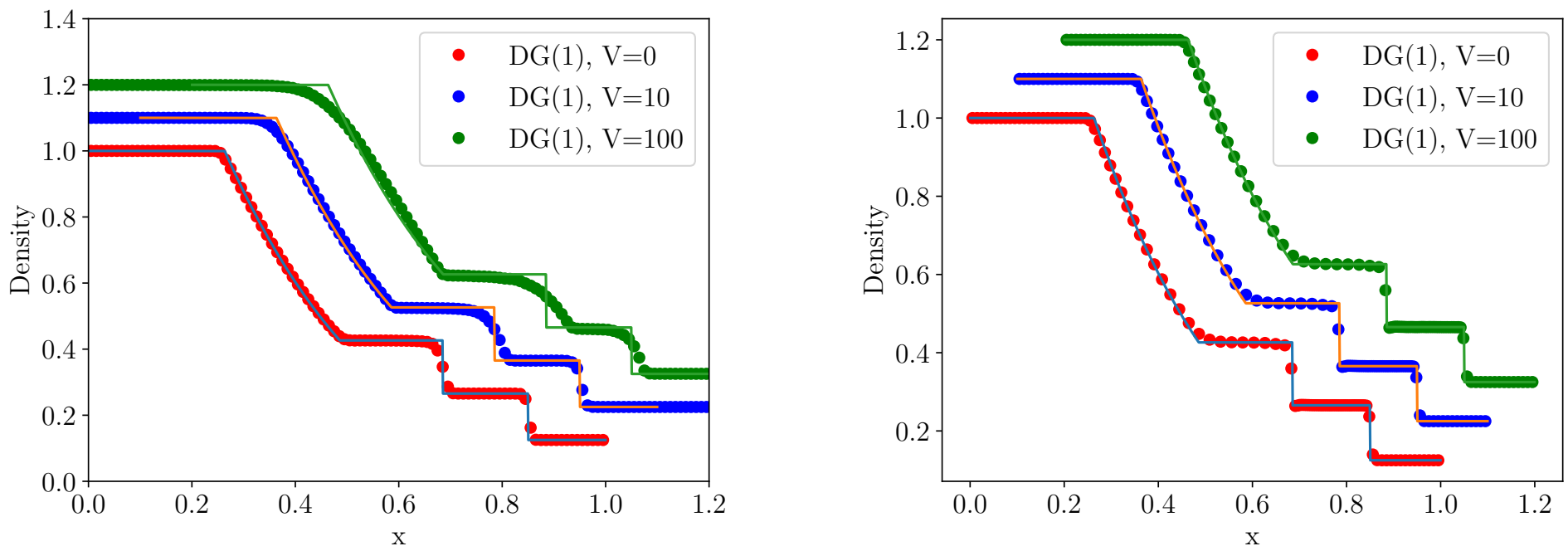

(a) (b)

Figure 5.4.4: Effect of coordinate frame motion on Sod problem using Roe flux, 100 cells and TVD limiter: (a) static mesh (b) moving mesh

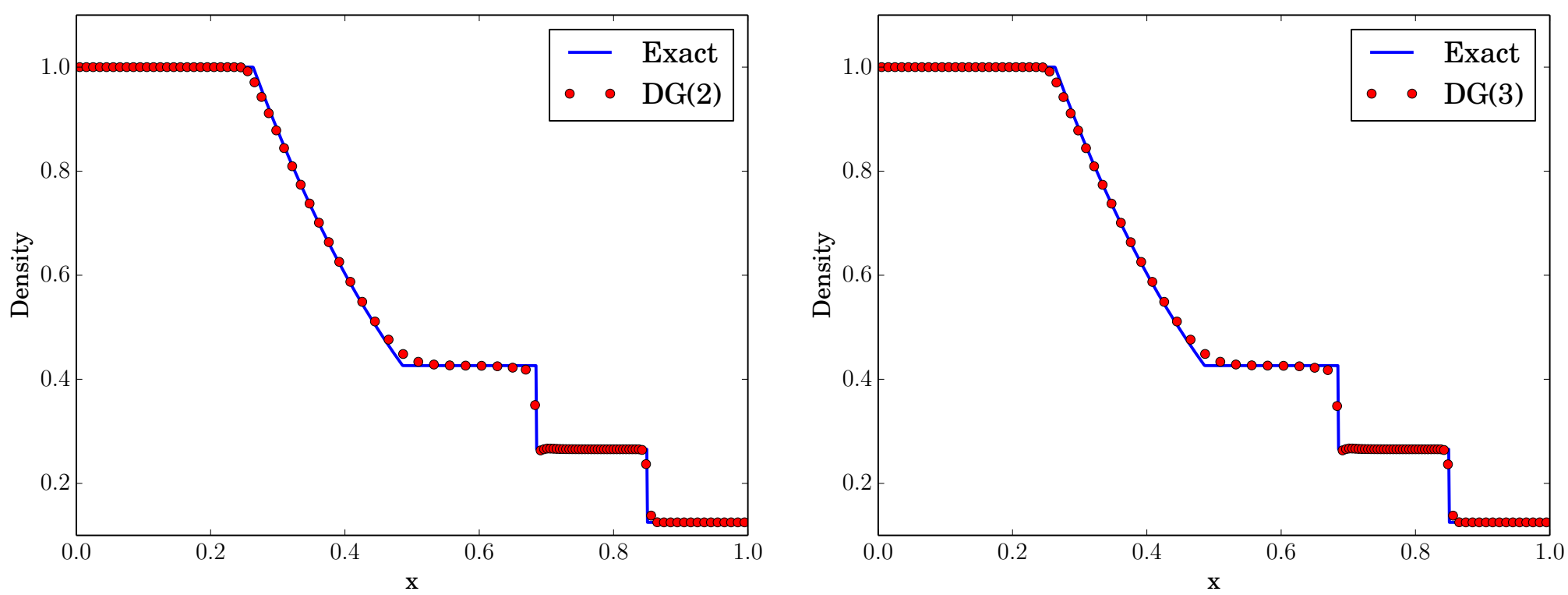

(a) (b)

Figure 5.4.5: Sod problem on moving mesh using Roe flux, 100 cells and TVD limiter: (a) Degree = 2 (b) Degree = 3

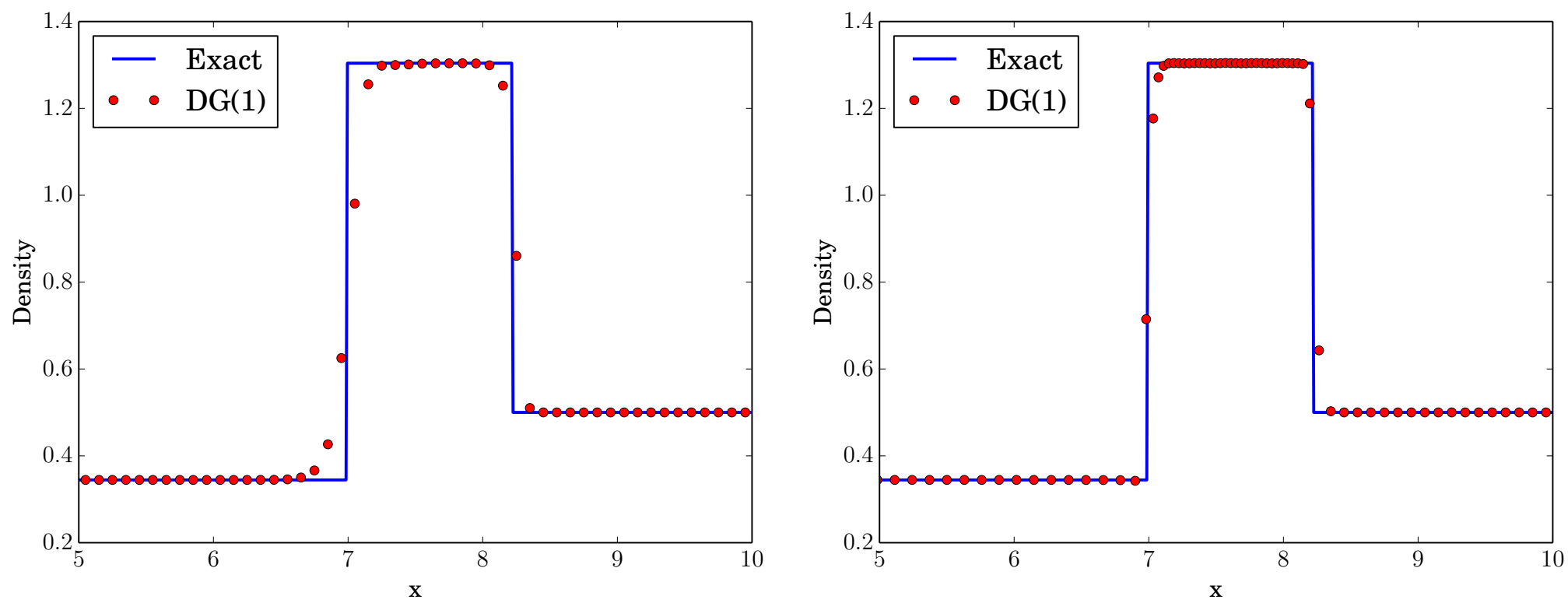

Figure 5.5.1: Lax problem using HLLC flux, 100 cells and TVD limiter: (a) static mesh (b) moving mesh

## 5.5 Lax problem

The initial condition is given by

$$(\rho, v, p) = \begin{cases} (0.445, 0.698, 3.528) & \text{if } x < 0 \\ (0.5, 0.0, 0.571) & \text{if } x > 0 \end{cases}$$

The computational domain is $[-10, +10]$ and we compute the solution up to a final time of $T = 1.3$. This problem has a strong shock and a contact wave that is difficult to resolve accurately. The zoomed view of density is shown at the final time in figure (5.5.1), and we observe the moving mesh results are more accurate for the contact wave, which is the first discontinuity in the figure. The second discontinuity is a shock which is equally well resolved in both cases. We can observe that the grid is automatically clustered in the region between the contact and shock wave, but no explicit grid adaptation was used in this simulation.

## 5.6 Shu-Osher problem

The initial condition is given by [**Shu1988439**]

$$(\rho, v, p) = \begin{cases} (3.857143, 2.629369, 10.333333) & \text{if } x < -4 \\ (1 + 0.2\sin(5x), 0.0, 1.0) & \text{if } x > -4 \end{cases}$$

which involves a smooth sinusoidal density wave which interacts with a shock. The domain is $[-5, +5]$ and the solution is computed up to a final time of $T = 1.8$. The solutions are shown in figure (5.6.1a)-(5.6.1b) on static and moving meshes using 200 cells and TVD limiter. The moving mesh scheme is considerably more accurate in resolving the sinusoidal wave structure that arises after interaction with the shock. In figure (5.6.1c) we compute the solution on static mesh with TVB limiter and the parameter $M = 100$ in equation (**??**). In this case the solutions on static mesh are more accurate compared to the case of TVD limiter but still not as good as the moving mesh results. The moving mesh result has more than 200 cells in the interval $[-5, +5]$ at the final time since cells enter the domain from the left side. Hence in figure (5.6.1d), we show the static mesh results with 300 cells and using TVB limiter. The results are further improved for the static mesh case but still not as accurate as the moving mesh case. The choice of parameters in the TVB limiter is very critical but we do not have a rigorous algorithm to choose a good value for this. Hence it is still advantageous to use the moving mesh scheme which gives improved solutions even with TVD limiter.

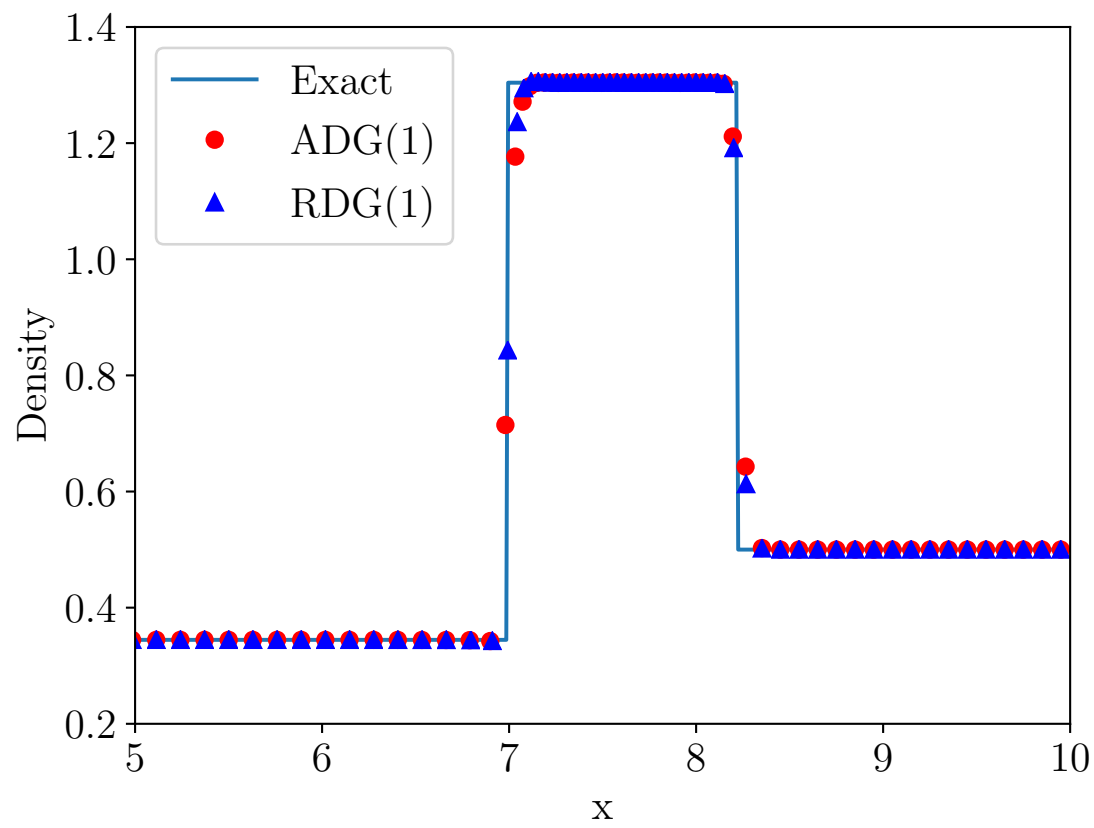

Figure 5.5.2: Lax problem using HLLC flux, 100 cells and TVD limiter. ADG : Average Velocity, RDG : Linearized Riemann Velocity

The above results show the ALE method is very accurate in terms of the cell averages. In figure (5.6.2), we show a zoomed view of density and pressure, where we also plot the linear polynomial solution. The slope of the solution is not accurately predicted with the Roe scheme and there are spurious contact discontinuities as the pressure and velocity are nearly continuous. This behaviour is observed with all contact preserving fluxes like Roe, HLLC and HLL-CPS but not with the Rusanov flux. Due to the almost Lagrangian character of the scheme, the eigenvalue corresponding to the contact wave, $\lambda_2 = v - w$, is nearly zero, which leads to loss of dissipation in the corresponding characteristic field. If a spurious contact wave is generated during the violent dynamics, then this wave will be preserved by the scheme leading to wrong solutions. We modify the Roe scheme by preventing this eigenvalue from becoming too small or zero, which is similar to the approach used for the entropy fix. The eigenvalue $|\lambda_2|$ used in the dissipative part of the Roe flux is determined from

$$|\lambda_2| = \begin{cases} |v-w| & \text{if } |v-w| > \delta = \alpha c \\ \frac{1}{2}(\delta + |v-w|^2/\delta) & \text{otherwise} \end{cases}$$

With this modification and using $\alpha = 0.1$, the solution on moving mesh is shown in figure (5.6.3) and we do not observe the spurious contact discontinuities which arise with the standard Roe flux, while at the same time, the solution accuracy compares favourably with the previous results that did not use the eigenvalue fix.

We next compute the solutions using quadratic polynomials. Figure (5.6.4) shows the results obtained with the TVD limiter which shows the dramatically better accuracy that is achieved on moving mesh compared to static mesh. In figure (5.6.5) we perform the same computation with a WENO limiter taken from [**Zhong:2013:SWE:2397205.2397446**]. The static mesh results are now improved over the case of TVD limiter but still not as good as the moving mesh results in terms of capturing the extrema. In figure (5.6.6) we show a zoomed view of the results on moving mesh with TVD and WENO limiter. We see that the TVD limiter is also able to capture all the features and is almost comparable to the WENO limiter.

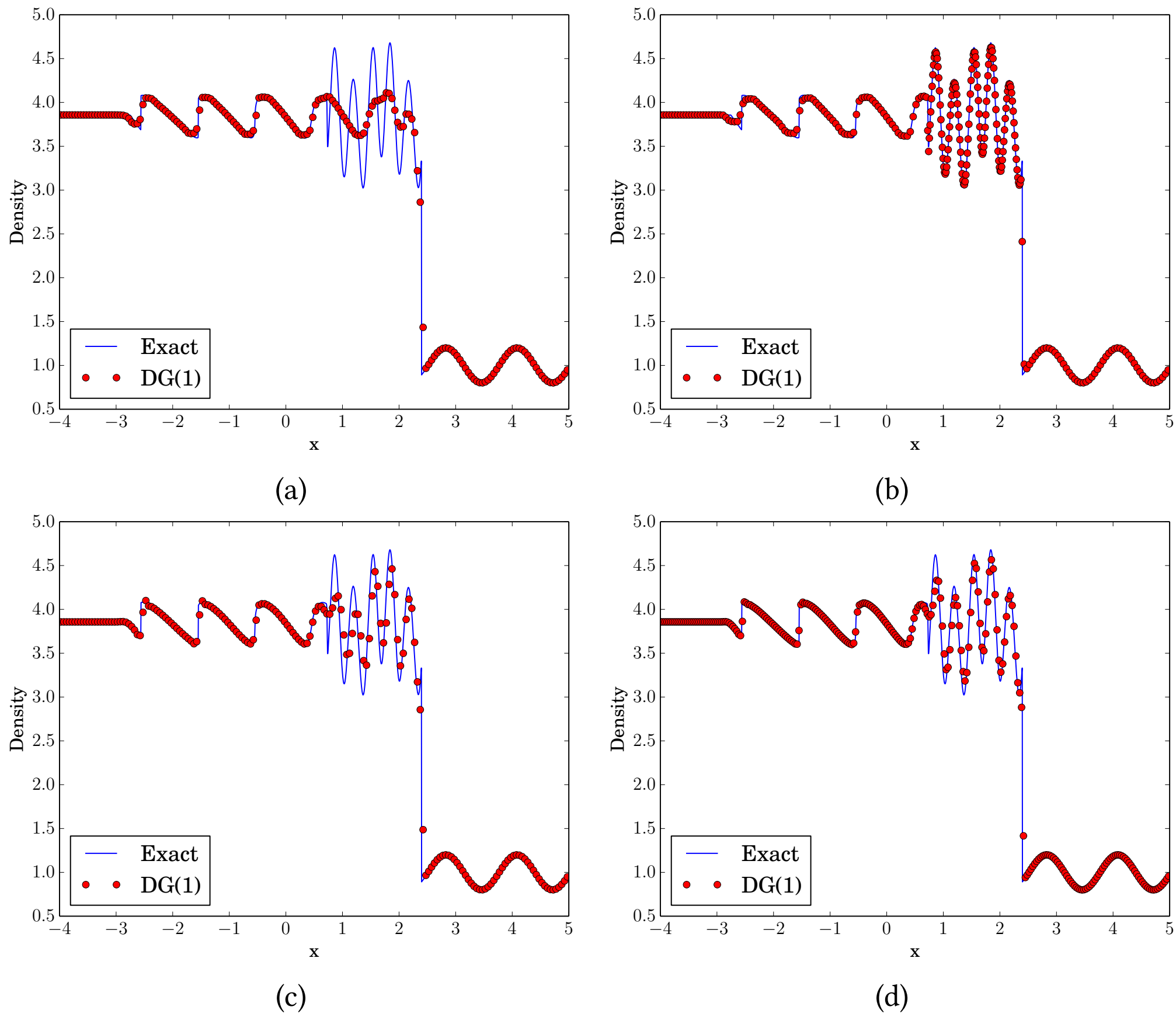

Figure 5.6.1: Shu-Osher problem using Roe flux: (a) static mesh, 200 cells, $M = 0$ (b) moving mesh, 200 cells, $M = 0$ (c) static mesh, 200 cells, $M = 100$ (d) static mesh, 300 cells, $M = 100$

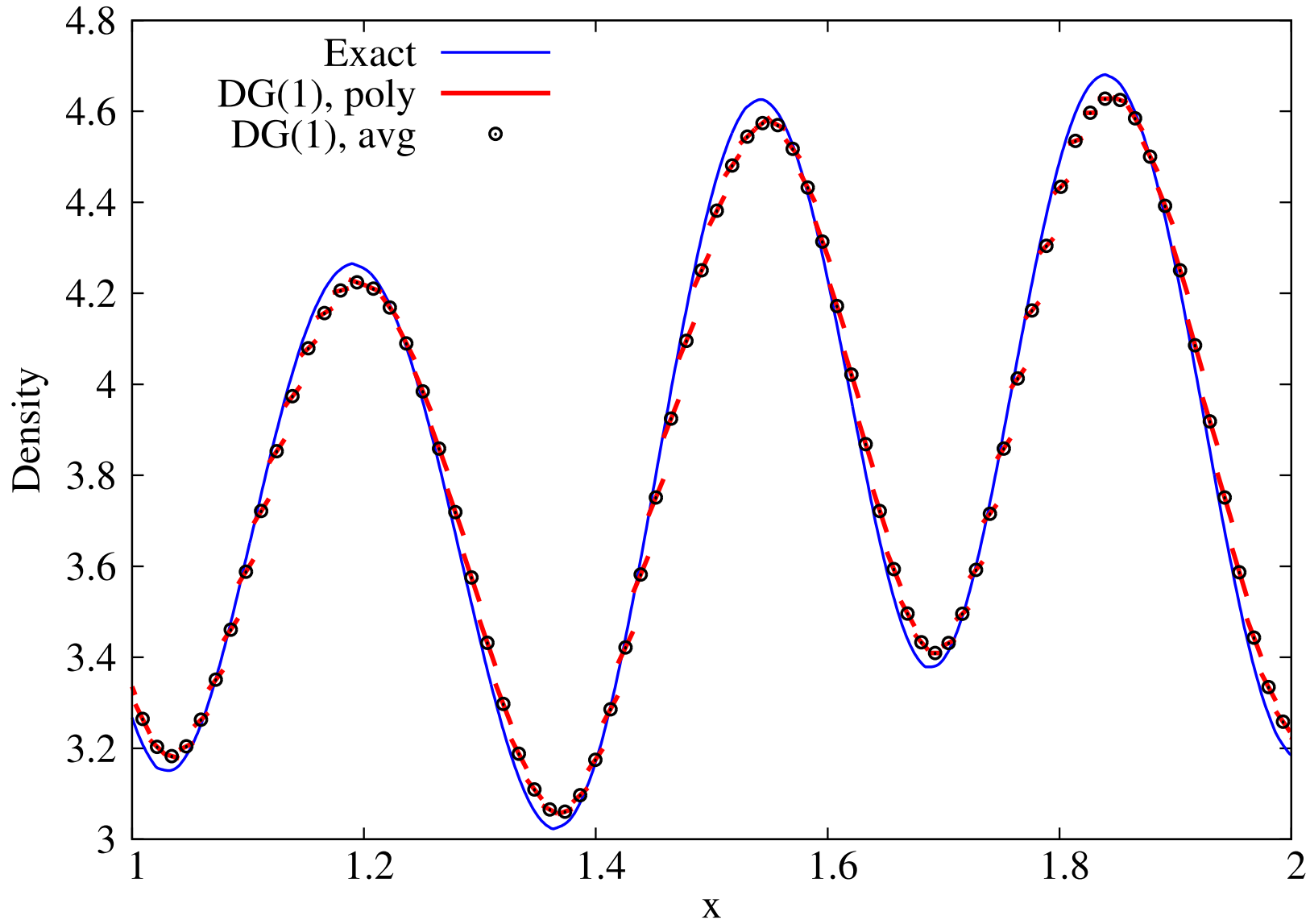

Figure 5.6.2: Shu-Osher problem using Roe flux on moving mesh

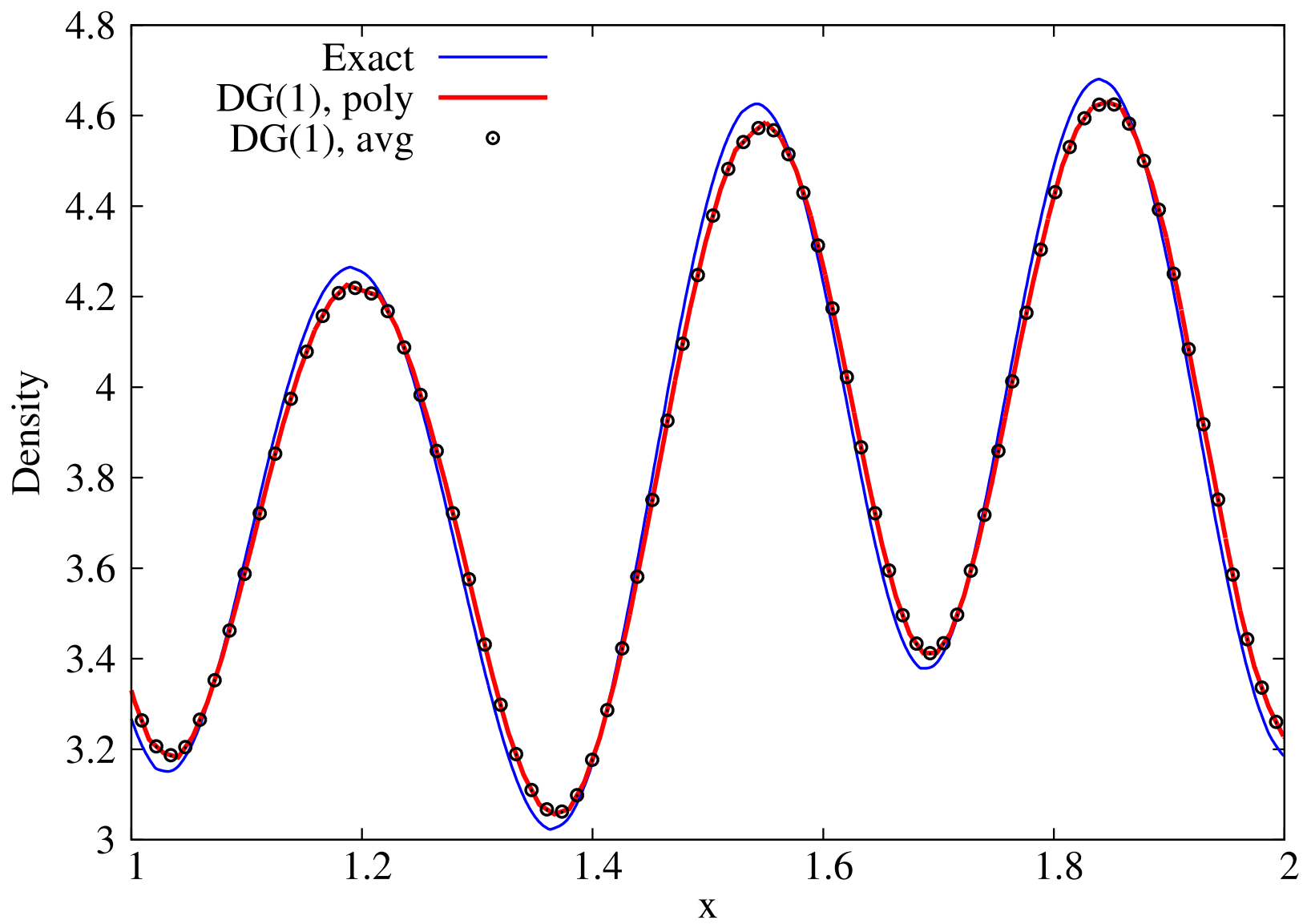

Figure 5.6.3: Shu-Osher problem using modified Roe flux on moving mesh

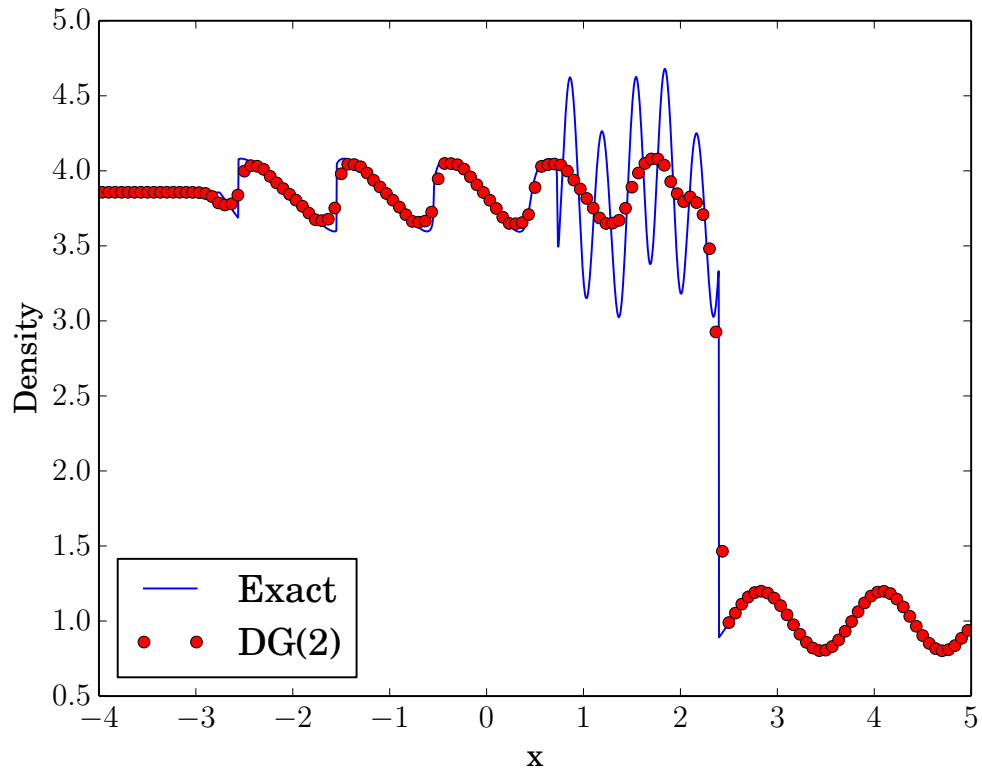

(a)

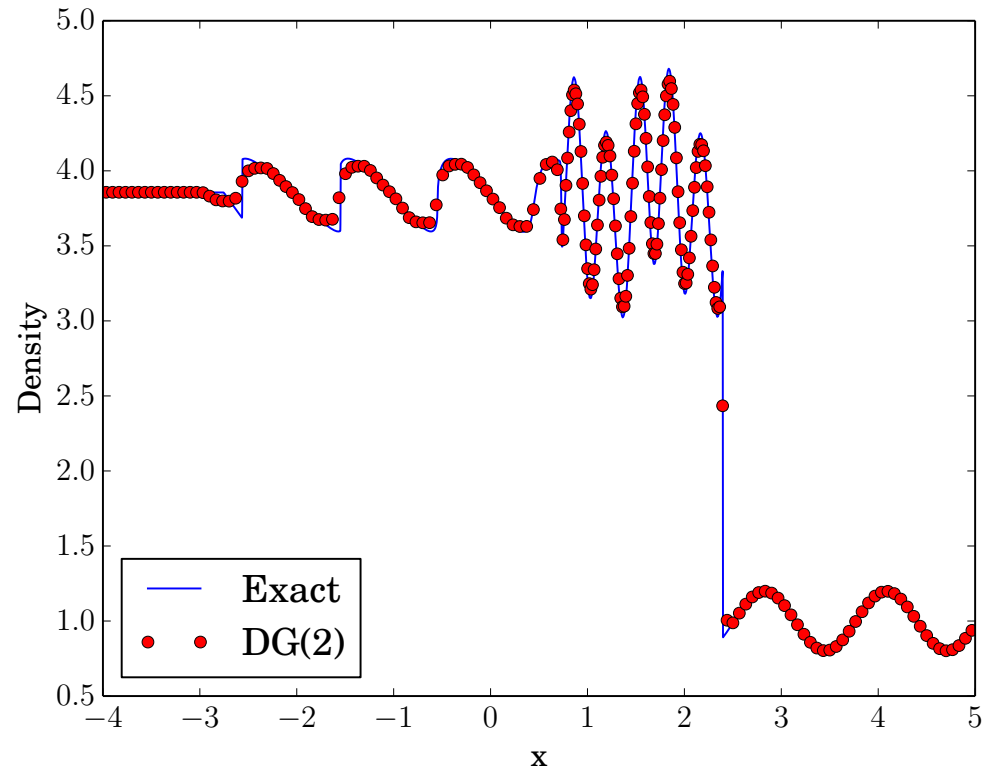

(b)

Figure 5.6.4: Shu-Osher problem using modified Roe flux, TVD limiter, quadratic polynomials and 150 cells. (a) static mesh, (b) moving mesh

## 5.7 Titarev-Toro problem

Titarev-Toro problem is an extension of the Shu-Osher problem [**Titarev2014**] to test a severely oscillatory wave interacting with a shock wave. It aims to test the ability of higher-order methods to capture the extremely high frequency waves. The initial condition is given by

$$(\rho, v, p) = \begin{cases} (1.515695, 0.523346, 1.805), & -5 < x \leq -4.5 \\ (1 + 0.1\sin(20\pi x), 0, 1), & -4.5 < x \leq 5 \end{cases} \tag{5.7.1}$$

The computation is carried out on a mesh of 1000 cells with the final time $T = 5$ and the density at this final time is shown in Figures (5.7.1), (5.7.2). The fixed mesh is not able to resolve the high frequency oscillations due to dissipation in the fluxes and the TVD limiter, but the ALE scheme gives an excellent resolution of these high frequency oscillations. Note that the ALE scheme also uses the same TVD limiter but it is still able to resolve the solution to a very degree of accuracy. This result again demonstrates the superior accuracy that can be achieved by using a nearly Lagrangian ALE scheme in problems involving interaction of shocks and smooth flow structures.

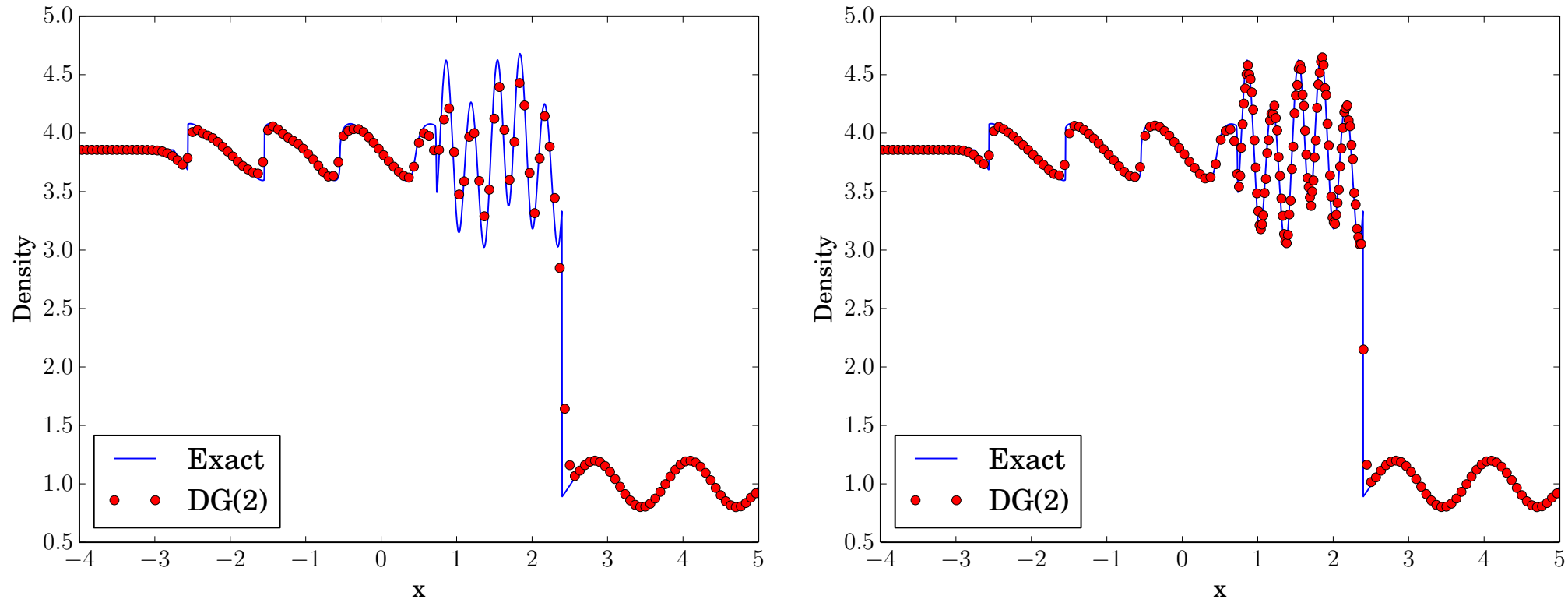

(a) (b)

Figure 5.6.5: Shu-Osher problem using modified Roe flux, WENO limiter, quadratic polynomials and 150 cells. (a) static mesh, (b) moving mesh

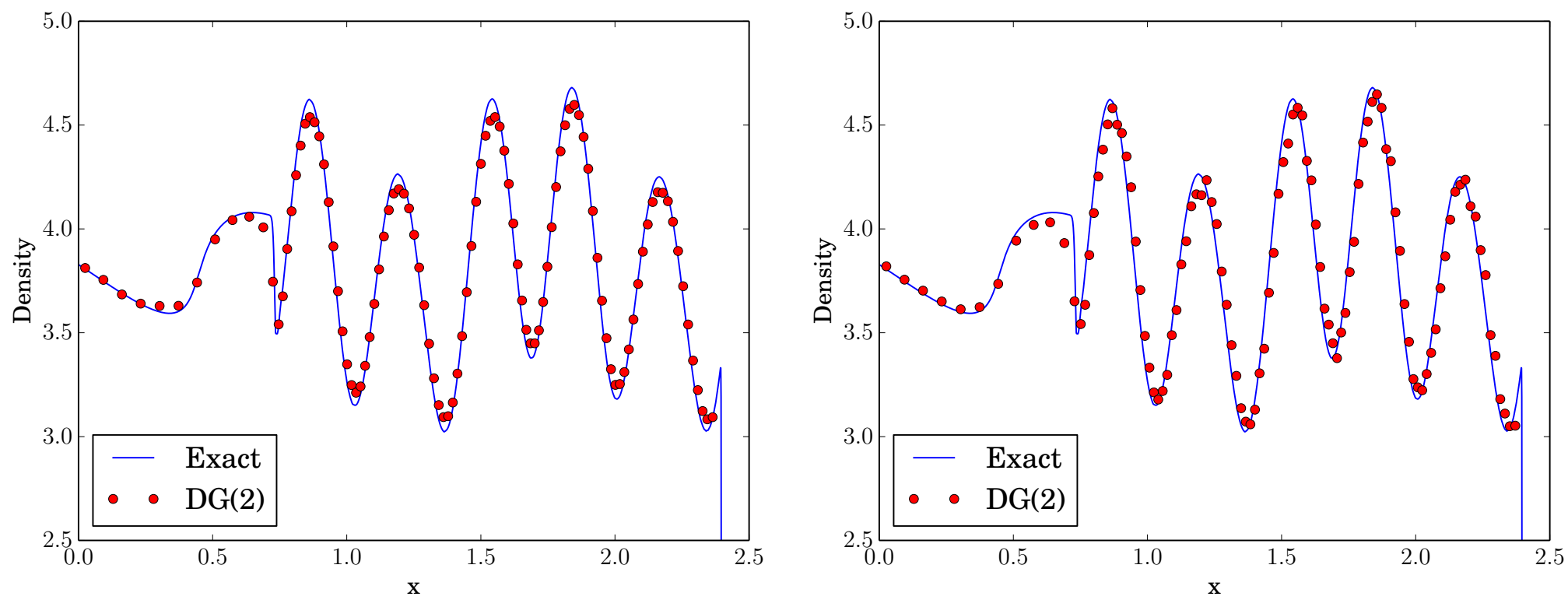

(a) (b)

Figure 5.6.6: Shu-Osher problem using modified Roe flux, moving mesh, quadratic polynomials and 150 cells. (a) TVD limiter, (b) WENO limiter

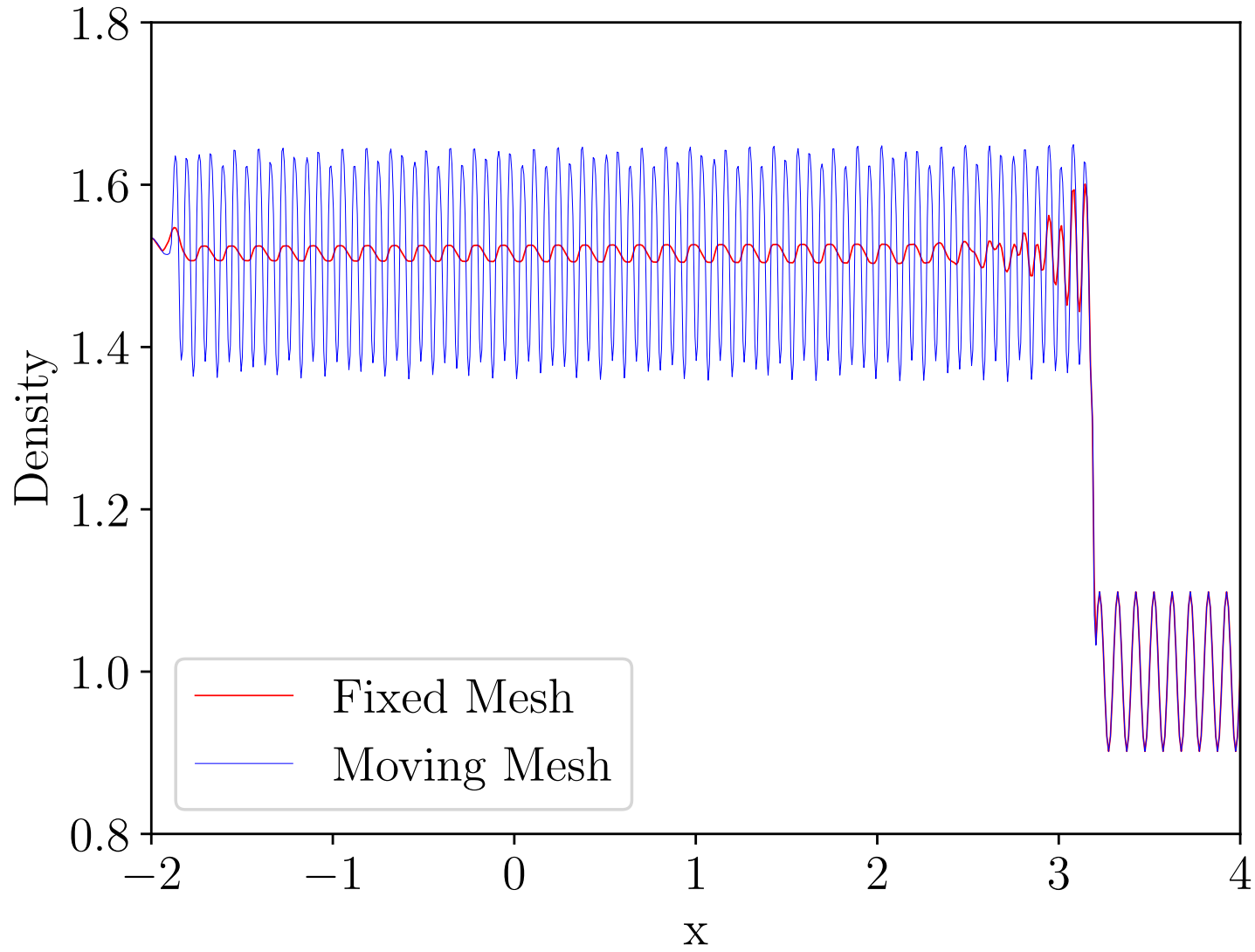

Figure 5.7.1: Titarev Problem with HLLC flux, 1000 cells and TVD limiter

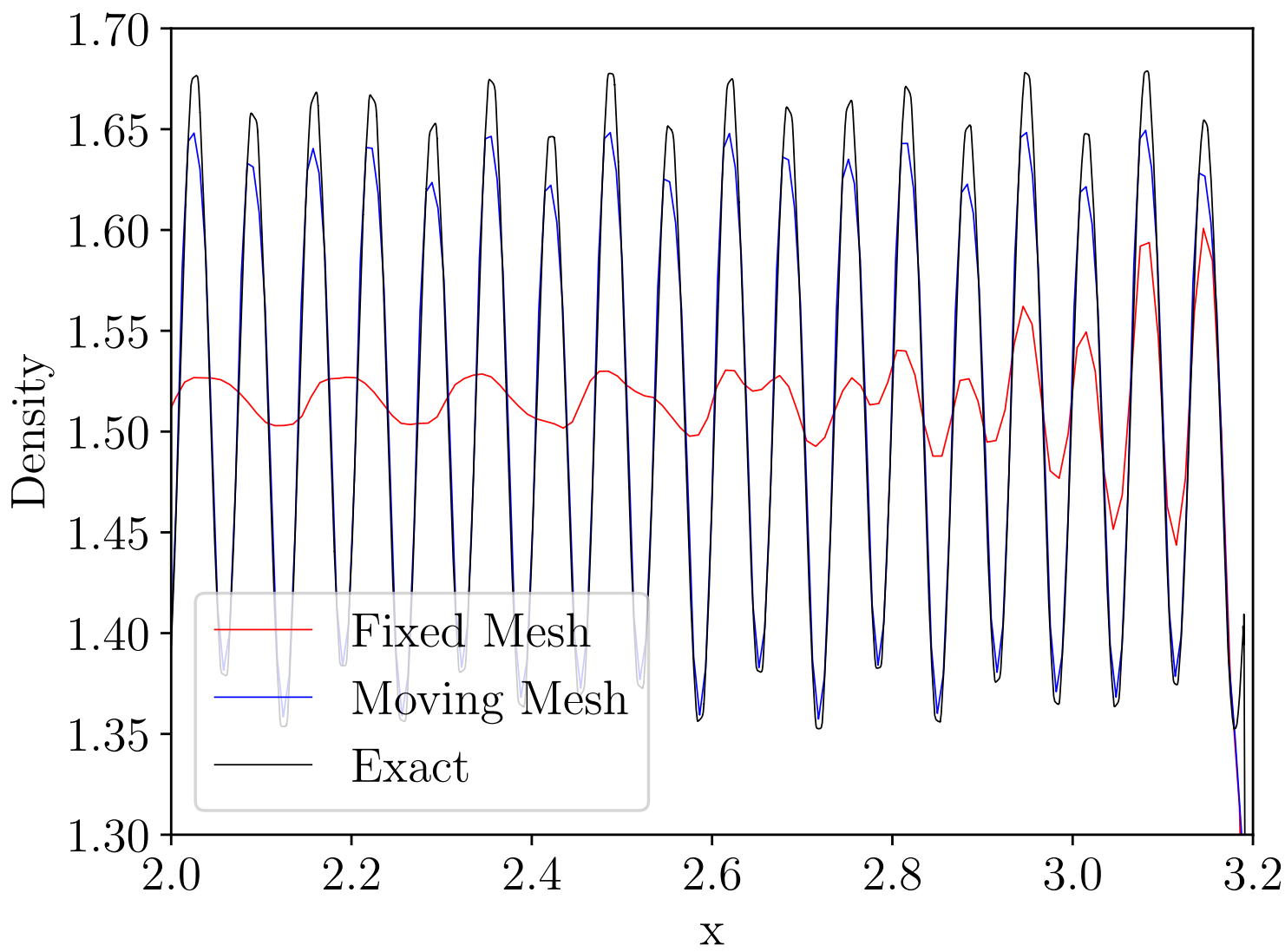

Figure 5.7.2: Titarev Problem with HLLC flux, 1000 cells and TVD limiter. (Zoomed Version)

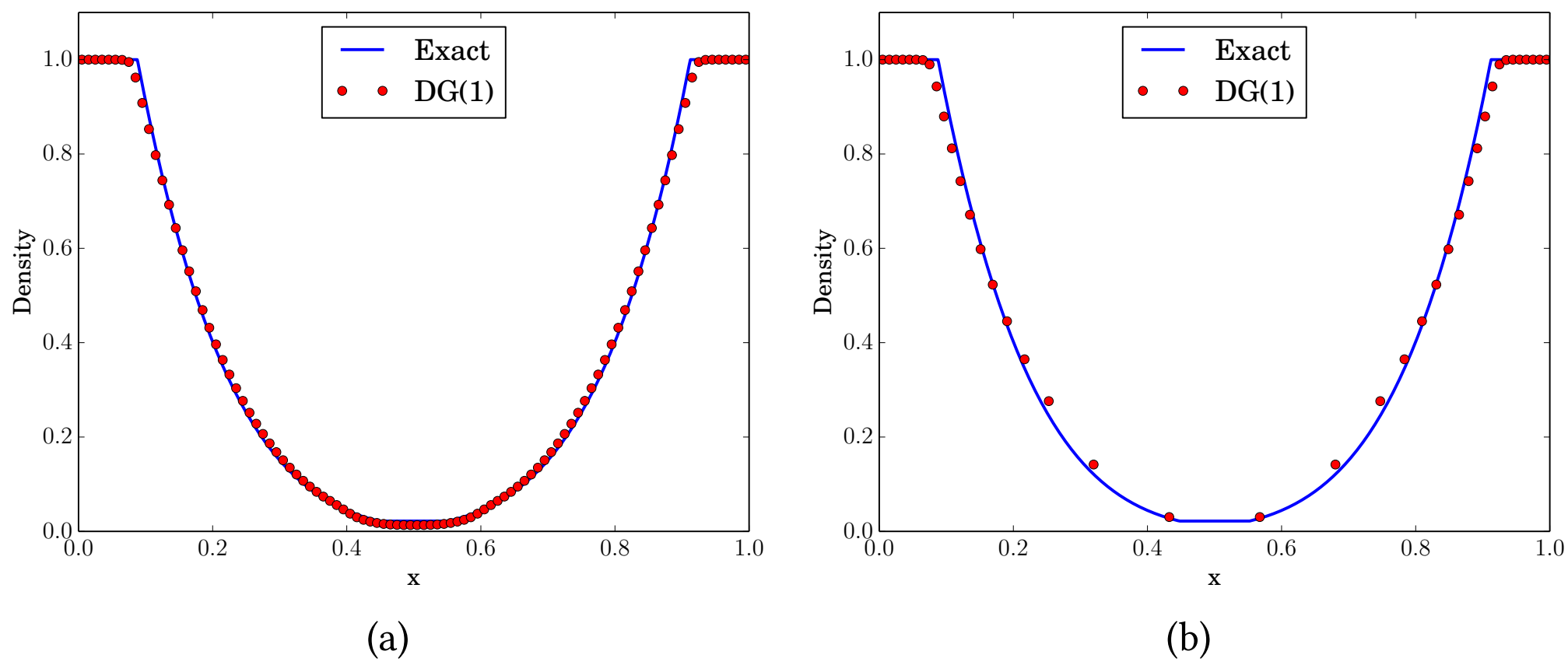

Figure 5.8.1: 123 problem using HLLC flux and 100 cells: (a) static mesh, (b) moving mesh

## 5.8 123 problem

The initial condition is given by [**torobook**]

$$(\rho, v, p) = \begin{cases} (1.0, -2.0, 0.4) & x < 0.5 \\ (1.0, +2.0, 0.4) & x > 0.5 \end{cases}$$

The computational domain is $[0, 1]$ and the final time is $T = 0.15$. The density using 100 cells is shown in figure (5.8.1) with static and moving meshes. The mesh motion does not significantly improve the solution compared to the static mesh case since the solution is smooth. On the contrary, the mesh becomes rather coarse in the expansion region, though the solution is still well resolved. However, severe expansion may lead to very coarse meshes which may be undesirable. To prevent very coarse cells, we switch on the mesh refinement algorithm as described before and use the upper bound on the mesh size as $h_{\max} = 0.05$. The resulting solution is shown in figure (5.8.2) where the number of cells has increased to 108 at the time shown. The central expansion region is now resolved by more uniformly sized cells compared to the case of no grid refinement.

## 5.9 Blast problem

The initial condition is given by [**Woodward1984115**]

$$(\rho, v, p) = \begin{cases} (1.0, 0.0, 1000.0) & x < 0.1 \\ (1.0, 0.0, 0.01) & 0.1 < x < 0.9 \\ (1.0, 0.0, 100.0) & x > 0.9 \end{cases}$$

with a domain of $[0, 1]$ and the final time is $T = 0.038$. A reflective boundary condition is used at $x = 0$ and $x = 1$. A mesh of 400 cells is used for this simulation and in case of moving mesh, we perform grid adaptation with $h_{min} = 0.001$ since some cells be-come very small during the collision of the two shocks. The positivity preserving limiter of [**Zhang:2010:PHO:1864819.1865066**] is applied together with TVD limiter and HLLC flux. The static mesh results shown in figure (5.9.1a) indicate too much numerical viscosity in the contact wave around $x = 0.6$. This wave is more accurately resolved in the moving mesh scheme as seen in figure (5.9.1b) which is an advantage due to the ALE scheme and is a very good indicator of the scheme accuracy as this is a very challenging feature to compute accurately. We next compute the same problem using quadratic polynomials with all other parameters being as before. The solutions are shown in figure (5.9.2) and indicate that the Lagrangian moving mesh scheme is more accurate in resolving the contact discontinuity. The higher polynomial degree does not show any major improvement in the solution compared to the linear case, which could be a consequence of the strong shock interactions present in this problem, see figure (4.11-4.12) in [**Zhong:2013:SWE:2397205.2397446**] and figure (3.7) in [**Zhu2016110**] in comparison to current results.

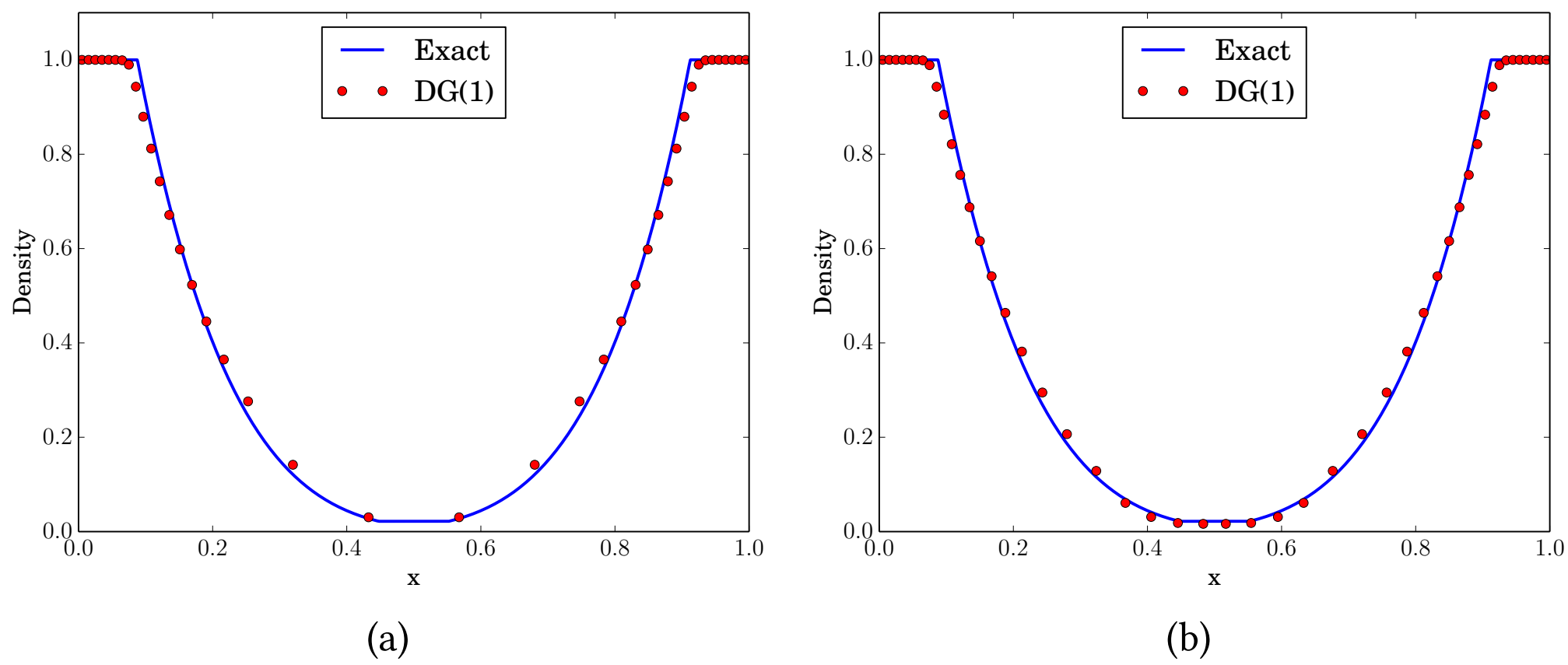

Figure 5.8.2: 123 problem using HLLC flux and grid refinement: (a) static mesh, (b) moving mesh with mesh adaptation ($h_{\max} = 0.05$) leading to 108 cells at final time

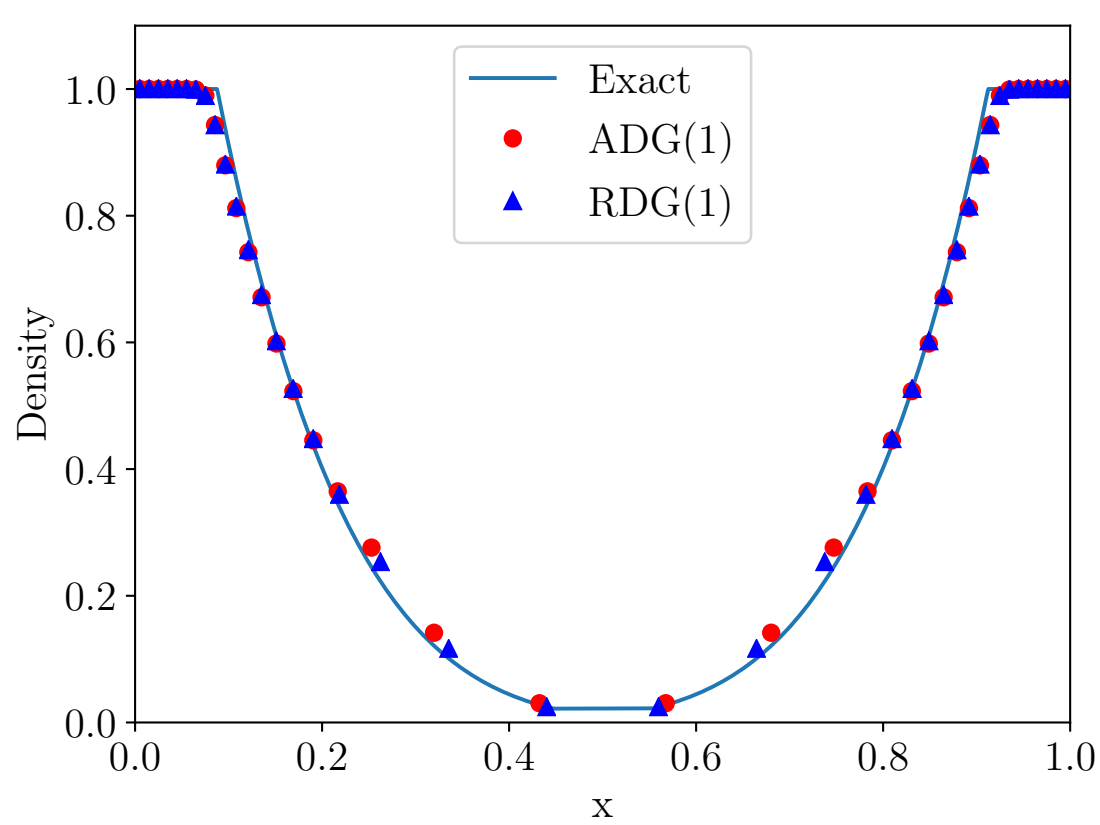

Figure 5.8.3: 123 problem using HLLC flux, 100 cells and TVD limiter. ADG : Average Velocity, RDG : Linearized Riemann Velocity

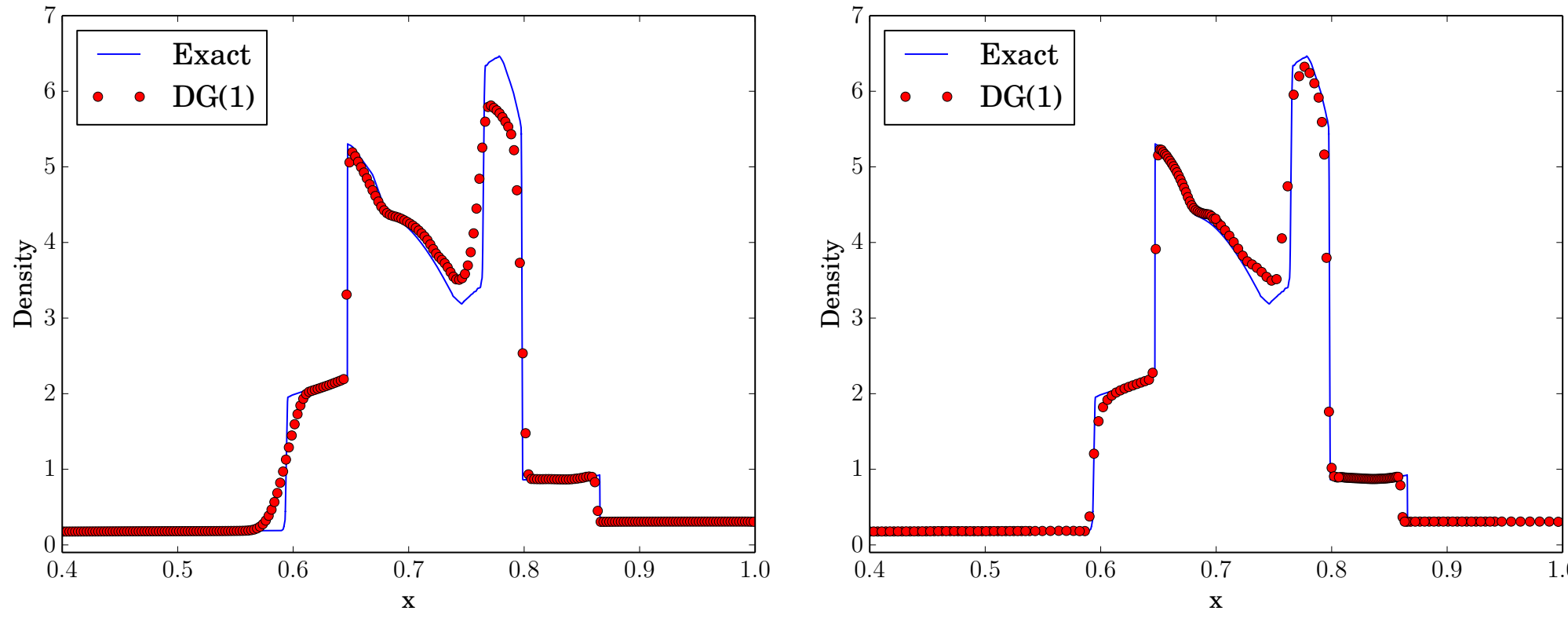

(a) (b)

Figure 5.9.1: Blast problem using HLLC flux and 400 cells. (a) static mesh, (b) moving mesh with adaptation ($h_{min} = 0.001$) leading to 303 cells at final time.

## 5.10 Le Blanc shock tube test case

The Le Blanc shock tube test case is an extreme shock tube problem where the initial discontinuity separates a region of high energy and density from one of low energy and density. This is a much more severe test than the Sod problem and hence more challenging for numerical schemes. The computational domain is $0 \le x \le 9$ and is filled with an ideal gas with $\gamma = 5/3$. The gas is initially at rest and we perform the simulation up to a time of $T = 6$ units. The initial discontinuity is at $x = 3$ and the initial condition is given by

$$(\rho, v, p) = \begin{cases} (1.0, 0.0, 0.1) & \text{if } x < 3 \\ (0.001, 0.0, 10^{-7}) & \text{if } x > 3 \end{cases} \tag{5.10.1}$$

Note that both the density and pressure have a very large jump in the initial condition. The solution that develops from this initial condition consists of a rarefaction wave moving to the left and a contact discontinuity and a strong shock moving to the right. In Figure (5.10.1), we show the comparison of the internal energy profile at final time between a fixed mesh solution and moving mesh solutions with two different mesh velocities as described before. Most methods tend to generate a very large spike in the internal energy in the contact region, e.g., compare with Figure (11) in [**Shashkov2005**], while the present ALE method here is able to give a better profile. We plot the pressure profile in Figure (5.10.2) which shows that the ALE scheme is able to better represent the region around the contact wave as compared to fixed mesh method.

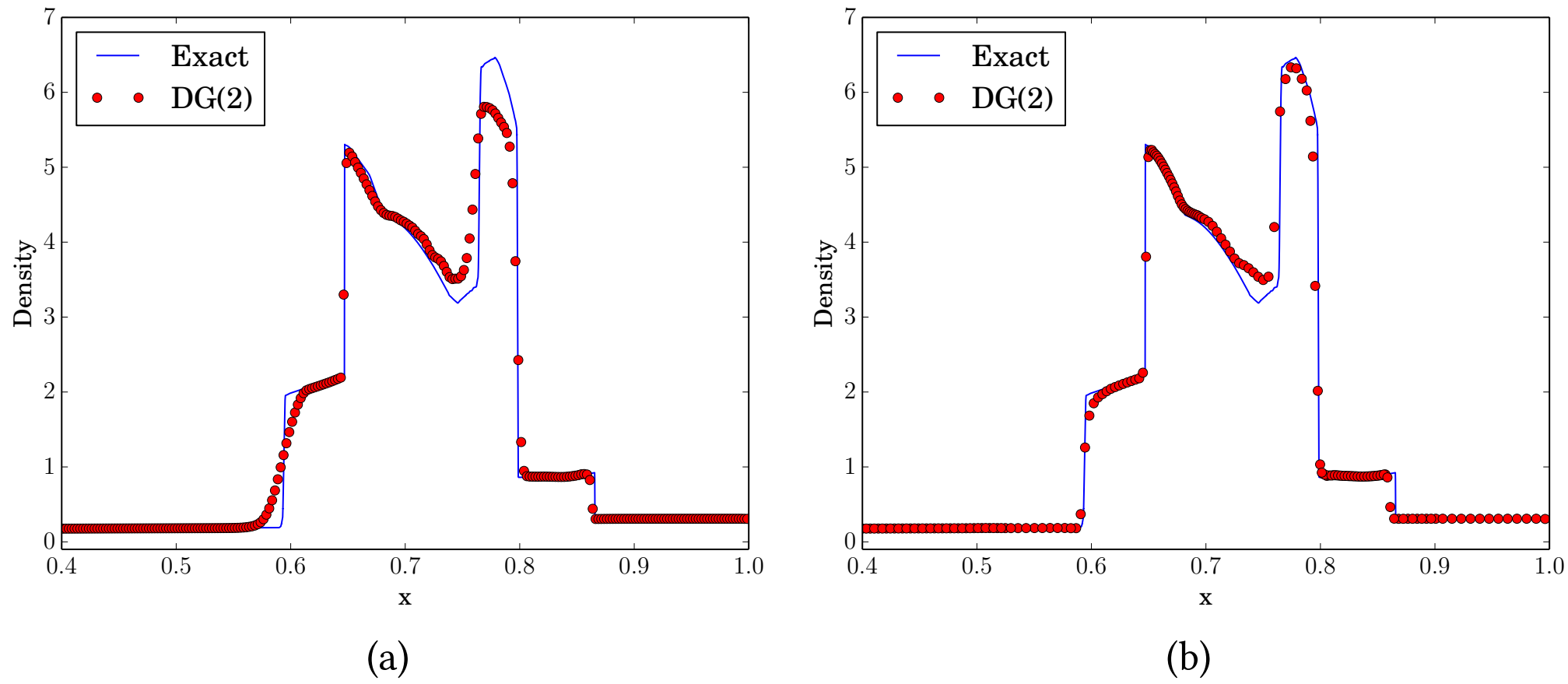

Figure 5.9.2: Blast problem using HLLC flux, quadratic polynomials and 400 cells. (a) static mesh, (b) moving mesh with adaptation ($h_{min} = 0.001$) leading to 293 cells at final time.

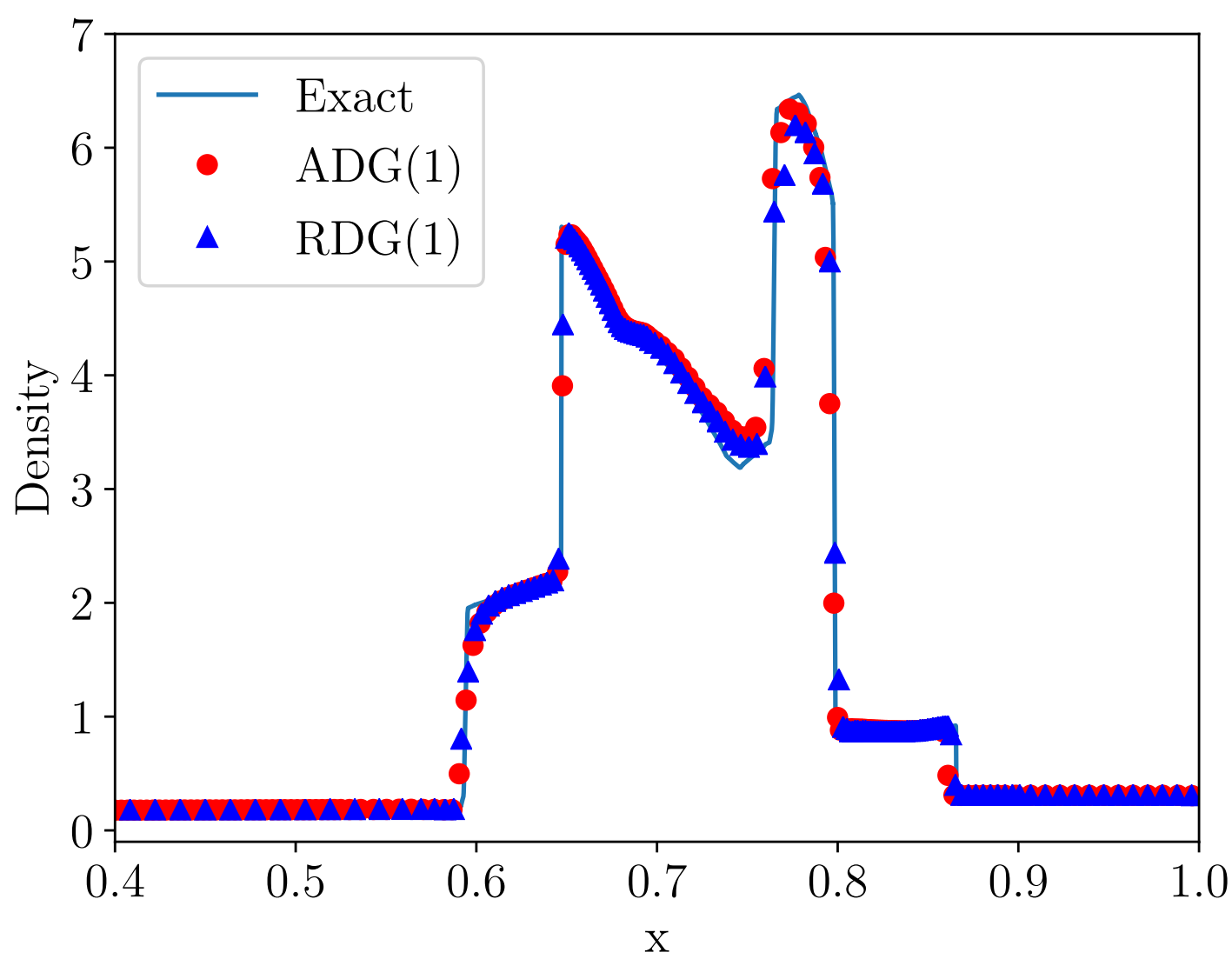

Figure 5.9.3: Blast problem using HLLC flux, 100 cells and TVD limiter. ADG : Average Velocity, RDG : Linearized Riemann Velocity

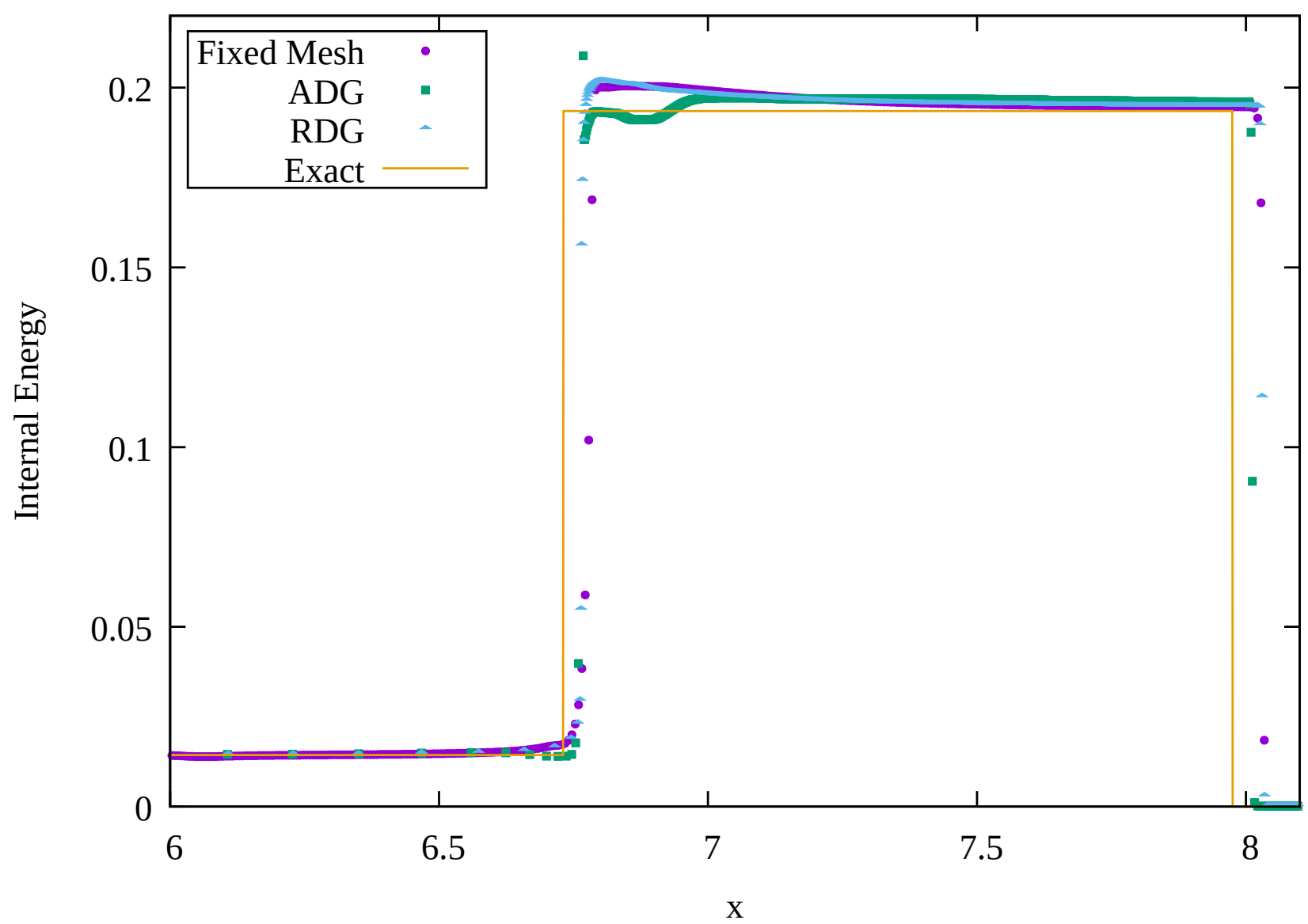

Figure 5.10.1: Internal energy for Le Blanc Shock Tube with Rusanov flux, 1400 cells and TVD limiter, ADG : Average Velocity, RDG : Linearized Riemann Velocity

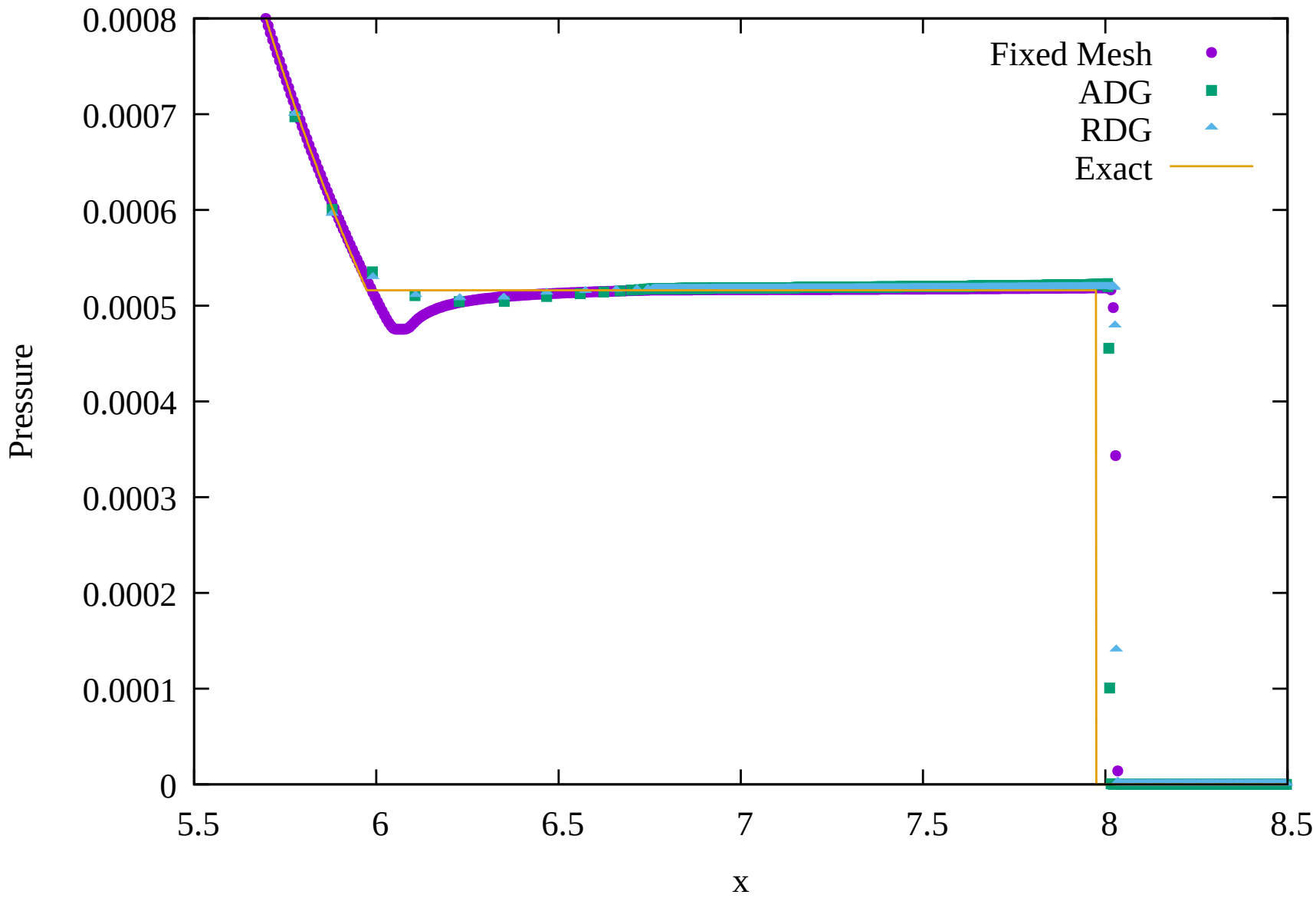

Figure 5.10.2: Pressure for Le Blanc Shock Tube with Rusanov flux, 1400 cells and TVD limiter, ADG : Average Velocity, RDG : Linearized Riemann Velocity

Chapter

# 6 Two Dimensional Applications

The extension to two dimensions involves two aspects that need to be addressed. The first issue is how to handle the grid motion and the second is how to formulate the ALE-DG scheme. The second part is a natural generalization of the DG scheme we have described for the 1-D case in this paper, except that we have to construct basis functions on triangles and perform some numerical quadrature. We now consider the two dimensional Euler equations written as

$$\frac{\partial u}{\partial t} + \frac{\partial f(u)}{\partial x} + \frac{\partial g(u)}{\partial y} = 0 \tag{6.0.1}$$

where

$$u = \begin{bmatrix} \rho \\ \rho u \\ \rho v \\ E \end{bmatrix} \qquad f(u) = \begin{bmatrix} \rho u \\ p + \rho u^2 \\ \rho uv \\ (E+p)u \end{bmatrix} \qquad g(u) = \begin{bmatrix} \rho v \\ \rho uv \\ p + \rho v^2 \\ (E+p)v \end{bmatrix} \tag{6.0.2}$$

$$p = (\gamma - 1)\left[E - \frac{1}{2}\rho(u^2 + v^2)\right] \tag{6.0.3}$$

## 6.1 Isentropic Vertex

The test case we consider involves an isentropic vortex that is advecting with constant velocity and is a smooth solution for which error norms can be calculated. The test is carried out on a square domain $[-10, 10] \times [-10, 10]$ with periodic boundary conditions. The initial conditions

is an isentropic vortex

$$T = 1 - \frac{(\gamma-1)\beta^2}{8\gamma\pi^2} e^{1-r^2} \tag{6.1.1}$$

$$\rho = T^{\frac{1}{\gamma-1}} \tag{6.1.2}$$

$$u = u_\infty - \frac{\beta}{2\pi} y e^{\frac{1-r^2}{2}} \tag{6.1.3}$$

$$v = v_\infty - \frac{\beta}{2\pi} y e^{\frac{1-r^2}{2}} \tag{6.1.4}$$

$$p = \rho^\gamma \tag{6.1.5}$$

with $u_\infty = 1, v_\infty = 0, \gamma = 1.4, \beta = 10$. As the solution evolves in time, the mesh becomes quite deformed because the vortex is continually shearing the mesh, which can lead to degenerate meshes, as shown in figure 6.1.1. With the ALE DG method, we are able to maintain a good mesh quality even after the vortex has rotated 20 times around its center as shown in figures 6.1.2. As the vortex is translating, we plot the solution in a window centered at the vortex center. We can see that the method maintains its high order of accuracy from the convergence rates of the error shown in table 6.1.2; using linear basis functions yields second order convergence while quadratic basis functions lead to third order convergence. We can also see that the error in the moving mesh method evolves considerably slower than the error in the fixed mesh mesh.

| $N$ | $k=1$ | | $k=2$ | |
|---|---|---|---|---|
| | Error | Rate | Error | Rate |
| 50x50 | 2.230e-03 | | 1.762E-04 | |
| 100x100 | 5.987E-04 | 1.945 | 2.305E-05 | 2.934 |
| 200x200 | 1.498E-04 | 1.998 | 2.973E-06 | 2.955 |
| 400x400 | 3.786E-05 | 1.984 | 3.762E-07 | 2.982 |
| 800x800 | 9.617E-06 | 1.977 | 3.474E-08 | 2.991 |

Table 6.1.1: Isentropic Vortex in 2D: Order of accuracy study on two dimensional static mesh

| $N$ | $k=1$ | | $k=2$ | |
|---|---|---|---|---|
| | Error | Rate | Error | Rate |
| 50x50 | 1.110e-03 | | 1.762E-04 | |
| 100x100 | 2.847E-04 | 1.963 | 2.305E-05 | 2.934 |
| 200x200 | 7.231E-05 | 1.977 | 2.973E-06 | 2.955 |
| 400x400 | 1.863E-05 | 1.956 | 3.762E-07 | 2.982 |
| 800x800 | 9.617E-06 | 1.977 | 3.474E-08 | 2.991 |

Table 6.1.2: Isentropic Vortex in 2D: Order of accuracy study on two dimensional moving mesh

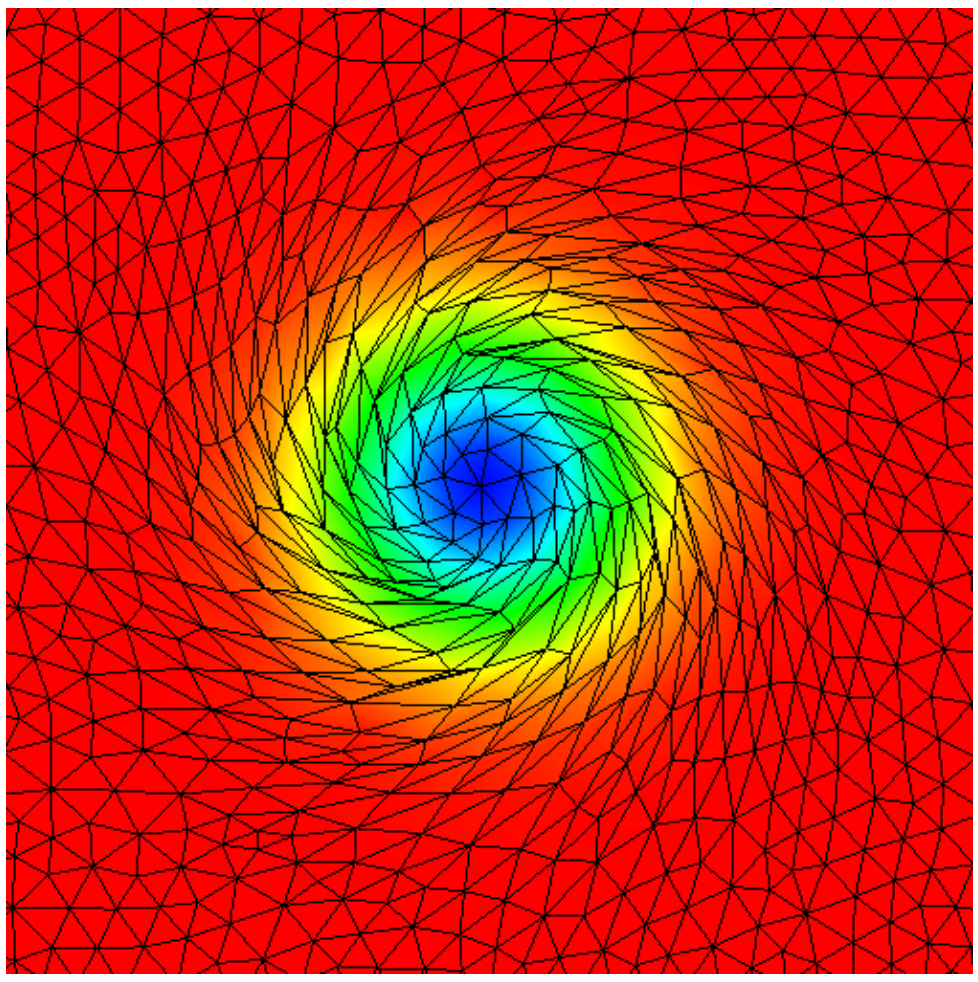

Figure 6.1.1: Isentropic Vortex in 2-D: Skewed Mesh without Remeshing $t = 2.660534$

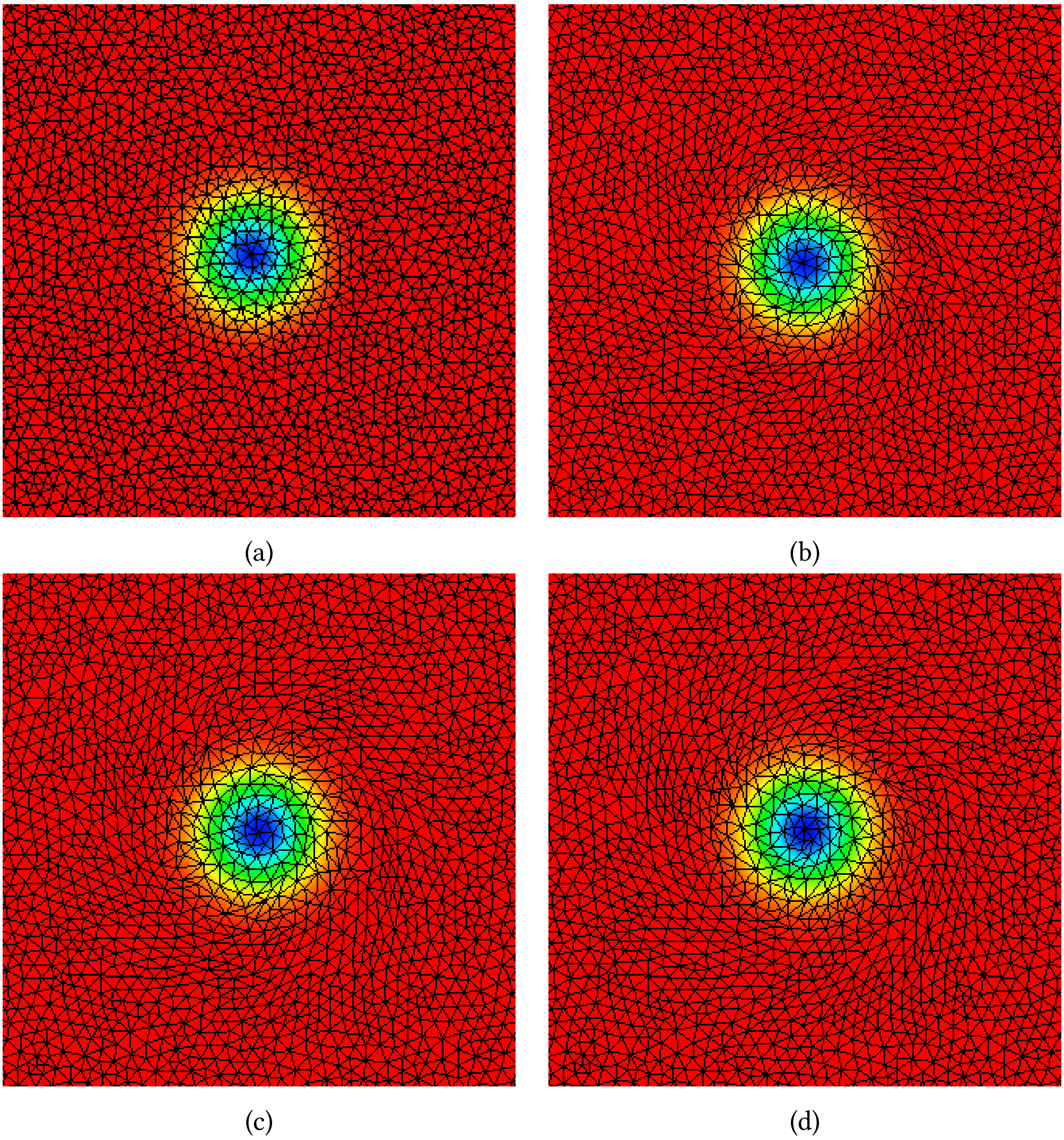

Figure 6.1.2: Isentropic vortex in 2-D: Mesh and pressure solution at various times (a) $t = 0$ (b) $t = 25$ (c) $t = 50$ (d) $t = 100$

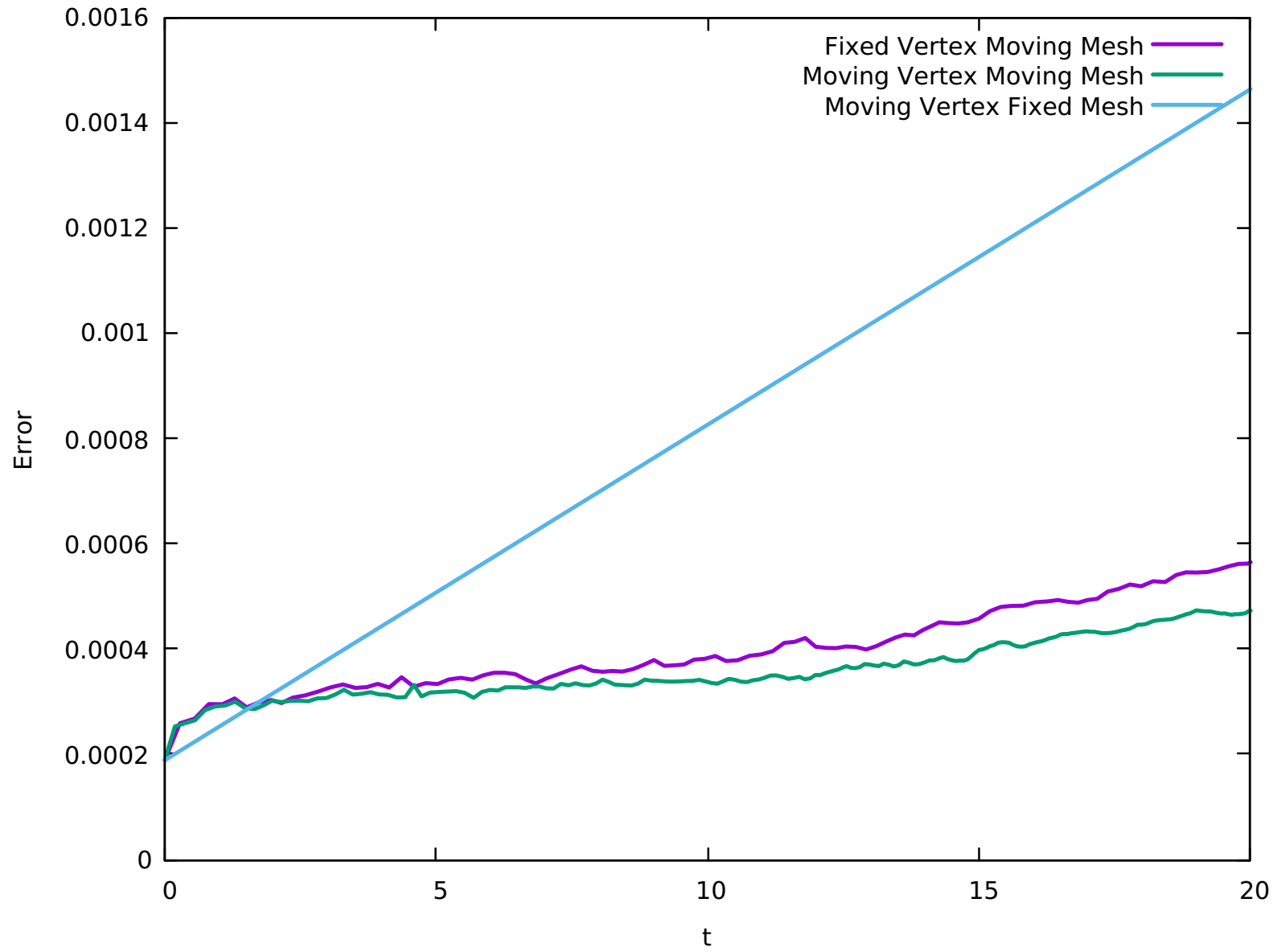

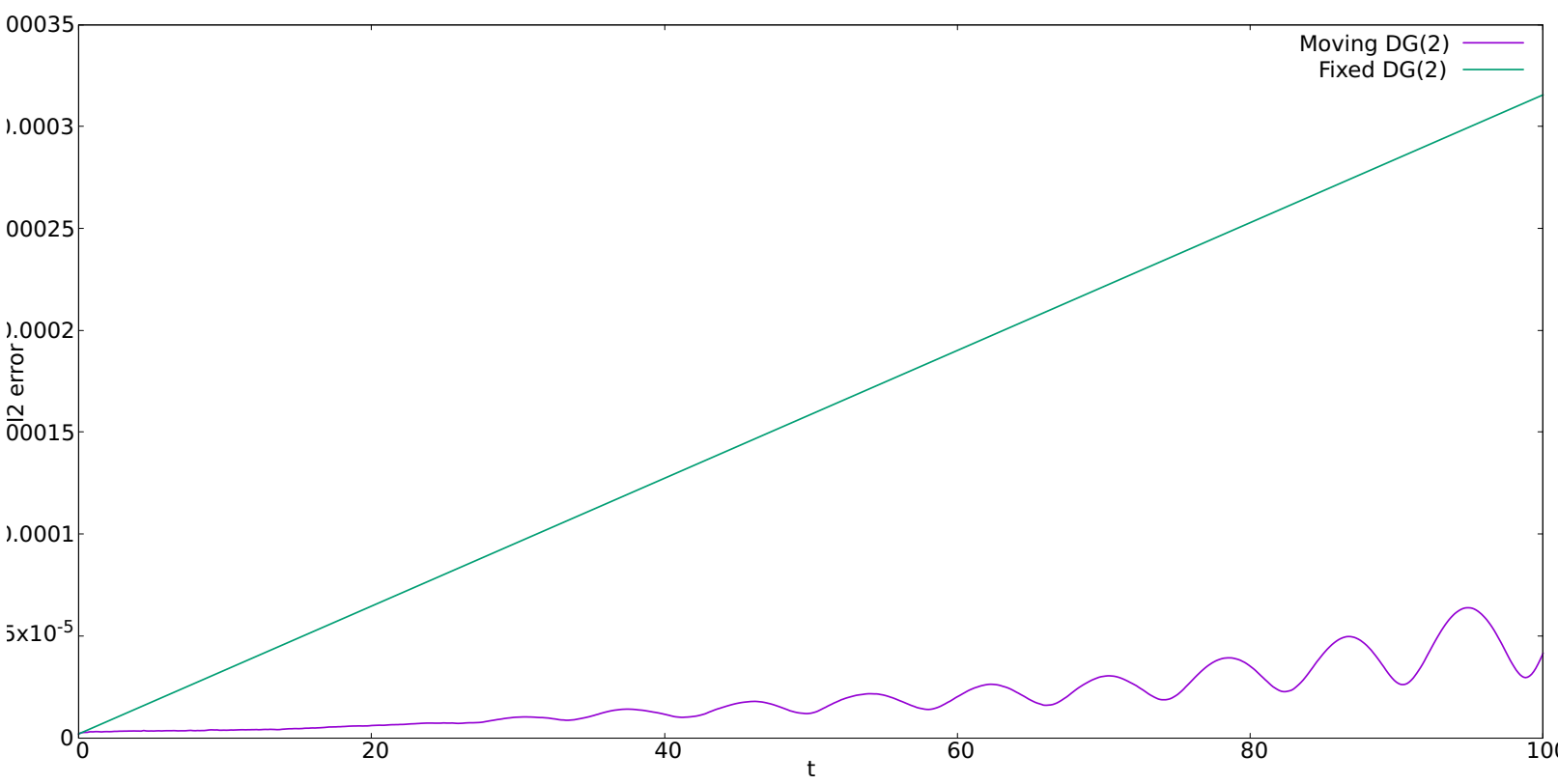

## 6.1.1 Isentropic Vortex with Non-Uniform Mesh

Due to the self-adaption nature of the moving mesh method, it is possible to use a non-uniform mesh to simulate a translating isentropic vortex. We start with a mesh which is fine near the vortex and coarse elsewhere. The vortex itself moves over the course of the simulation and we show that the mesh tracks the vortex. This allows us to concentrate the mesh in places with more features and thus reducing the computational cost of the code.

In comparison with uniform meshes, the non-uniform mesh only requires 1000 cells in comparison to uniform mesh's 10,000 cells to achieve a similar accuracy of error of $1.6e-4$. This represents a significant reduction in computational costs.

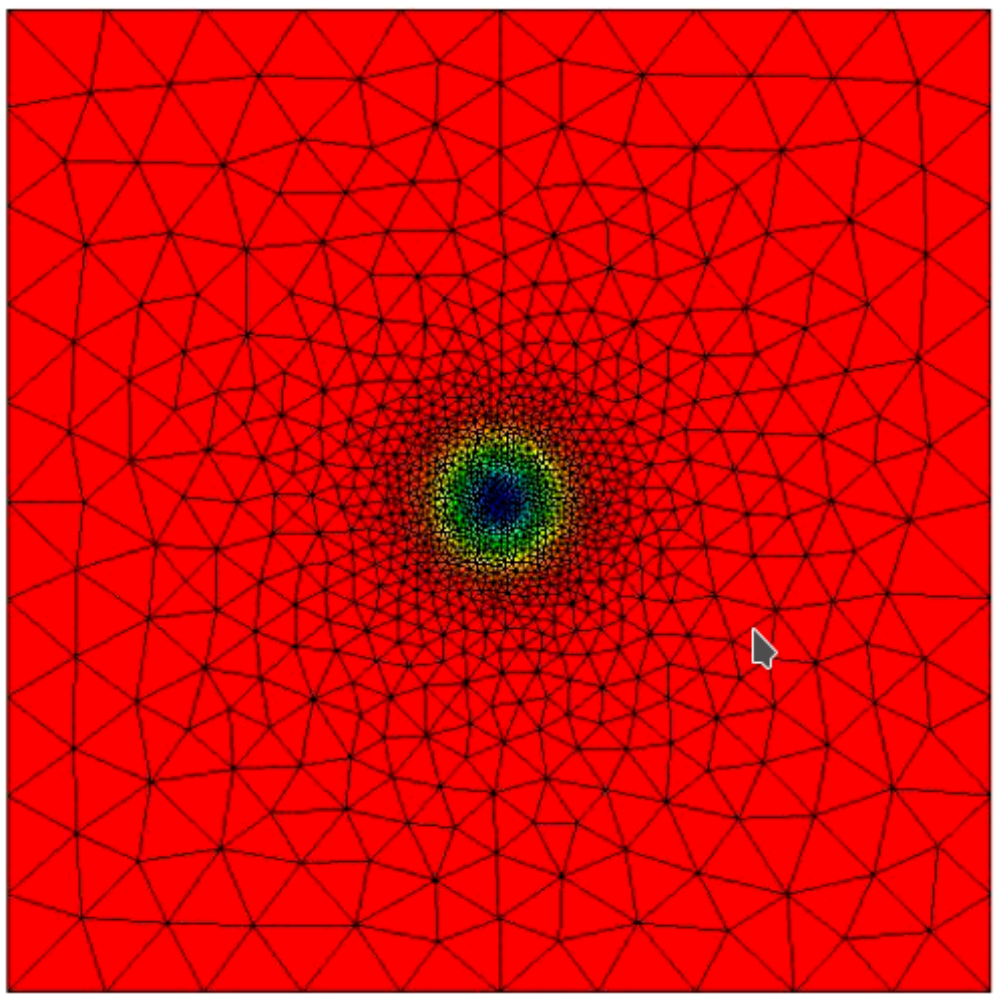

Figure 6.1.3: Isentropic Vortex with Non-Uniform Mesh, $T = 0$

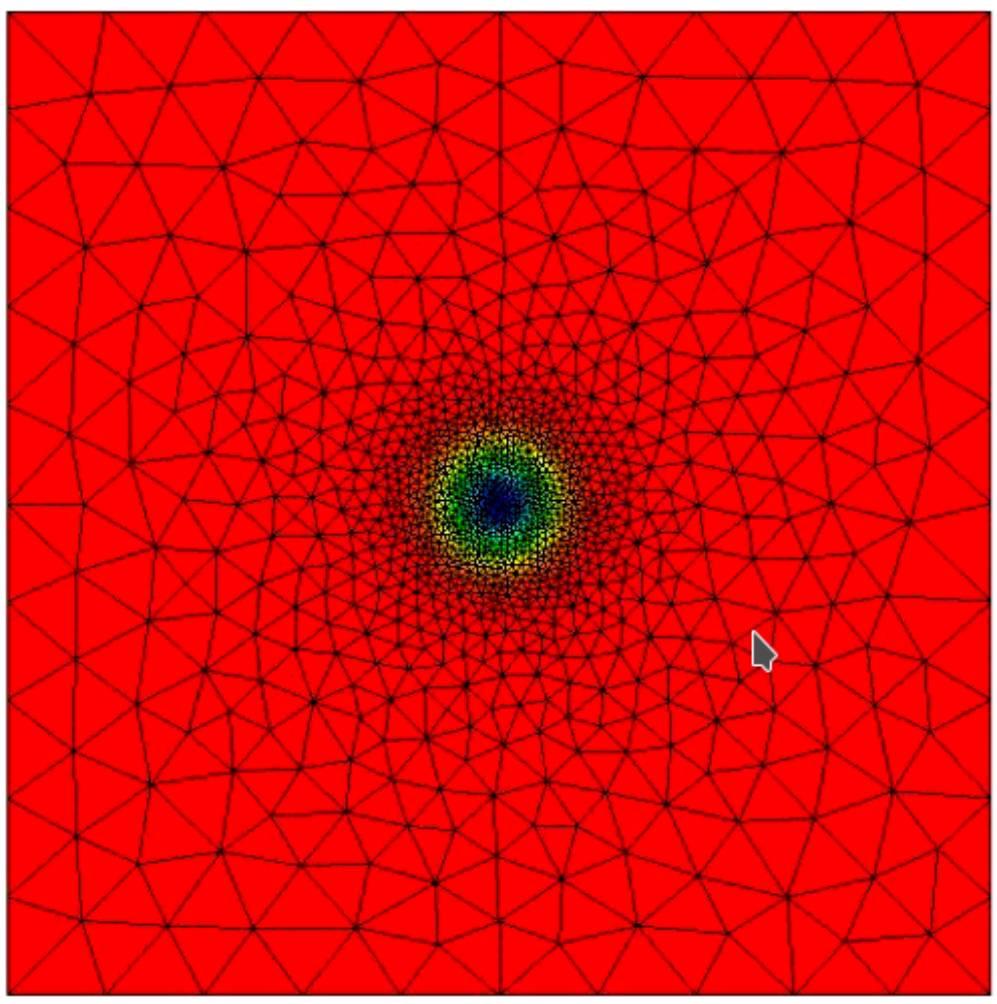

Figure 6.1.4: Isentropic Vortex with Non-Uniform Mesh, $T = 100$

### 6.1.2 Comparison of Computational Costs between Methods

The moving of the mesh and grid adaptation introduces some computational overhead. In case of smooth solutions, due to the absence of sharp discontinuities, the quality of solution itself does not significantly change. Therefore, it is useful to look at the improvement in error against the increased computational cost for the moving mesh method. For that purpose, we have analysed the time spent by the scheme in different sections of the code.

We observe that the velocity smoothing algorithm is the most significant cost for the moving mesh algorithm and reduces the effectiveness of moving mesh methods for smooth solutions.

However, the problem can be resolved by using non-uniform meshes as shown in . Another method to resolve the computational cost would be to improve the velocity smoothing algorithm in manner which will be detailed in REF CONCLUSIONS.

## 6.2 Sod Shock Tube Problem

A two dimensional Sod shock tube problem is the same problem as the one-dimensional problem described in REFER TO ONED but in a two dimensional domain. In this problem, we have a domain of $[0, 1] \times [0, 0.1]$ with periodic boundaries at the top and the bottom. We have chosen a transfinite mesh for the problem in order to accurately capture the initial condition. There are total of $200x40$ cells in a configuration of $100$ blocks along the $x$ and $20$ block along the $y$ direction. The final time of the simulation is $0.2$.

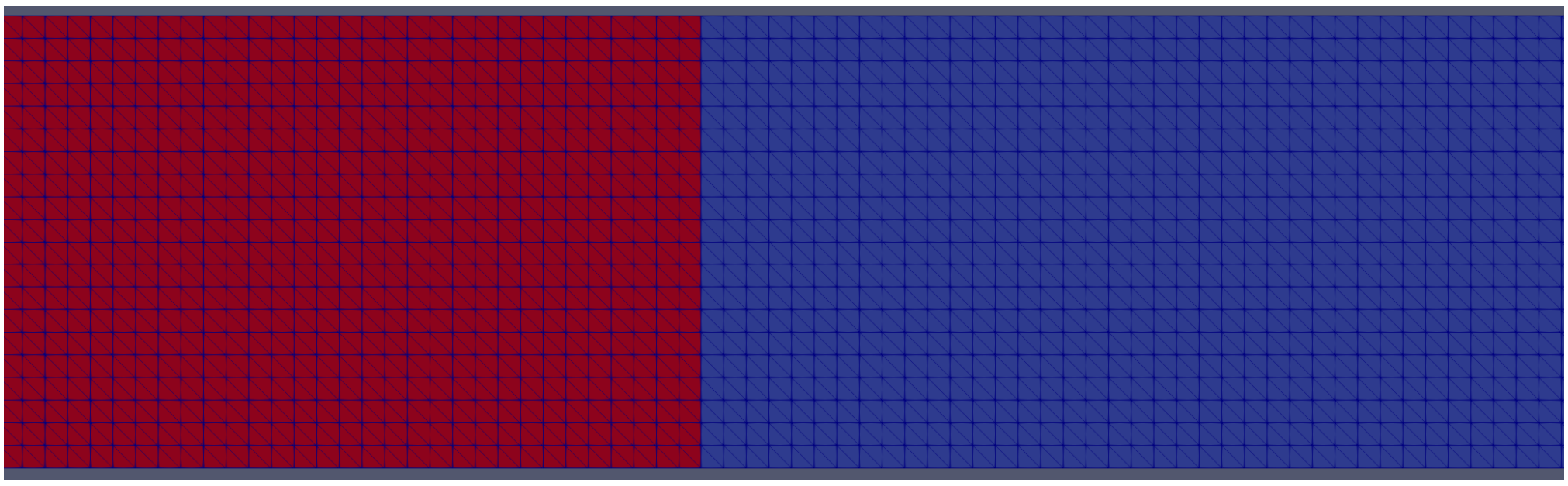

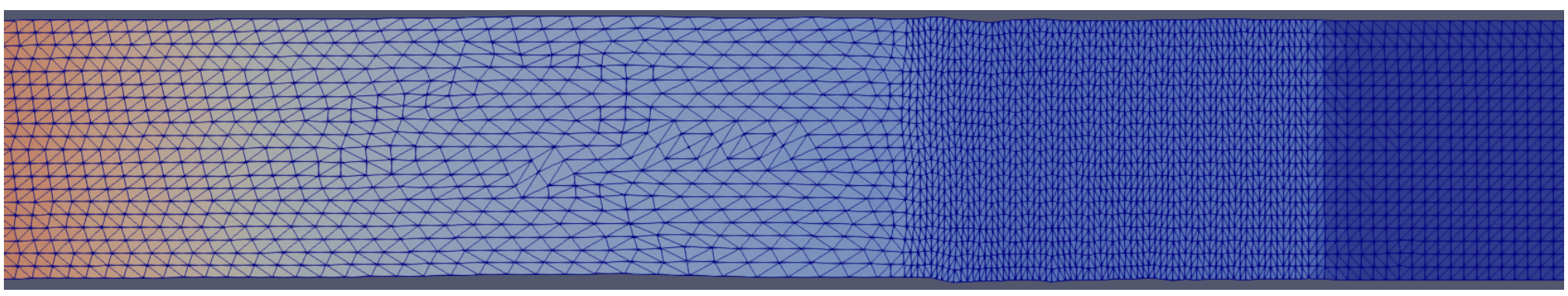

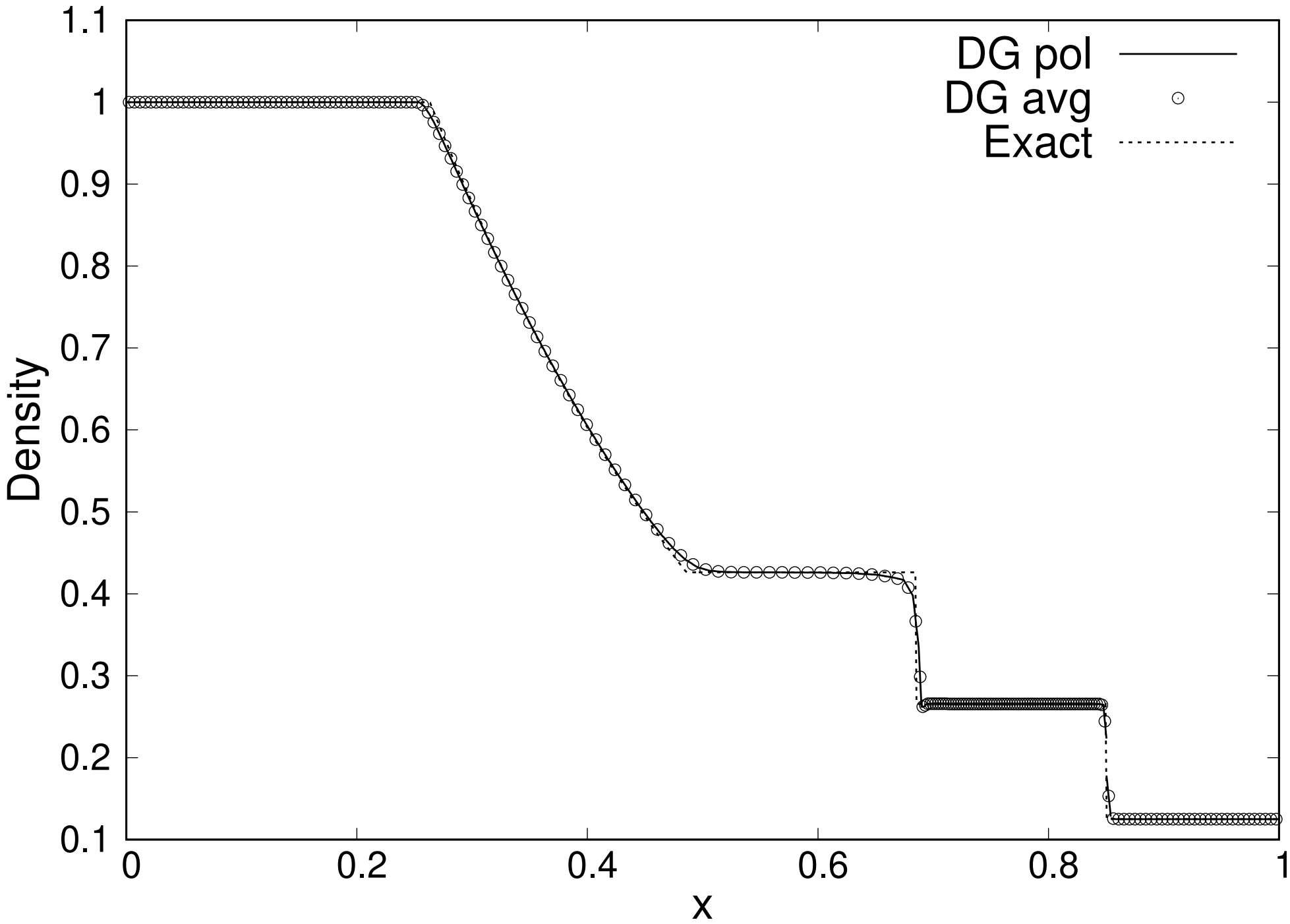

## 6.3 Sedov Taylor Blast Wave Solution

Another famous hydrodynamical test problem is the strong explosion, which creates a circular shock in two dimensions, known as Sedov-Taylor blast wave. This phenomenon occurs if a large amount of energy is released in a small region and on a short time scale.

The test problem can be modelled by injecting a large amount of energy into a point of fluid at rest. The initial density is $\rho = 1$ everywhere, velocity is zero and the total energy is equal to $10^{-12}$ everywhere but in the center of the domain where it is a constant $E_0$. The concentration of energy in the center generates a shock wave that propogates radially with time. A rapidly expanding rarefaction region is developed in the wake of the shock wave with a smooth solution in the region.

Because of its nearly-vacuum conditions, it is used to test whether the scheme is positivity-preserving. The problem is also interesting from the perspective of moving-mesh methods. The simultaneous existence of rapidly expanding and rapidly contracting regions present an opportunity to test the self-adapting properties of the mesh in the region of shock. We are also interested in looking at how well the method maintains the radial symmetry of both the solution as well as the mesh.

In the light of above considerations, we have two test cases here, one with a uniform transfinite mesh and one with uniform unstructured mesh. The uniform mesh has a domain of $[-1, 1] \times [-1, 1]$ with $xxxxx$ cells. The final time of the simulation is $T$.

As we can see, the solution is more accurate in the case of moving meshes. We can also see that the mesh tracks the shock accurately and leaves a coarse mesh in the rarefaction region. The radial symmetry of the solution is also maintained. We further observe that the unstructured mesh maintains the symmetry of the solution better than the transfinite mesh.

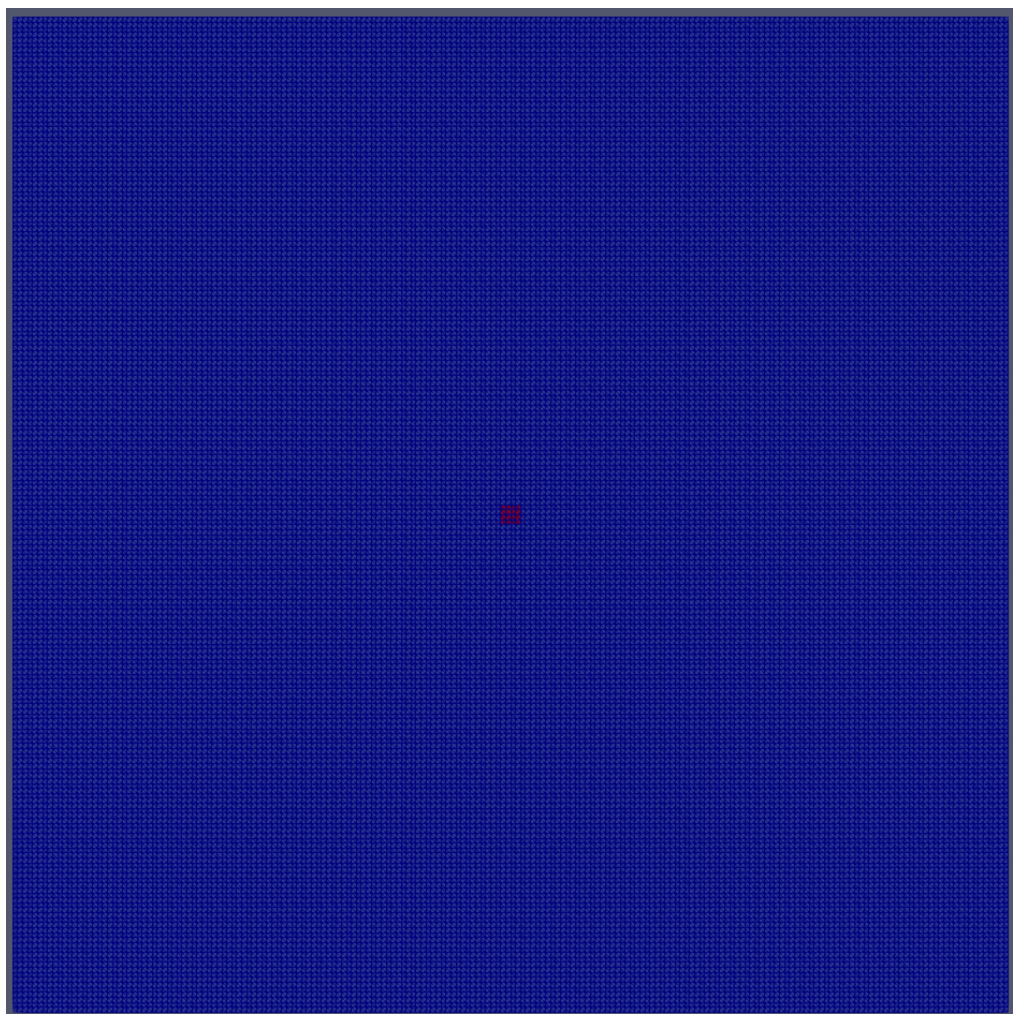

Figure 6.3.1: Transfinite Mesh: Pressure at Initial Time

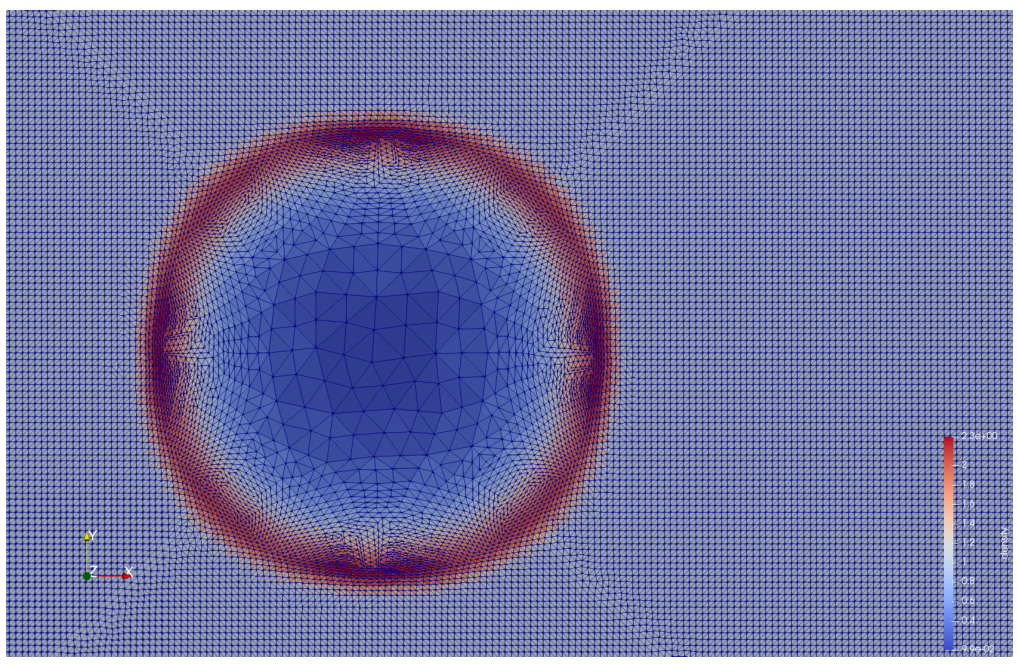

Figure 6.3.2: Transfinite Mesh: Solution at Final Time

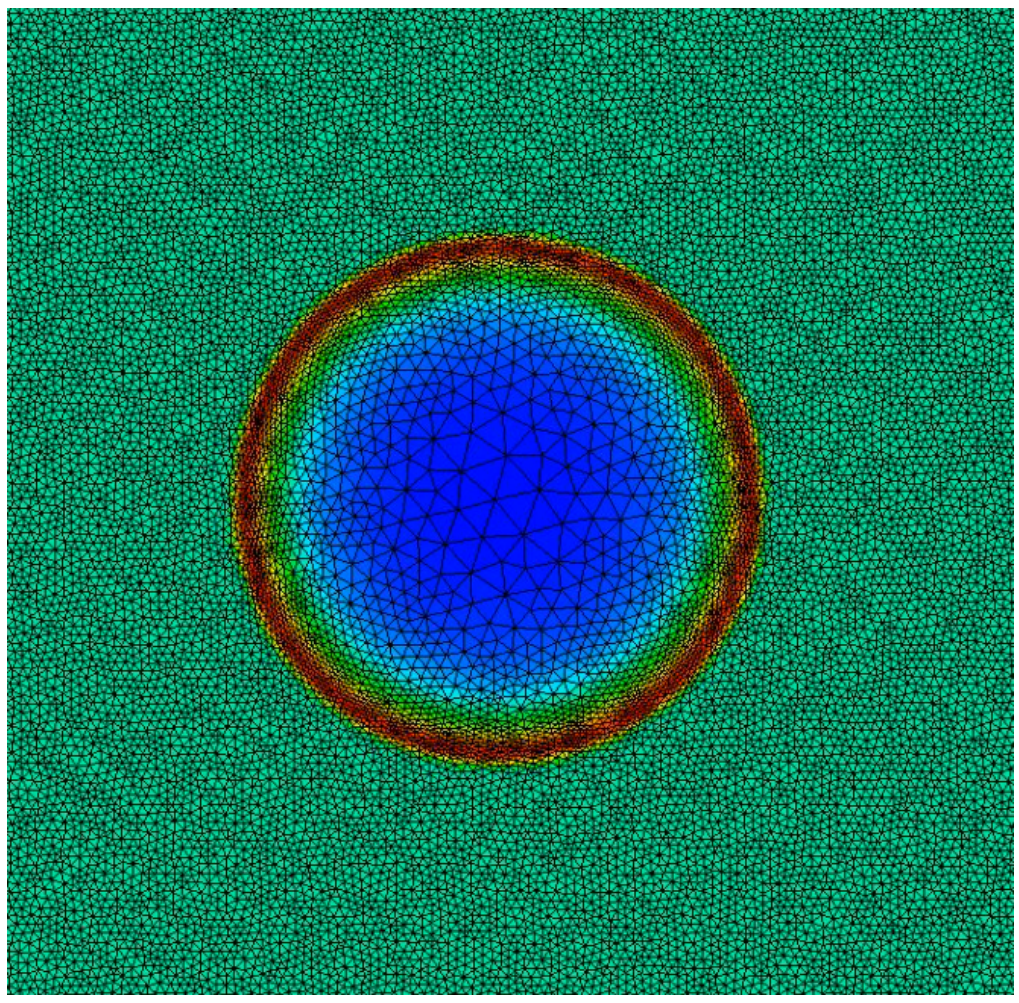

Figure 6.3.3: Unstructured Mesh: Solution at Final Time

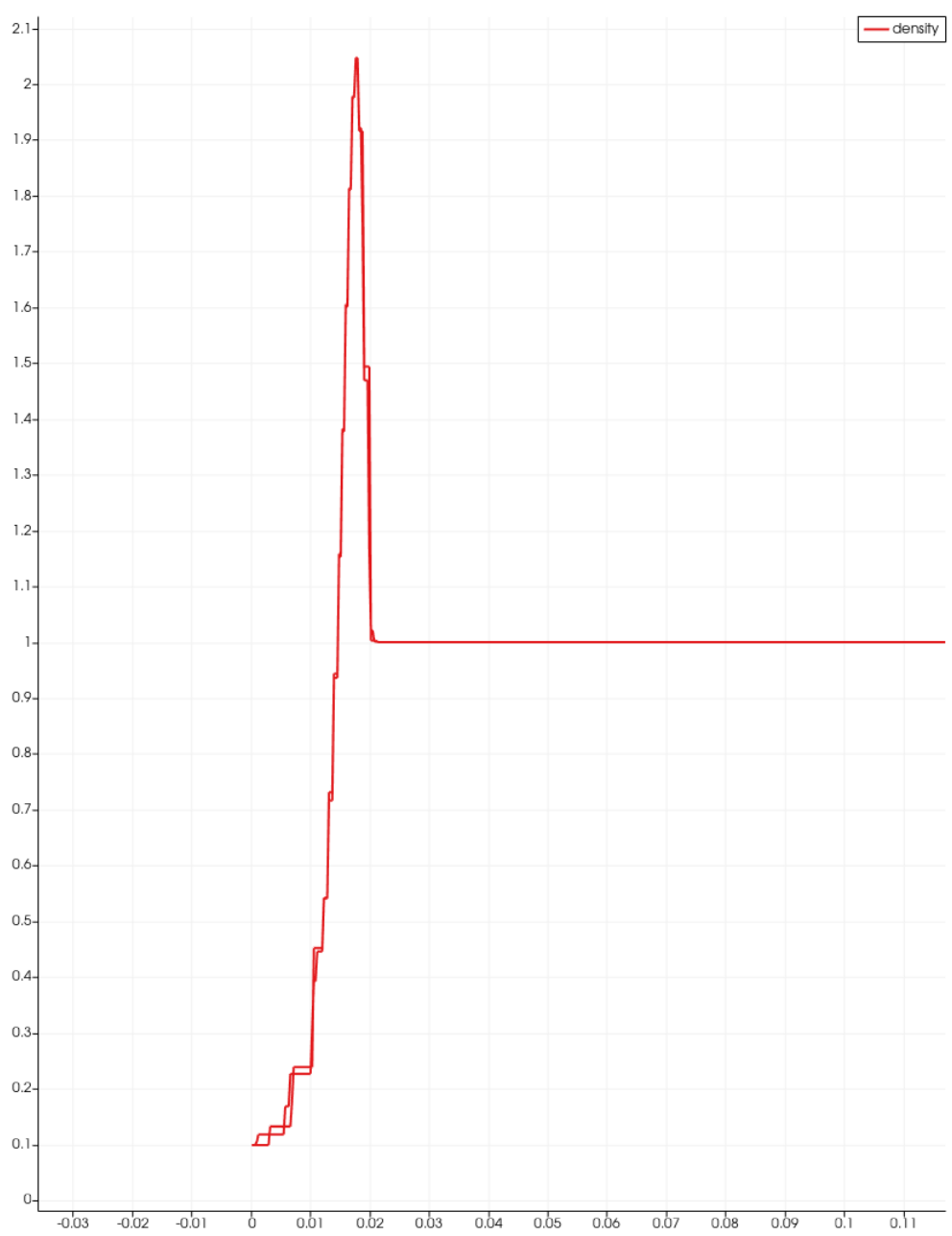

Figure 6.3.4: Radial Disparity in Solutions

## 6.4 Titarev-Toro Test Case

Titarev-Toro problem is an extension of the Shu-Osher problem [**Titarev2014**] to test a severely oscillatory wave interacting with a shock wave. It aims to test the ability of higher-order

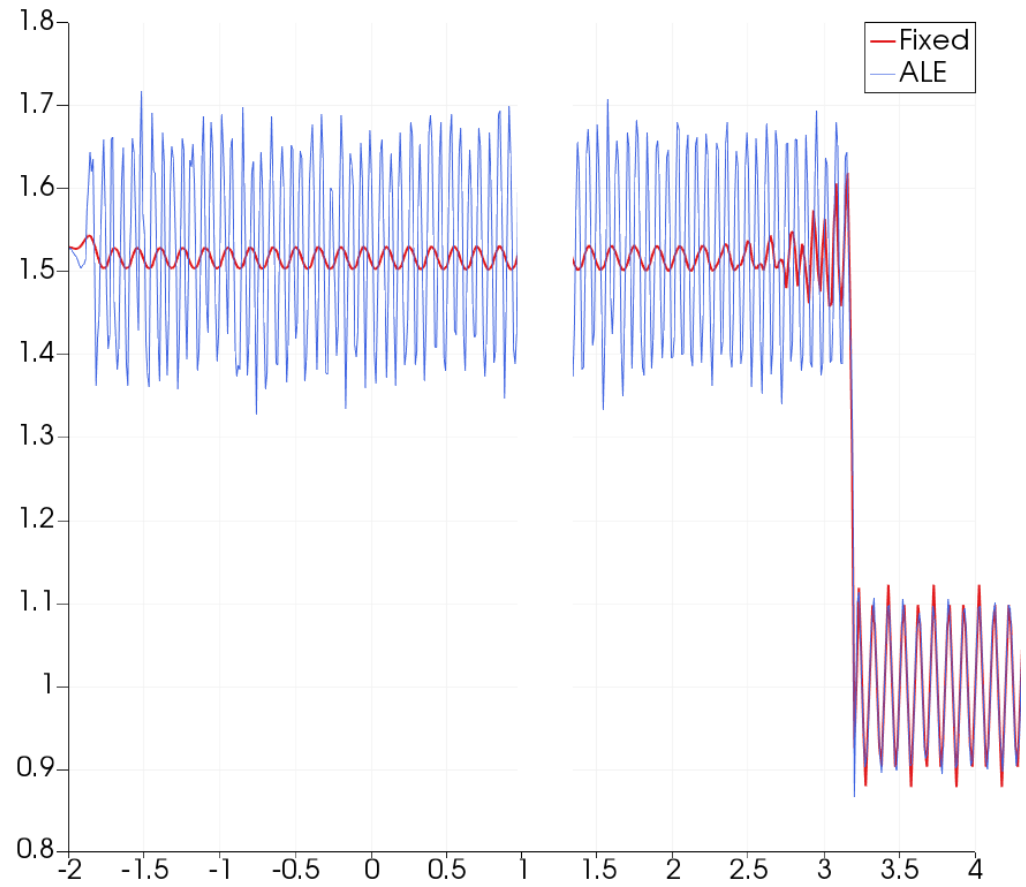

Figure 6.4.1: Titarev Toro Simulation Cross Section

methods to capture the extremely high frequency waves. The initial condition is given by

$$(\rho, v, p) = \begin{cases} (1.515695, 0.523346, 1.805), & -5 < x \leq -4.5 \\ (1 + 0.1\sin(20\pi x), 0, 1), & -4.5 < x \leq 5 \end{cases} \quad (6.4.1)$$

Just like in the one dimensional case, in the two dimensions, the fixed mesh is not able to resolve the high frequency oscillations due to dissipation in the fluxes and the TVD limiter, but the ALE scheme gives an excellent resolution of these high frequency oscillations. Note that the ALE scheme also uses the same TVD limiter but it is still able to resolve the solution to a very degree of accuracy. This result again demonstrates the superior accuracy that can be achieved by using a nearly Lagrangian ALE scheme in problems involving interaction of shocks and smooth flow structures.